\documentclass[a4paper,10pt]{article}

\usepackage{amsmath,amssymb,amsthm,mathtools,mathrsfs}
\renewcommand{\Re}{\operatorname{Re}}
\renewcommand{\Im}{\operatorname{Im}}
\usepackage{booktabs,array,enumitem}
\usepackage{indentfirst}
\usepackage[a4paper,textwidth=164mm,textheight=247mm,centering]{geometry}
\usepackage{verbatim}
\usepackage{float}
\usepackage{tikz}
\usetikzlibrary{arrows.meta,patterns,calc,positioning}
\usepackage{microtype}
\usepackage{xcolor}
\usepackage{caption}
\usepackage[hidelinks]{hyperref}
\allowdisplaybreaks

\newcommand{\R}{\mathbb R}
\newcommand{\C}{\mathbb C}
\newcommand{\ri}{\mathrm i}
\newcommand{\re}{\mathrm e}
\newcommand{\rd}{\,\mathrm d}
\newcommand{\one}{\mathbf 1}
\newcommand{\supp}{\operatorname{supp}}
\newcommand{\Res}{\operatorname{Res}}
\newcommand{\pv}{\operatorname{p.v.}}
\newcommand{\Lcal}{\mathcal L}
\newcommand{\Lap}{\mathscr{L}}
\newcommand{\Acal}{\mathcal A}
\newcommand{\Bcal}{\mathcal B}
\newcommand{\Om}{\Omega}
\newcommand{\ang}[1]{\langle #1\rangle}
\newtheorem{theorem}{Theorem}[section]
\newtheorem{proposition}[theorem]{Proposition}
\newtheorem{lemma}[theorem]{Lemma}
\newtheorem{corollary}[theorem]{Corollary}
\theoremstyle{definition}

\theoremstyle{remark}
\newtheorem{remark}[theorem]{Remark}
\numberwithin{equation}{section}

\title{Pointwise Wave--Damping for Linearized Vlasov--Poisson:\\
Compactly Supported Radial Equilibria\\
and Relativistic Applications}
\author{Haina Li, Quoc-Hung Nguyen, and Zhilong Xue}
\hypersetup{
pdfauthor={Haina Li, Quoc-Hung Nguyen, and Zhilong Xue},
pdftitle={Pointwise Wave-Damping for Linearized Vlasov-Poisson: Compactly Supported Radial Equilibria and Relativistic Applications},
pdfsubject={Weighted pointwise estimates for the linearized Vlasov-Poisson system},
pdfkeywords={Landau damping; Vlasov-Poisson system; dispersive estimates; Green function; survival threshold}
}
\date{}

\begin{document}
\maketitle

\begin{abstract}
We establish weighted pointwise estimates for the response Green function and the complete electric field of the linearized Vlasov--Poisson system in 
\(\R_x^3\times\R_v^3\) around compactly supported radial equilibria. Both decompose into a dispersive wave component with
\(t^{-3/2}\) decay in $L_x^\infty$ and a faster-decaying Landau-damped component. The wave estimate includes the survival threshold, where the oscillatory branch meets the continuous spectrum, while the damped estimates resolve the spatial fronts associated with the velocity boundaries. The equilibrium class includes shells with independent inner and outer orders of vanishing, down to first order, and profiles positive throughout a velocity ball. The proof combines a joint treatment of pole and continuum contributions with a low-frequency normal form for the electric field. The electric-field estimates hold for real, neutral perturbations with finite weighted regularity and imply uniform scattering along free characteristics. An exact momentum--velocity correspondence for suitable filled-ball profiles yields analogous estimates and scattering for a class of nonnegative relativistic equilibria with polynomial momentum tails.
\end{abstract}

\par\smallskip
{\small
\par
}

\section{Introduction and main results}
\label{sec:introduction}

The Vlasov--Poisson system provides a kinetic description of
collisionless plasmas in which charged particles interact through a
self-consistent electrostatic field. Understanding the relaxation of
perturbations of a spatially homogeneous equilibrium is a basic problem
in this theory. Landau damping describes the decay of macroscopic
electric perturbations in the absence of collisional dissipation,
and has motivated extensive work on both the linearized and nonlinear
dynamics \cite{Lan46,MV11}.
In the whole space, this relaxation also involves spatial dispersion.
Long-wavelength oscillations may persist at the level of individual
Fourier modes while their superposition decays in physical space.
For the unscreened Coulomb interaction, the coexistence of these
oscillations with a more rapidly damped response leads to a
wave--damping description of the evolution \cite{BMM22,Ngu26}.

Compact velocity support introduces additional features into the spectral and spatial behavior of the linearized response. The velocity boundaries determine endpoints of the cut contribution in the response representation, and their orders of vanishing influence
the decay near the corresponding spatial fronts. A further difficulty occurs at the survival threshold, where an oscillatory pole meets the continuum and the two contributions must be analyzed together. Velocity shells also have an inner boundary, whose order of vanishing may differ from that at the outer boundary. We consider both such shells and profiles positive throughout a velocity ball.

\subsection{The equation, equilibria, and initial data}

The classical and relativistic Vlasov--Poisson systems describe
collisionless particles coupled to a self-consistent electrostatic field:
\[
 \partial_t f+b(v)\cdot\nabla_xf+E\cdot\nabla_vf=0,
 \qquad E=-\nabla_x\phi,
 \qquad -\Delta_x\phi=\int_{\R^3}f(t,x,v)\,\rd v-\rho_{\rm bg}.
\]
Here $\rho_{\rm bg}$ is a spatially uniform neutralizing background
density. We use $\ang v=(1+|v|^2)^{1/2}$ and the same convention
for tuples of variables. The transport velocity is $b(v)=v$ in the
classical model, where $v$ denotes velocity, and
$b(v)=v/\ang v$ in the relativistic model, where $v$ denotes momentum
and the speed of light is normalized to one.

Let $f^0=f^0(v)$ be a spatially homogeneous equilibrium, with
$\rho_{\rm bg}=\int_{\R^3}f^0(v)\,\rd v$.
Writing $f=f^0+h$ and neglecting the quadratic term
$E\cdot\nabla_v h$, we obtain the general linearized system
\[
\begin{cases}
\partial_t h+b(v)\cdot\nabla_xh
    +E(t,x)\cdot\nabla_v f^0(v)=0,\\
E=-\nabla_x\phi,\qquad
-\Delta_x\phi=\rho,\qquad
\displaystyle\rho(t,x)=\int_{\R^3}h(t,x,v)\,\rd v,\\
h(0,x,v)=h_{\rm in}(x,v).
\end{cases}
\]

We mainly consider the classical case $b(v)=v$, with a radial
equilibrium $f^0(v)$ specified below. The corresponding linearized
problem reads
\begin{equation}
\label{eq:linear-VP}
\begin{cases}
\partial_t h+v\cdot\nabla_xh+E(t,x)\cdot\nabla_v f^0(v)=0,
 &(t,x,v)\in\R_+\times\R^3\times\R^3,\\
E=-\nabla_x\phi,\qquad -\Delta_x\phi=\rho,
 &\displaystyle\rho(t,x)=\int_{\R^3}h(t,x,v)\,\rd v,\\
h(0,x,v)=h_{\rm in}(x,v).
\end{cases}
\end{equation}
The equilibrium is radial,
\begin{equation}
\label{eq:radial-equilibrium}
f^0(v)=\mu\!\left(\frac{|v|^2}{2}\right),
\qquad
\int_{\R^3}f^0(v)\,\rd v=1.
\end{equation}

We allow either a velocity shell or a filled velocity ball.  The shell
hypotheses are
\begin{enumerate}[label=$(\mathbf{H}\theenumi)$,ref=$\mathbf{H}\theenumi$]
	\item\label{H1} There exists \(\kappa\in(0,1)\) such that
	\(r\mapsto\mu(r^2/2)\) is smooth and strictly positive on
	\((\kappa,1)\), and vanishes on
	\([0,\kappa]\cup[1,\infty)\).
	\item\label{H2} There exist
	\(c_1,c_2\in(0,\frac12\min\{\kappa,1-\kappa\})\) and
	\({\gamma_1,\gamma_\kappa\geq1}\) such that
	\begin{itemize}[label=$\bullet$]
		\item \(s\mapsto
		\mu(s^2/2)/(1-s)^{{\gamma_1}}\), initially defined on
		\([1-c_1,1)\), extends to a strictly positive
		\(C^\infty\) function on \([1-c_1,1+c_1]\);
		\item \(s\mapsto
		\mu(s^2/2)/(s-\kappa)^{{\gamma_\kappa}}\), initially defined on
		\((\kappa,\kappa+c_2]\), extends to a strictly positive
		\(C^\infty\) function on
		\([\kappa-c_2,\kappa+c_2]\).
	\end{itemize}
\end{enumerate}
Alternatively, we consider the filled-ball case ($\kappa=0$):
\begin{enumerate}[label=$(\mathbf{H0})$,ref=$\mathbf{H0}$]
\item\label{H0}
The function $\mu$ is strictly positive on $[0,1/2)$ and vanishes on
$[1/2,\infty)$.  It is smooth on $[0,1/2)$ and has a smooth extension
through zero.  There exist $c_1\in(0,1/2)$ and $\gamma_1\ge1$ such that
\[
 s\longmapsto\frac{\mu(s^2/2)}{(1-s)^{\gamma_1}},
 \qquad s\in[1-c_1,1),
\]
extends to a strictly positive smooth function on $[1-c_1,1+c_1]$.
\end{enumerate}

Two typical families are the shell profiles
\[
 f^0(v)=c_{\kappa,\gamma_\kappa,\gamma_1}
 (1-|v|^2)_+^{\gamma_1}(|v|^2-\kappa^2)_+^{\gamma_\kappa},
 \qquad 0<\kappa<1,\quad \gamma_\kappa,\gamma_1\ge1,
\]
and the filled-ball profiles
\[
 f^0(v)=c_\gamma(1-|v|^2)_+^\gamma,
 \qquad c_\gamma=\frac{\Gamma(\gamma+5/2)}
 {\pi^{3/2}\Gamma(\gamma+1)},\qquad \gamma\ge1.
\]
Here $\Gamma$ is Euler's gamma function, $z_+=\max\{z,0\}$, and
the positive constant in the shell family is chosen to give unit mass.  The first family satisfies
\eqref{H1}--\eqref{H2}; the second satisfies \eqref{H0} with
$\kappa=0$ and $\gamma_1=\gamma$.

For integers $M\ge11$, set
\begin{equation}
\label{eq:data-norm}
[h_{\rm in}]_{M}
:=
\sum_{|\alpha|+|\beta|\leq M}
\sup_{x,v\in\R^3}
\ang{x,v}^{M+6}
\left|\partial_x^\alpha\partial_v^\beta h_{\rm in}(x,v)\right|.
\end{equation}
For the electric-field and scattering conclusions we assume that
$h_{\rm in}$ is real-valued and neutral:
\begin{equation}
\label{eq:neutrality}
\iint_{\R^3\times\R^3}h_{\rm in}(x,v)\,\rd x\,\rd v=0.
\end{equation}

\subsection{The Green function and the wave--damping decomposition}

We use the spatial Fourier convention
\begin{equation}
\label{eq:fourier-transform-intro}
\widehat f(\xi)=\int_{\R^3}\re^{-\ri x\cdot\xi}f(x)\,\rd x,
\qquad
f(x)=\frac{1}{(2\pi)^3}\int_{\R^3}
\re^{\ri x\cdot\xi}\widehat f(\xi)\,\rd\xi,
\end{equation}
and the same convention in the velocity variable.  For a spatial
multiplier $a(D_x)$, we use $\widehat{a(D_x)f}=a\widehat f$ and denote
the inverse Fourier transform of $a$ by $\check a$.
Let
\begin{equation}
\label{eq:free-density-field-intro}
\rho^{\rm free}(t,x)
:=\int_{\R^3}h_{\rm in}(x-tv,v)\,\rd v,
\qquad
E^{\rm free}:=-\nabla_x(-\Delta_x)^{-1}\rho^{\rm free}.
\end{equation}
Here $h^{\rm free}(t,x,v)=h_{\rm in}(x-tv,v)$ solves the free transport equation 
$\partial_t h^{\rm free}+v\cdot\nabla_xh^{\rm free}=0$, and
$\rho^{\rm free}=\int_{\R^3}h^{\rm free}\,\rd v$.
Thus, if $\widehat h_{\rm in}(\xi,\eta)$ denotes the Fourier transform
of the initial datum in both variables, then
$\widehat{\rho^{\rm free}}(t,\xi)=\widehat h_{\rm in}(\xi,t\xi)$.

The density equation associated with \eqref{eq:linear-VP} is the
Volterra equation
\begin{equation}
\label{eq:volterra-intro}
\widehat\rho(t,\xi)
=\widehat{\rho^{\rm free}}(t,\xi)
+\int_0^t\widehat K(t-s,\xi)\widehat\rho(s,\xi)\,\rd s,
\qquad
\widehat K(t,\xi)=-t\widehat{f^0}(t\xi).
\end{equation}
The kernel $\widehat K$ describes the contribution of the
self-consistent field to the density equation.  We call the resolvent
kernel of this Volterra equation, written in physical space, the
\emph{response Green function} $G$.  It accounts for the induced response
beyond the free-transport source.  We take \(G\) to be
\emph{causal}, meaning that it vanishes for \(t<0\); at time \(t\),
its convolutions use only the history \(0\le s\le t\).
For each spatial frequency, \(\widehat G\) is defined by
\begin{equation}
\label{eq:green-volterra-definition}
\widehat G(t,\xi)
=\widehat K(t,\xi)+\int_0^t\widehat K(t-s,\xi)\widehat G(s,\xi)\,\rd s.
\end{equation}
We use $*_x$ for spatial convolution and $*_t$ for causal time
convolution.  The symbol $*_{t,x}$ denotes their combination; for functions,
whenever the integral exists, it is given by
\[
(F*_{t,x}H)(t,x)
:=\int_0^t\int_{\R^3}F(t-s,x-y)H(s,y)\,\rd y\,\rd s.
\]
Section~\ref{subsec:volterra-construction-v3} constructs this resolvent \eqref{eq:green-volterra-definition}
by a convergent Volterra series on every finite time interval.
The same section derives the density representation below by
substitution in \eqref{eq:volterra-intro} and finite-time uniqueness.
Applying the Poisson multiplier, which commutes with spatial
convolution, gives the field representation:
\begin{equation}
\label{eq:green-field-representation-intro}
\rho=\rho^{\rm free}+G*_{t,x}\rho^{\rm free},
\qquad
E=E^{\rm free}+G*_{t,x}E^{\rm free}.
\end{equation}
The full density solution kernel is
\(\delta_{t=0}\delta_{x=0}+G\); the identity term supplies the free
density, while the estimates below concern the response \(G\).

The wave--damping decomposition is guided by the spectral description
in Sections~\ref{sec:dispersion}--\ref{sec:spectral-decomposition}.
For a fixed spatial frequency $k=|\xi|$, the Laplace transform of
$\widehat G(\cdot,\xi)$ agrees on $\Re z>0$ with $R(z,k)$ from
\eqref{eq:rho-resolvent-section2}; see
Corollary~\ref{cor:fixed-k-boundedness-section3}.
For $0<k<k_*$, its continuation has two simple imaginary-axis poles
$\pm\ri\tau_{\rm p}(k)$, with temporal frequency $\tau_{\rm p}(k)$
and survival threshold $k_*>0$.
Proposition~\ref{prop:pole-branch-section2} gives
$\tau_{\rm p}(k)>k$ and shows that, as $k$ increases to $k_*$,
the poles meet the outer continuum endpoints $\tau=\pm k$.
We write $\varrho_\pm(k)$ for the residues of $R(z,k)$ at these poles.
Their inverse contributions are
$\varrho_\pm(k)\re^{\pm\ri t\tau_{\rm p}(k)}$.

Near this contact the pole and the continuum concentration must be
retained together.  Away from it we separate the smooth pole wave
from the remaining continuum.  To implement this distinction, choose
a smooth partition depending only on the spatial frequency \(k\),
\begin{equation}
\label{eq:k-partition-intro}
\chi_\ell(k)+\vartheta(k)+\chi_h(k)=1,
\end{equation}
where
$\supp\chi_\ell\subset[0,k_*)$ and
$\supp\chi_h\subset(k_*,\infty)$, both at a positive distance from
$k_*$, and $\vartheta=1$ on a fixed neighborhood of $k_*$.  Moreover,
$\chi_\ell=1$ at sufficiently low frequencies and $\chi_h=1$
sufficiently far above the threshold.  The separated supports are
specified in the construction leading to \eqref{eq:three-blocks-section3}.  We then set
\begin{equation}
\label{eq:green-wave-LD-intro}
\widehat G_{\rm wave}(t,\xi)
:=\chi_\ell(k)\sum_{\sigma=\pm}\varrho_\sigma(k)
\re^{\ri\sigma t\tau_{\rm p}(k)}
+\vartheta(k)\widehat G(t,\xi),
\qquad
G_{\rm LD}:=G-G_{\rm wave},
\qquad k=|\xi|.
\end{equation}
The pole expression is understood by its continuous extension at
$k=0$.

For the electric field, two integrations in time of the smooth pole
response expose a current-time boundary term that combines with the
free density at low frequencies.  To state this normal form,
introduce the first two temporal primitives of the smooth pole term
in \eqref{eq:green-wave-LD-intro}:
for \(j=1,2\), define
\begin{equation}
\label{eq:pole-primitives-intro}
\widehat K_j(t,\xi)
:=\chi_\ell(k)\sum_{\sigma=\pm}
\frac{\varrho_\sigma(k)}
{(\ri\sigma\tau_{\rm p}(k))^j}
\re^{\ri\sigma t\tau_{\rm p}(k)},
\qquad k=|\xi|.
\end{equation}
We set
\begin{align}
\rho_{\rm wave}
&:=K_1(t)*_x\rho^{\rm free}(0)+K_2(t)*_x\partial_t\rho^{\rm free}(0)
+K_2*_{t,x}\partial_t^2\rho^{\rm free}+\vartheta(D_x)\rho,
\label{eq:rho-wave-intro}\\
\rho_{\rm LD}
&:=\bigl[1-\vartheta(D_x)-\widehat K_1(0,D_x)\bigr]\rho^{\rm free}
+G_{\rm LD}*_{t,x}\rho^{\rm free},
\label{eq:rho-LD-intro}\\
E_{\rm wave}
&:=-\nabla_x(-\Delta_x)^{-1}\rho_{\rm wave},
\qquad
E_{\rm LD}:=-\nabla_x(-\Delta_x)^{-1}\rho_{\rm LD}.
\label{eq:electric-components-intro}
\end{align}
The normal form \eqref{eq:pole-normal-form-section5}, proved in
Proposition~\ref{prop:pole-normal-form-section5}, gives
\(\rho=\rho_{\rm wave}+\rho_{\rm LD}\), and hence
\begin{equation}
\label{eq:electric-sum-intro}
E=E_{\rm wave}+E_{\rm LD}.
\end{equation}
The current-time boundary term cancels the order-zero free-density
contribution at low frequencies:
\[
\chi_\ell(k)-\widehat K_1(0,k)=k^2b_0(k),
\]
with \(b_0\) smooth, by the sum rule
\eqref{eq:pole-sum-rule-section5}.
This $k^2$ factor removes the Coulomb singularity in the local
term.  The pole-subtracted response supplies a second such factor,
$R(z,k)-P(z,k)=k^2\widetilde{\mathcal R}(z/k,k)$ at low frequency;
see \eqref{eq:R-minus-P-slit-factor-section3}.  All cutoffs in
\eqref{eq:rho-wave-intro}--\eqref{eq:rho-LD-intro} depend only on \(k\),
and the entire threshold block belongs to \(\rho_{\rm wave}\).

\subsection{Classical estimates and scattering}

Under \eqref{H1}--\eqref{H2}, set
$\gamma_{\min}=\min(\gamma_\kappa,\gamma_1)$.
Under \eqref{H0}, set $\kappa=0$ and $\gamma_{\min}=\gamma_1$.
All endpoint sums below run over the positive elements of $\{\kappa,1\}$,
so that only the outer front remains in the filled-ball case.
\begin{theorem}
\label{thm:main-green}
Assume \eqref{eq:radial-equilibrium} and either the shell hypotheses
\eqref{H1}--\eqref{H2} or the filled-ball hypothesis \eqref{H0}.
Let $G$ be the response defined by \eqref{eq:green-volterra-definition},
$G_{\rm wave}$ the pole--threshold component in
\eqref{eq:green-wave-LD-intro}, and $G_{\rm LD}=G-G_{\rm wave}$.
The full response has finite propagation speed:
\begin{equation}
\label{eq:green-finite-propagation-main}
 \supp G(t,\cdot)\subset\{x\in\R^3:|x|\le t\},\qquad t>0.
\end{equation}
For every integer $N\ge0$,
\begin{equation}
\label{eq:main-green-wave-weighted}
 |G_{\rm wave}(t,x)|
 \le C_N\ang t^{-3/2}\frac{\ang t^N}{\ang{t,x}^N},
 \qquad t>0.
\end{equation}
The Landau-damped component satisfies
\begin{equation}
\label{eq:main-green-LD-envelope}
 |G_{\rm LD}(t,x)|
 \le \frac{C\bigl(1+t^{-2}\one_{\{\kappa t\le|x|\le t\}}\bigr)}
 {\ang{t,x}^{6}}
 \left[1+\sum_{\substack{a\in\{\kappa,1\}\\a>0}}
 \left(\frac{\ang{t,x}}{\ang{|x|-at}}\right)^{(4-\gamma_a)_+}\right],
 \qquad t>0.
\end{equation}
\end{theorem}

The corresponding bounds for the complete electric field use the data norm \eqref{eq:data-norm}.

\begin{theorem}
\label{thm:main-electric}
Assume \eqref{eq:radial-equilibrium} and either the shell hypotheses
\eqref{H1}--\eqref{H2} or the filled-ball hypothesis \eqref{H0}.
Let $M_*=11$.  Suppose that $h_{\rm in}$ is real-valued, satisfies
the neutrality condition \eqref{eq:neutrality}, and obeys
\[
 [h_{\rm in}]_{M_*}
 =\sum_{|\alpha|+|\beta|\le M_*}\sup_{x,v\in\R^3}
 \ang{x,v}^{M_*+6}
 |\partial_x^\alpha\partial_v^\beta h_{\rm in}(x,v)|<\infty.
\]
Then the complete electric field of \eqref{eq:linear-VP} decomposes as $E=E_{\rm wave}+E_{\rm LD}$,
with the components defined by
\eqref{eq:rho-wave-intro}--\eqref{eq:electric-components-intro}.
For every $t\ge0$,
\begin{equation}
\label{eq:main-electric-wave}
 |E_{\rm wave}(t,x)|
 \le C\frac{\ang t^{3/2}}{\ang{t,x}^3}[h_{\rm in}]_{M_*}, 
\end{equation}
while the Landau-damped component satisfies the front estimate\footnote{The leading spherical term leaves a radial factor
$k^{-1}$ when $\gamma_a=2$. Integrating it over frequencies
between scales $\ang{t,x}^{-1}$ and $\ang{|x|-at}^{-1}$
produces the logarithm displayed in
\eqref{eq:main-electric-front-v3}.}
\begin{equation}\label{eq:main-electric-front-v3}
\begin{aligned}
 |E_{\rm LD}(t,x)|
 \le \frac{C[h_{\rm in}]_{M_*}}{\ang{t,x}^{4}}
 \biggl[1+\sum_{\substack{a\in\{\kappa,1\}\\a>0}}
 \biggl\{\left(\frac{\ang{t,x}}{\ang{|x|-at}}\right)^{(2-\gamma_a)_+}+\one_{\{\gamma_a=2\}}
 \log\!\left(2+\frac{\ang{t,x}}{\ang{|x|-at}}\right)\biggr\}\biggr],
\end{aligned}
\end{equation}
and the spatially weighted decay estimate
\begin{equation}\label{eq:main-electric-gamma-v3}
 |E_{\rm LD}(t,x)|
 \le \frac{C[h_{\rm in}]_{M_*}}{\ang{t,x}^{4}}
 \left[\ang t^{(2-\gamma_{\min})_+}
 +\one_{\{\gamma_{\min}=2\}}\log(2+\ang t)\right],
 \qquad t\ge0.
\end{equation}
Moreover, for $1<p\le\infty$ and $t\ge0$,
\begin{equation}\label{eq:main-electric-global}
 \|E(t)\|_{L_x^p}+\|E_{\rm wave}(t)\|_{L_x^p}
 \le C_p\ang t^{-3/2+3/p}[h_{\rm in}]_{M_*},
\end{equation}
and at the weak $L^1$ endpoint,
\begin{equation}\label{eq:E-weak-L1-section6}
 \|E(t)\|_{L_x^{1,\infty}}+\|E_{\rm wave}(t)\|_{L_x^{1,\infty}}
 \le C\ang t^{3/2}[h_{\rm in}]_{M_*},\qquad t\ge0,
\end{equation}
where $\|f\|_{L^{1,\infty}}:=\sup_{\lambda>0}\lambda
|\{x:|f(x)|>\lambda\}|$.
\end{theorem}

The cubic spatial weight in \eqref{eq:main-electric-wave} cannot
in general be improved under neutrality alone; see
Remark~\ref{rem:wave-dipole-v4}. The case $p=\infty$ of \eqref{eq:main-electric-global} is
time-integrable and yields uniform scattering along free characteristics.

\begin{corollary}
\label{cor:main-scattering}
Under the assumptions of Theorem~\ref{thm:main-electric}, the final state
\begin{equation}
\label{eq:h-infinity}
h_\infty(x,v)
=h_{\rm in}(x,v)
-\int_0^\infty
E(s,x+sv)\cdot\nabla_v f^0(v)\,\rd s
\end{equation}
is well defined and
\begin{equation}
\label{eq:main-scattering}
\operatorname*{ess\,sup}_{x,v}
|h(t,x+tv,v)-h_\infty(x,v)|
\leq C\ang t^{-1/2}[h_{\rm in}]_{M_{\ast}},
\qquad t\geq0.
\end{equation}
\end{corollary}

The front factors in \eqref{eq:main-green-LD-envelope} and
\eqref{eq:main-electric-front-v3} distinguish neighborhoods of the rays
\(|x|=at\) from the interior and the far field. The velocity boundaries
of the equilibrium determine these rays; the initial perturbation need
not have the same velocity support. Each front bound uses its own
order \(\gamma_a\), while the wave component includes the entire
pole--continuum contact region.

Section~\ref{sec:consequences} proves the scattering statement and
the relativistic consequences below.

\subsection{Relativistic consequences}
\label{subsec:relativistic-results-intro}

We relate the response kernels of the classical and relativistic
linearized Vlasov--Poisson systems by an exact change from momentum
to velocity.  The classical estimates then apply to the matched
backgrounds and transformed data below.  The proofs are given in
Section~\ref{subsec:relativistic-correspondence}.
For a classical equilibrium $f^0(u)=\mu(|u|^2/2)$ satisfying \eqref{H0} and
\eqref{eq:radial-equilibrium}, define the radial momentum profile
\begin{equation}\label{eq:relativistic-general-profile}
 \Phi(v):=-\int_{|v|/\ang v}^{1}(1-s^2)s\mu'(s^2/2)\,\rd s,
 \qquad v\in\R^3.
\end{equation}
We assume that $\Phi\ge0$; this holds, for example, if
$s\mapsto\mu(s^2/2)$ is nonincreasing.  The profile is smooth at the
origin and has a polynomial tail of order $\ang v^{-2\gamma_1-2}$,
and hence finite mass, as verified in Section~\ref{subsec:relativistic-correspondence}.
The restriction to \eqref{H0} is intrinsic to this construction:
for a genuine shell, the same matching formula produces a negative
profile at small momenta, as shown in that section.
Taking $b(v)=v/\ang v$ and $f^0(v)=\Phi(v)$ in the general
linearization above, and writing $h_{\rm rel}$ and $E_{\rm rel}$
for the perturbation and its electric field, we obtain the
linearized relativistic Vlasov--Poisson system
\begin{equation}\label{eq:linear-relativistic-VP-intro}
\begin{cases}
\displaystyle
\partial_t h_{\rm rel}+\frac{v}{\ang v}\cdot\nabla_xh_{\rm rel}
   +E_{\rm rel}\cdot\nabla_v\Phi(v)=0,\\[1mm]
E_{\rm rel}=-\nabla_x\phi_{\rm rel},\qquad
-\Delta_x\phi_{\rm rel}=\rho_{\rm rel},\qquad
\displaystyle\rho_{\rm rel}(t,x)=\int_{\R^3}h_{\rm rel}(t,x,v)\,\rd v,\\[1mm]
h_{\rm rel}(0,x,v)=h_{{\rm rel,in}}(x,v).
\end{cases}
\end{equation}
We will show that its response Green function coincides exactly
with that of the classical system \eqref{eq:linear-VP} for the
corresponding radial equilibrium $f^0(u)=\mu(|u|^2/2)$. For $\xi\ne0$, define
\[
 \widehat K_{\rm rel}(t,\xi):=\frac{\ri\xi}{|\xi|^2}\cdot
 \int_{\R^3}\re^{-\ri t\xi\cdot v/\ang v}\nabla_v\Phi(v)\,\rd v.
\]
The causal response $G_{\rm rel}$ is defined in spatial Fourier
variables by
\begin{equation}\label{eq:relativistic-volterra-intro}
 \widehat G_{\rm rel}(t,\xi)=\widehat K_{\rm rel}(t,\xi)
 +\int_0^t\widehat K_{\rm rel}(t-s,\xi)
                 \widehat G_{\rm rel}(s,\xi)\,\rd s,
 \qquad t\ge0,
\end{equation}
and is zero for $t<0$.  The notation $G$ continues to denote the
classical response for the matched equilibrium $f^0$.

\begin{theorem}\label{thm:relativistic-green}
Assume \eqref{H0} and \eqref{eq:radial-equilibrium}, and suppose that
$\Phi$ in \eqref{eq:relativistic-general-profile} is nonnegative.
Then
\[
 \widehat K_{\rm rel}(t,\xi)=\widehat K(t,\xi)=-t\widehat f^0(t\xi),
 \qquad G_{\rm rel}=G,
\]
where $\widehat K$ is defined in \eqref{eq:volterra-intro}, and the
kernel identity at $\xi=0$ is understood by continuity.
All conclusions of Theorem~\ref{thm:main-green} therefore hold for
$G_{\rm rel}$ with $\kappa=0$ and only the outer front $|x|=t$.
The wave--damping decomposition is
\eqref{eq:green-wave-LD-intro} for the classical equilibrium $f^0$;
finite propagation holds for the full response.
\end{theorem}

For the initial-data correspondence, set
\[
 v(u):=\frac{u}{\sqrt{1-|u|^2}},\qquad
 J(u):=(1-|u|^2)^{-5/2},\qquad |u|<1.
\]
\begin{theorem}\label{thm:relativistic-electric}
Assume \eqref{H0} and \eqref{eq:radial-equilibrium}, and suppose that
$\Phi$ in \eqref{eq:relativistic-general-profile} is nonnegative.
Let the zero extension of
\begin{equation}\label{eq:relativistic-data-intro}
 g_{\rm in}(x,u):=J(u)h_{{\rm rel,in}}(x,v(u)),\qquad |u|<1,
\end{equation}
be real-valued and satisfy the neutrality condition
\eqref{eq:neutrality} and $[g_{\rm in}]_{11}<\infty$, with
\eqref{eq:data-norm} and \eqref{eq:neutrality} applied to $g_{\rm in}$
in the variables $(x,u)$.
The field $E_{\rm rel}$ of \eqref{eq:linear-relativistic-VP-intro}
equals the field of \eqref{eq:linear-VP} with velocity variable $u$,
equilibrium $f^0(u)$ and initial datum $g_{\rm in}$.
Its decomposition
$E_{\rm rel}=E_{{\rm rel},\rm wave}+E_{{\rm rel},\rm LD}$ is defined by
\eqref{eq:rho-wave-intro}--\eqref{eq:electric-components-intro}
using the free density generated by $g_{\rm in}$, and
\[
 |E_{{\rm rel},\rm wave}(t,x)|
 \le C\frac{\ang t^{3/2}}{\ang{t,x}^3}[g_{\rm in}]_{11}.
\]
The damped component satisfies
\eqref{eq:main-electric-front-v3} and \eqref{eq:main-electric-gamma-v3},
and the field and its wave component satisfy
\eqref{eq:main-electric-global} and \eqref{eq:E-weak-L1-section6},
with $\kappa=0$, $\gamma_{\min}=\gamma_1$, only $a=1$ in the front sum,
and $[g_{\rm in}]_{11}$ on the right-hand sides.
Moreover, there is a final state $h_{{\rm rel},\infty}$ such that
\[
 \operatorname*{ess\,sup}_{x,v}
 \left|h_{\rm rel}\!\left(t,x+\frac{tv}{\ang v},v\right)
       -h_{{\rm rel},\infty}(x,v)\right|
 \le C\ang t^{-1/2}[g_{\rm in}]_{11},\qquad t\ge0.
\]
\end{theorem}

\begin{corollary}\label{cor:relativistic-power-profiles}
For $\gamma\ge1$, the nonnegative relativistic equilibrium
\begin{equation}\label{eq:relativistic-profile-section6}
 \Phi_\gamma(v)=\frac{\gamma}{\gamma+1}c_\gamma
                    \ang v^{-2\gamma-2}
\end{equation}
is matched to the normalized profile
$f^0(u)=c_\gamma(1-|u|^2)_+^\gamma$, with $c_\gamma$ as defined after
\eqref{H0}.  Its response satisfies
Theorem~\ref{thm:relativistic-green}.  Under the transformed-data
assumptions of Theorem~\ref{thm:relativistic-electric}, its electric
field and solutions satisfy the corresponding field and scattering
conclusions.
\end{corollary}

\subsection{Relation with earlier work}

The whole-space wave--damping picture combines phase mixing with
spatial dispersion.  Its spectral background includes the analyses
of Case, Penrose, and Degond \cite{Case59,Deg86,Pen60}, while Glassey
and Schaeffer \cite{GS94,GS95} studied the sensitivity of long-time
decay to the equilibrium in unbounded space.  For Maxwellian
backgrounds, Bedrossian, Masmoudi, and Mouhot \cite{BMM22} separated
a weakly damped Klein--Gordon-type long-wave component from a faster
Landau-damped component and identified the corresponding long-wave
Euler--Poisson dynamics.  Han-Kwan, T.~T. Nguyen, and Rousset
\cite{HKNR21L} established physical-space dispersive estimates for
a class of stable analytic equilibria.  Ionescu, Pausader, Wang,
and Widmayer \cite{IPWW23} analyzed the linearized response for a
broad class of analytic equilibria, including explicit formulas for
generalized Poisson profiles.

For radial equilibria positive throughout a velocity ball, with
finite or infinite maximal speed, T.~T. Nguyen \cite{Ngu26} studied
the survival threshold and obtained Green-function and electric-field
decay estimates under finite smoothness and decay assumptions.
His main electric-field result gives time-decay estimates in spatial
$L_x^p$ norms, with $2\le p<\infty$ for the oscillatory component.

In the compactly supported setting, our equilibrium class allows
both a central velocity hole and lower boundary regularity than
in \cite{Ngu26}. The inner and outer vanishing orders are independent
and may be as low as one, so the equilibrium need not be $C^1$
across the velocity boundaries. This low-order boundary behavior
requires explicit control of endpoint singularities in the spectral
analysis and the spatial estimates.

Our main estimates are pointwise in the spatial variable: they retain
explicit spatial weights and quantify the dependence on the distances
to the fronts $|x|=at$, with decay rates determined by the individual
boundary orders. Theorem~\ref{thm:main-green} establishes these
Green-function bounds uniformly through pole--continuum contact.
For the complete electric field, Theorem~\ref{thm:main-electric}
gives corresponding weighted pointwise estimates, including a
$t^{-3/2}$ wave bound in $L_x^\infty$ with a cubic spatial weight
and no moment cancellation beyond neutrality.
This finer spatial resolution uses the data norm
\eqref{eq:data-norm}, with eleven total mixed derivatives and
joint weight $\ang{x,v}^{17}$.

For relativistic plasmas, Young \cite{Young15} studied linear damping
in periodic and unbounded space. Our relativistic results follow
from the exact correspondence proved in
Section~\ref{subsec:relativistic-correspondence}.

The nonlinear theory provides a broader context.  On the torus,
Mouhot and Villani \cite{MV11} proved Landau damping near stable
homogeneous equilibria in analytic regularity, with extensions to
certain Gevrey classes.  Bedrossian, Masmoudi, and Mouhot
\cite{BMM16} developed a proof using paraproducts and time-response
estimates, and Grenier, T.~T. Nguyen, and Rodnianski \cite{GNR21}
gave a resolvent-based approach.  Ionescu, Pausader, Wang, and
Widmayer \cite{IPWW24G} reached the weighted Gevrey-3 endpoint and
constructed nonlinear wave and scattering operators.  Bedrossian's
echo constructions \cite{Bed21} show why finite Sobolev regularity
alone does not support a comparable general damping theory.

In unbounded space with screening, Bedrossian, Masmoudi, and Mouhot
\cite{BMM18} established nonlinear damping for localized perturbations
with finite Sobolev regularity.  Han-Kwan, T.~T. Nguyen, and Rousset
\cite{HKNR21S} developed a physical-space approach, and Huang,
Q.-H. Nguyen, and Xu \cite{HNX25S} obtained sharp estimates in
dimensions $d\ge3$ without a time-logarithmic loss.  The latter
authors also treated the two-dimensional massless-electron model
\cite{HNX25M}.  Wang, Xiao, and Xiong \cite{WXX26} considered the
two-species screened system with large initial distributions under
a structural quasi-neutrality condition.  Screening changes the
low-frequency force law, so these results have a different linear
response from the unscreened system considered here.

For the unscreened nonlinear problem, Ionescu, Pausader, Wang, and
Widmayer \cite{IPWW24} proved asymptotic stability of the Poisson
equilibrium in three dimensions.  Q.-H. Nguyen, Wei, and Zhang
\cite{NWZ25} gave a new proof, and Ionescu and Q.-H. Nguyen
\cite{IN26} established damping and free-transport scattering near
the two-dimensional Poisson equilibrium.  These results concern small,
smooth, localized perturbations; neutrality is also imposed in the
two-dimensional result.

Other aspects of Vlasov--Poisson dynamics include the classical
existence and regularity theory of Pfaffelmoser and of Lions and
Perthame \cite{LP91,Pfa92}, the low-regularity local theory of Jeong
and Tae \cite{JT25}, and the discontinuity results of Chen, Jeong,
Q.-H. Nguyen, and Tae in $L_x^pL_v^\infty$ \cite{CJNT26}.
In one-dimensional periodic settings, the classical BGK waves
\cite{BGK57} and the traveling and quasi-periodic electron layers
of Roulley and Franzoi--Roulley \cite{FR26,Rou25} illustrate coherent
dynamics distinct from the relaxation studied here.

\subsection{Main ideas and organization}

Our proof combines Laplace--Fourier analysis of the causal Volterra
response with spatially resolved oscillatory integral estimates.
We first establish weighted pointwise bounds for the Green function,
including the velocity fronts and the pole--continuum transition,
and then estimate the complete electric field through a low-frequency
normal form and weighted estimates for the free-transport source.
These two stages require different arguments: the first concerns the
spectral structure of the response, while the second must also retain
the cancellations and space--time dependence of the source.

The most delicate spectral issue is uniformity at the survival
threshold.  Near contact, the continuum density can concentrate, and
its limiting mass must be compared with that of the approaching pole.
We keep both contributions in a common spectral measure.  Its time
oscillations recover the Green function, while its Cauchy transform
recovers the resolvent; weak convergence and uniform total-variation
bounds identify the same causal response at contact.  Rather than
requiring a uniform time-frequency contour deformation through the
contact point, we localize jointly in the spatial frequency $k$ and an
endpoint coordinate and integrate first in $k$.  The transverse Cauchy
inversion produces complex dispersion roots, with their full amplitudes
and nonnegative damping.

Two mechanisms then control the spatial inverse.  In the conic region,
group-velocity oscillations are combined with continuum damping.
Bounded-variation estimates control the amplitudes, while the endpoint
order determines where oscillation or damping supplies the decay.
In particular, near the limiting group velocity, low-order continuum
endpoints are treated by damping rather than by a smooth stationary-phase
expansion.  In the interior region $0\le|x|\le v_0t$, with $v_0>0$
fixed and small, the pole and continuum share a holomorphic phase and
amplitude.  A Cauchy--Green deformation estimates their combined
boundary integral before angular integration.  Keeping the spherical
representation avoids the factor $|x|^{-1}$ in radial inversion and
covers $x=0$.

Away from contact, the challenge is to retain the spatial fronts
rather than only global time decay.  Taking absolute values too early
in the spatial Fourier integral would lose the distance to $|x|=at$.
We first transform the velocity-endpoint variable, converting its order
of vanishing into decay in $tk$ and a remaining radial frequency power.
Spatial inversion then relates that power to the front distance.
For the Green function, pairing the positive and negative endpoints
before taking absolute values removes the borderline logarithm at
fourth-order vanishing.  For the electric field, the source amplitudes
also depend on the angular variable.  A spherical expansion retains
this dependence and leads to the logarithm at a second-order endpoint.

Passing to the complete electric field requires two low-frequency
cancellations.  A direct convolution estimate does not display these
cancellations at the singularity of the Poisson multiplier.
Following the time-integration-by-parts approach in
\cite[Lemma~2.5]{BMM22} and
\cite[Proposition~3.3 and Section~4.2]{Ngu26}, we integrate the smooth
pole convolution in time.  Its current-time boundary term combines
with the free density to give
$\chi_\ell(k)-\widehat K_1(0,k)=k^2b_0(k)$.
The pole-subtracted response independently factors as
$R(z,k)-P(z,k)=k^2\widetilde{\mathcal R}(z/k,k)$.
These factors compensate for the Coulomb singularity in the local and
history terms, respectively.  Neutrality controls the first initial
pole amplitude in polar coordinates, and the polar dispersive estimate
gives its cubic spatial weight without higher moment cancellation.

The remaining response term must be estimated together with the full
finite source history.  Frequency moments and physical derivatives
enter this step simultaneously, and differentiating the moving upper
limit produces additional boundary terms.  The Fourier-ray estimates
in Appendix~\ref{app:electric-source} control these terms with the finite
index $M_*=11$.  Frequencies separated from zero are handled by weighted
spatial convolution estimates, which preserve the spatial decay while
integrating over the source history.

Section~\ref{sec:dispersion} establishes the dispersion geometry and
boundary regularity, including the filled-ball center.
Section~\ref{sec:spectral-decomposition} constructs the spectral
representation used in the Green-function estimates of
Section~\ref{sec:green}.  Section~\ref{sec:electric} proves the complete
electric-field bounds, and Section~\ref{sec:consequences} proves
scattering and the relativistic applications.
Appendices~\ref{app:endpoints}--\ref{app:electric-source} collect the
oscillatory, transverse-inversion, and source estimates.

\paragraph{Notation.}
The relations $A\lesssim B$ and $A\asymp B$
denote an upper bound up to a constant and two-sided comparability,
respectively.  Unless specified otherwise, constants may depend on
the fixed equilibrium, its endpoint orders, the chosen cutoffs, and
the finite derivative orders used in the estimate.
For a positive endpoint distance $u$, we write
$E(u)=O_{\rm con}(u^\beta(1+|\log u|)^b)$ through order $J$ if
\[
 |D_u^jE(u)|\le C_j u^\beta(1+|\log u|)^b,
 \qquad D_u=u\partial_u,\quad 0\le j\le J.
\]
The identity
\[
 u^m\partial_u^m=\prod_{j=0}^{m-1}(D_u-j)
\]
converts conormal bounds into ordinary endpoint derivative bounds,
with each derivative lowering the endpoint power by one and retaining
the displayed logarithmic factor.  Signed endpoint coordinates are
understood separately on their two open sides.

\section{Dispersion function and survival threshold}
\label{sec:dispersion}

For the classical system \eqref{eq:linear-VP}, this section fixes the
transform normalization and determines the dispersion geometry used
in the sequel.  In particular, the
oscillatory branch below consists only of the imaginary-axis poles
present before their contact with the velocity continuum.
Throughout Sections~\ref{sec:dispersion}--\ref{sec:electric}, we treat
both \eqref{H1}--\eqref{H2} and \eqref{H0}, setting $\kappa=0$ in
the latter case.  Lemma~\ref{lem:A-endpoint-section2} treats the
velocity center together with the other interior points.

\subsection{Volterra reduction and the Laplace--Fourier transform}
\label{subsec:volterra-construction-v3}

In addition to \eqref{eq:fourier-transform-intro}, we use
\[
\widehat h_{\rm in}(\xi,\eta)
:=\iint_{\R^3\times\R^3}
\re^{-\ri(x\cdot\xi+v\cdot\eta)}
h_{\rm in}(x,v)\,\rd x\,\rd v .
\]
Under either \eqref{H1}--\eqref{H2} or \eqref{H0},
$f^0\in W^{1,\infty}(\R^3)$, and $\nabla_v f^0$ is interpreted almost
everywhere.  With the Newton-potential field
$E=-\nabla_x(-\Delta_x)^{-1}\rho$, we use Duhamel's formula for
\eqref{eq:linear-VP}, valid for almost every $(x,v)$:
\begin{equation}
\label{eq:duhamel-section2}
h(t,x,v)
=h_{\rm in}(x-tv,v)
+\int_0^t
\nabla_x\phi(s,x-v(t-s))\cdot\nabla_vf^0(v)\,\rd s .
\end{equation}
Equivalently, after taking the Fourier transform in $x$ only,
\[
\widehat h(t,\xi,v)
=\re^{-\ri t v\cdot\xi}\widehat h_{\rm in}(\xi,v)
+\frac{\ri\xi}{|\xi|^2}\cdot\nabla_vf^0(v)
\int_0^t\re^{-\ri(t-s)v\cdot\xi}
\widehat\rho(s,\xi)\,\rd s .
\]
Here $\widehat h_{\rm in}(\xi,v)$ denotes the Fourier transform in $x$
only.  The sign follows from
\(\widehat E=-\ri\xi|\xi|^{-2}\widehat\rho\) and from the minus sign
introduced when the force term is integrated along free
characteristics.  Integrating this identity in \(v\) gives, for
\(\xi\ne0\),
\begin{equation}
\label{eq:volterra-density}
\widehat\rho(t,\xi)
=\widehat h_{\rm in}(\xi,t\xi)
+\int_0^t \widehat K(t-s,\xi)\widehat\rho(s,\xi)\,\rd s ,
\end{equation}
where
\begin{equation}
\label{eq:kernel-K}
\widehat K(t,\xi)
:=\frac{\ri\xi}{|\xi|^2}\cdot
\widehat{\nabla_vf^0}(t\xi)
=-t\,\widehat{f^0}(t\xi),\qquad t\geq0.
\end{equation}
The last equality uses
\(\widehat{\nabla_vf^0}(\eta)=\ri\eta\widehat{f^0}(\eta)\).  It is at
this point that the sign of the Poisson coupling enters the scalar
Volterra equation.
The first term on the right of \eqref{eq:volterra-density} is the
freely transported density:
\[
\widehat{\rho^{\rm free}}(t,\xi)
=\widehat h_{\rm in}(\xi,t\xi).
\]
Since \(f^0\) is radial, \(\widehat K(t,\xi)\) depends on \(\xi\) only through
\(k=|\xi|\); below we write \(\widehat K(t,k)\).

Equation~\eqref{eq:volterra-density} is the Volterra equation
\eqref{eq:volterra-intro} with its kernel written explicitly.  We first construct its response through
\eqref{eq:green-volterra-definition}, before taking any time transform.
For fixed \(k\), write
\[
(A*_tB)(t,k):=\int_0^t A(t-s,k)B(s,k)\,\rd s.
\]
The $n$-fold causal time convolution is defined by
\[
(\widehat K)^{*1}=\widehat K,\qquad
(\widehat K)^{*(n+1)}=\widehat K*_t(\widehat K)^{*n}.
\]
Iteration of \(\widehat G=\widehat K+\widehat K*_t\widehat G\) leads to
\(\widehat G=\sum_{n=1}^{\infty}(\widehat K)^{*n}\).
The nonnegativity of $f^0$ and the mass normalization
\eqref{eq:radial-equilibrium} give $|\widehat K(t,k)|\le t$.
The inductive step is
\[
|(\widehat K)^{*(n+1)}(t,k)|
\le\int_0^t(t-s)\frac{s^{2n-1}}{(2n-1)!}\,\rd s
=\frac{t^{2n+1}}{(2n+1)!}.
\]
Consequently,
\begin{equation}\label{eq:volterra-series-section2}
 |(\widehat K)^{*n}(t,k)|\le\frac{t^{2n-1}}{(2n-1)!},\qquad
 \widehat G(t,k)=\sum_{n\ge1}(\widehat K)^{*n}(t,k),\qquad
 |\widehat G(t,k)|\le\sinh t.
\end{equation}
The series converges uniformly on compact time intervals, so it can
be convolved term by term with \(\widehat K\) and satisfies
\eqref{eq:green-volterra-definition}.
For two solutions on \([0,T]\), iteration of their difference
equation bounds the difference by \(T^{2n}/(2n)!\) times its supremum.
Letting \(n\to\infty\) proves uniqueness.

To recover the density, the resolvent identity gives
\[
 \widehat K*_t\bigl(\widehat{\rho^{\rm free}}
       +\widehat G*_t\widehat{\rho^{\rm free}}\bigr)
 =\widehat G*_t\widehat{\rho^{\rm free}}.
\]
Thus $\widehat{\rho^{\rm free}}+\widehat G*_t\widehat{\rho^{\rm free}}$
solves \eqref{eq:volterra-density}.
Associativity is justified on the finite time simplex because
$\widehat{\rho^{\rm free}}$ is bounded and $\widehat K,\widehat G$
are continuous on bounded intervals.
The same uniqueness argument therefore yields
\begin{equation}\label{eq:rho-G-free-section2}
 \widehat\rho(t,\xi)=\widehat{\rho^{\rm free}}(t,\xi)
       +\int_0^t\widehat G(t-s,k)\widehat{\rho^{\rm free}}(s,\xi)\,\rd s,
       \qquad k=|\xi|.
\end{equation}
This proves the density formula in
\eqref{eq:green-field-representation-intro}.
Since $|\widehat{\rho^{\rm free}}|\le\|h_{\rm in}\|_{L^1_{x,v}}$,
the series bound \eqref{eq:volterra-series-section2} also gives
\begin{equation}\label{eq:exp-bound-rho}
 |\widehat\rho(t,\xi)|
 \le\|h_{\rm in}\|_{L^1_{x,v}}\cosh t.
\end{equation}

We use \(\Lap\) for the Laplace transform in time, with
\[
(\Lap f)(z):=\int_0^\infty\re^{-zt}f(t)\,\rd t,\qquad
\widetilde f(z)=(\Lap f)(z)
\]
where the integral converges.  For spatial Fourier transforms the
same convention means
\[
\widetilde\rho(z,\xi)
=(\Lap[\widehat\rho(\cdot,\xi)])(z),\qquad
\widetilde G(z,k)=(\Lap[\widehat G(\cdot,k)])(z);
\]
the spatial frequencies are parameters.  In particular, the scalar
dispersion function is the transform of \(\widehat K\):
\begin{equation}
\label{eq:dispersion-Laplace}
\Lcal(z,k)
:=(\Lap[\widehat K(\cdot,k)])(z)
=\int_0^\infty\re^{-zt}\widehat K(t,k)\,\rd t,\qquad \Re z>0.
\end{equation}
The bound \(|\widehat K(t,k)|\le t\) makes \(\Lcal(\cdot,k)\) analytic on
this half-plane and gives \(|\Lcal(z,k)|\le(\Re z)^{-2}\).
Our inverse convention is the Bromwich formula
\begin{equation}
\label{eq:bromwich-convention-section2}
 f(t)=(\Lap^{-1}\widetilde f)(t)
 =\frac1{2\pi\ri}\int_{\sigma-\ri\infty}^{\sigma+\ri\infty}
 \re^{zt}\widetilde f(z)\,\rd z,\qquad \sigma>0,\quad t>0,
\end{equation}
with the line chosen in the half-plane of convergence, first when the
integral is classical and otherwise by the usual limiting procedure.

The bounds \eqref{eq:volterra-series-section2} and
\eqref{eq:exp-bound-rho} justify the actual transforms of the response
and density on \(\Re z>1\).
Taking \(\Lap\) in \eqref{eq:volterra-density} gives
\begin{equation}
\label{eq:rho-resolvent-section2}
\widetilde\rho(z,\xi)
=Q(z,k)\widetilde{\rho^{\rm free}}(z,\xi),\qquad
Q(z,k)=\frac1{1-\Lcal(z,k)},\qquad R(z,k)=Q(z,k)-1.
\end{equation}
On the same half-plane, transforming the absolutely convergent
Volterra series gives
\[
\sum_{n\ge1}\Lcal(z,k)^n
=\frac{\Lcal(z,k)}{1-\Lcal(z,k)};
\]
hence
\begin{equation}\label{eq:far-right-response-v3}
 \widetilde G(z,k)=\int_0^\infty \re^{-zt}\widehat G(t,k)\,\rd t
       =\frac{\Lcal(z,k)}{1-\Lcal(z,k)}=R(z,k).
\end{equation}
Thus \(Q\) is the full density resolvent and \(R\) is its response
part.  The constant \(1\) in \(Q=1+R\) inverts to \(\delta_{t=0}\),
which supplies the free-density term in \eqref{eq:rho-G-free-section2}.

Bromwich inversion on a line \(\Re z>1\) recovers the response
already constructed in time.  The stability result in
Lemma~\ref{lem:right-half-plane-stability-section2} below extends
\(Q\) and \(R\) analytically to \(\Re z>0\).
For \(0<\sigma\le\Re z\le\sigma_0\),
\eqref{eq:dispersion-vertical-decay-section2} gives
\(R(z,k)=O_{\sigma,\sigma_0,k}(|\Im z|^{-2})\).
The horizontal connecting integrals therefore vanish, and the
Bromwich line can be moved to any \(\Re z=\sigma>0\):
\begin{equation}\label{eq:Green-fun-G-t-k}
 \widehat G(t,k)=\frac1{2\pi\ri}
       \int_{\sigma-\ri\infty}^{\sigma+\ri\infty}\re^{zt}R(z,k)\,\rd z.
\end{equation}
This is a deformation of the inverse contour for the same response.
Absolute convergence of its original time-transform integral on
the whole right half-plane is obtained later: spectral reconstruction
gives fixed-\(k\) boundedness in
Corollary~\ref{cor:fixed-k-boundedness-section3}, which also establishes
the corresponding actual Laplace transform of the density.

The mass normalization in
\eqref{eq:radial-equilibrium} also gives the continuous \(k=0\)
extensions
\begin{equation}
\label{eq:kzero-resolvent}
\Lcal(z,0)=-\frac1{z^2},\qquad
Q(z,0)=\frac{z^2}{z^2+1},\qquad
R(z,0)=-\frac1{z^2+1}.
\end{equation}
In particular, the two limiting pole weights are
\(\ri/2\) at \(z=\ri\) and \(-\ri/2\) at \(z=-\ri\), and the inverse
response $\widehat G(t,0)$ is \(-\sin t\) for $t>0$. Equivalently, its forward extension to the whole time axis is $\widehat G(t,0)=-\one_{\{t>0\}}\sin t$. This elementary case is a useful check on all
signs below.
\subsection{Radial reduction and boundary values}
\label{subsec:radial-reduction-v199}

Put \(F(r)=\mu(r^2/2)\), and define the one-dimensional marginal
\begin{equation}
\label{eq:mu0-section2}
\mu_0(u)
:=\int_{\R^2}F\bigl((|v'|^2+u^2)^{1/2}\bigr)\,\rd v'
=\pi\int_{u^2}^{\infty}\mu(r/2)\,\rd r .
\end{equation}
It is even, supported in \([-1,1]\), and
\begin{equation}
\label{eq:mu0-prime-section2}
\mu_0'(u)=-2\pi u\mu(u^2/2).
\end{equation}
In particular, $\mu_0$ is continuously differentiable and $\mu_0'$
is locally Lipschitz, also when $\gamma_{\min}=1$.
The following identities are justified by smooth approximation.

Aligning \(\xi\) with the third coordinate axis in
\eqref{eq:dispersion-Laplace}, say \(\xi=ke_3\), gives
\begin{align*}
\Lcal(z,k)=\int_0^\infty\int_\R
\re^{-(z+\ri ku)t}\frac{\ri}{k}\mu_0'(u)\,\rd u\,\rd t=\frac{\ri}{k}\int_\R\frac{\mu_0'(u)}{z+\ri ku}\,\rd u .
\end{align*}
Fubini's theorem is justified by \(\operatorname{Re}z>0\) and the
compact support of \(\mu_0'\).  Using
\(z+\ri ku=\ri k(u-\ri z/k)\), we obtain a Cauchy integral
representation of the same dispersion function defined in
\eqref{eq:dispersion-Laplace}:
\begin{equation}
\label{eq:analytic-cauchy-section2}
\Lcal(z,k)
=\frac1{k^2}\int_{\R}
\frac{\mu_0'(u)}{u-\ri z/k}\,\rd u,
\qquad k>0,\quad \operatorname{Re}z>0 .
\end{equation}
Since $\mu_0$ is compactly supported, an integration by parts in
\eqref{eq:analytic-cauchy-section2} gives the useful representation
\begin{equation}
	\label{eq:dispersion-decay-section2}
	\Lcal(z,k)
	=-\int_{\R}\frac{\mu_0(u)}{(z+\ri ku)^2}\,\rd u .
\end{equation}
In particular, for every $\sigma>0$,
\begin{equation}
	\label{eq:dispersion-vertical-decay-section2}
	|\Lcal(\sigma+\ri\tau,k)|
	\le
	\int_{\R}
	\frac{\mu_0(u)}
	{\sigma^2+(\tau+ku)^2}\,\rd u
	\le C_{\sigma,k}\ang{\tau}^{-2}.
\end{equation}
The constant is uniform for $k$ in compact subsets of $[0,\infty)$.
Indeed, when $|\tau|\ge2k$ the last integral is bounded by
$4|\tau|^{-2}$, while for bounded $\tau$ one uses
$|\Lcal|\le\sigma^{-2}$. 

To describe its boundary values, we separate the real principal-value
part from the jump across the velocity continuum.  For real \(s\), define
\begin{equation}
\label{eq:A-B-section2}
\Acal(s)
:=\pv\int_\kappa^1
\frac{u^2\mu(u^2/2)}{s^2-u^2}\,\rd u,
\qquad
\Bcal(s):=\pi s\mu(s^2/2),
\qquad
\Om(s):=4\pi\Acal(s)+2\pi\ri\Bcal(s).
\end{equation}
The integral in \eqref{eq:analytic-cauchy-section2} continues
\(\Lcal\) to the other side of its slit.  Write
\[
\Lcal_\pm(\tau,k)
:=\lim_{\eta\downarrow0}\Lcal(\pm\eta+\ri\tau,k),\qquad s=\tau/k,
\]
where the minus trace is taken in this continuation.
To compute the plus trace directly from
\eqref{eq:analytic-cauchy-section2}, put \(z=\eta+\ri\tau\) and
use the oddness of \(\mu_0'\) to change \(u\) to \(-u\):
\[
\Lcal(\eta+\ri\tau,k)
=\frac1{k^2}\int_\R
  \frac{\mu_0'(u)}{u-s+\ri\eta/k}\,\rd u.
\]
The Plemelj formula \cite[Part~I, Chapter~2]{Mus58}
\[
\frac1{x+\ri0}=\pv\frac1x-\ri\pi\delta_0
\]
therefore gives
\[
\Lcal_+(\tau,k)
=\frac1{k^2}\left[
\pv\int_\R\frac{\mu_0'(u)}{u-s}\,\rd u-\ri\pi\mu_0'(s)\right].
\]
Pairing \(u\) and \(-u\) in the principal value, and using
\eqref{eq:mu0-prime-section2}, yields
\[
\pv\int_\R\frac{\mu_0'(u)}{u-s}\,\rd u
=4\pi\pv\int_\kappa^1\frac{u^2\mu(u^2/2)}{s^2-u^2}\,\rd u
=4\pi\Acal(s),\qquad
-\ri\pi\mu_0'(s)=2\pi\ri\Bcal(s).
\]
Thus the two boundary traces are
\begin{equation}
\label{eq:Omega-boundary-section2}
\Lcal_+(\tau,k)=\frac{\Om(s)}{k^2},
\qquad
\Lcal_-(\tau,k)=\frac{\overline{\Om(s)}}{k^2}.
\end{equation}
These are the traces of the same Cauchy integral, with opposite
sides of the slit.  Since \(\Acal\) is even and \(\Bcal\) is odd,
\[
\Lcal_+(-\tau,k)
=\overline{\Lcal_+(\tau,k)}=\Lcal_-(\tau,k).
\]
The jump of the scalar response resolvent is supported on
\begin{equation}
\label{eq:physical-cut-section2}
\Sigma_k=[-k,-\kappa k]\cup[\kappa k,k],
\end{equation}
with a central gap $(-\kappa k,\kappa k)$ when $\kappa>0$.
Under \eqref{H0} the cut is $[-k,k]$.

We use the same subscripts for the boundary traces of \(Q\) and \(R\).
Since $R$ is the transform inverted in \eqref{eq:Green-fun-G-t-k},
its jump supplies the density of the continuum time oscillations.
We first compute this density from the two boundary fractions, and
then explain the reverse Cauchy-transform construction.  Set
\[
X(s,k):=k^2-4\pi\Acal(s),\qquad
Y(s):=2\pi\Bcal(s).
\]
Away from the endpoints of \(\Sigma_k\),
\begin{equation}
\label{eq:jump-density-section2}
\begin{split}
R_+(\tau,k)-R_-(\tau,k)
&=\frac{2\ri k^2Y(s)}{X(s,k)^2+Y(s)^2},\\
J_k(\tau)
:=\frac1{2\pi}\bigl(R_+(\tau,k)-R_-(\tau,k)\bigr)
&=\frac{2\ri k^2\Bcal(s)}
{\bigl(k^2-4\pi\Acal(s)\bigr)^2+
\bigl(2\pi\Bcal(s)\bigr)^2}.
\end{split}
\end{equation}
Indeed,
\begin{equation}\label{eq:resolvent-boundary-section2}
Q_+(\tau,k)=\frac{k^2}{X(s,k)-\ri Y(s)},
\qquad
Q_-(\tau,k)=\frac{k^2}{X(s,k)+\ri Y(s)},
\end{equation}
and subtracting these two fractions proves
\eqref{eq:jump-density-section2}.  The ordinary jump vanishes away
from \(\Sigma_k\) wherever the resolvent is regular; isolated poles
will contribute separately.  Since \(R=Q-1\), the two resolvents
have the same jump, and \(J_k\) is purely imaginary and odd.

The factor \(1/(2\pi)\) in \(J_k\) is fixed by our Bromwich
convention.  For a single pole,
\begin{equation}
\label{eq:bromwich-model-section2}
\left(\Lap^{-1}\Big[\frac{1}{z-\ri\omega}\Big]\right)(t)
=\frac1{2\pi\ri}\int_{\sigma-\ri\infty}^{\sigma+\ri\infty}
\frac{\re^{zt}}{z-\ri\omega}\,\rd z
=\re^{\ri\omega t},\qquad t>0.
\end{equation}
Hence a term \(a/(z-\ri\omega)\) contributes
\(a\re^{\ri\omega t}\): its spectral atom has weight \(a\), the
residue, at the temporal frequency \(\omega\).

For comparison, take a smooth compactly supported, possibly complex
test density \(J\) and define its Cauchy transform
\[
C_J(z):=\int_\R\frac{J(\tau)}{z-\ri\tau}\,\rd\tau.
\]
Here \(J\) is the density and \(C_J\) is the function it generates in
the transform domain.  The distributional kernel identity
\[
\frac1{\eta+\ri x}-\frac1{-\eta+\ri x}
\longrightarrow2\pi\delta_0(x),\qquad \eta\downarrow0,
\]
gives \((C_J)_+-(C_J)_-=2\pi J\).  Inverting each Cauchy kernel
with \eqref{eq:bromwich-model-section2} gives the causal inverse
\[
(\Lap^{-1}C_J)(t)
=\one_{\{t>0\}}\int_\R\re^{\ri t\tau}J(\tau)\,\rd\tau.
\]
Thus \(J\) is exactly the coefficient density of the time oscillations.

Applied to \eqref{eq:jump-density-section2}, this normalization
identifies \(J_k(\tau)\,\rd\tau\) as the continuum part of the
spectral measure.  An isolated pole at \(\ri\omega\) contributes
\((\Res_{z=\ri\omega}R)\delta_\omega\).
The reverse reconstruction of $R$ therefore requires the pole atoms
as well as the continuum density.  Proposition~\ref{prop:fixed-k-representation-section3}
gives its Cauchy representation away from $k_*$, and
Proposition~\ref{prop:terminal-reconstruction-section3} includes the
pole--continuum contact.

The endpoint behavior of \(\Acal\) determines the contact of the pole
with the continuum.  We record its interior regularity, endpoint
expansion, and large-$s$ behavior.  For the latter, write
\begin{equation}\label{eq:radial-moments-section2}
m_{2j}:=\int_\kappa^1 u^{2j}\mu(u^2/2)\,\rd u,\qquad j\ge1,
\qquad 4\pi m_2=1.
\end{equation}
These are the radial moments of the equilibrium; the last identity
is the mass normalization \eqref{eq:radial-equilibrium}.
For the endpoint expansions let
\[
\mu_{\rm out}(r)
:=\frac{\mu(r^2/2)}{(1-r)^{\gamma_1}}
\quad\hbox{near }r=1,
\qquad
\mu_{\rm in}(r)
:=\frac{\mu(r^2/2)}{(r-\kappa)^{\gamma_\kappa}}
\quad\hbox{near }r=\kappa,
\]
using the smooth positive extensions in \eqref{H2}, or the outer
extension in \eqref{H0}.  All statements below involving the inner
endpoint $\kappa$ or its order $\gamma_\kappa$ apply only when
$\kappa>0$.

\begin{lemma}
\label{lem:A-endpoint-section2}
Under \eqref{H1}--\eqref{H2} or \eqref{H0}, the function $\Acal$
is bounded and even.  Both $\Acal$ and $\Bcal$ are smooth away
from the velocity endpoints: $\pm\kappa,\pm1$ in the shell case,
and only $\pm1$ under \eqref{H0}.  The function $\Acal$ is
{\(C^{\gamma_{\min}}\) when
\(\gamma_{\min}\notin\mathbb N\), and \(C^{\gamma_{\min}-}\) when
\(\gamma_{\min}\in\mathbb N\)}.  Here \(C^{\gamma_{\min}-}\)
means membership in every strictly lower H\"older class \(C^s\),
\(0<s<\gamma_{\min}\).  For \(s>1\),
\begin{equation}
\label{eq:A-monotone-section2}
\Acal(s)>0,\qquad \Acal'(s)<0,\qquad
\lim_{s\to\infty}\Acal(s)=0.
\end{equation}
If \({\gamma_1}\notin\mathbb N\), then, near \(s=1\),
\begin{equation}
\label{eq:A-outer-noninteger-section2}
\Acal(s)+\frac{s\mu_{\rm out}(s)}2
\left[
\frac{\pi}{\sin(\pi{\gamma_1})}\one_{\{s\geq1\}}(s-1)^{\gamma_1}
+\frac{\pi\cos(\pi{\gamma_1})}{\sin(\pi{\gamma_1})}
\one_{\{s<1\}}(1-s)^{\gamma_1}
\right]
\quad\hbox{is smooth}.
\end{equation}
If \({\gamma_1}=m\in\mathbb N\), then
\begin{equation}
\label{eq:A-outer-integer-section2}
\Acal(s)+\frac{s\mu_{\rm out}(s)}2
(1-s)^m\log|s-1|
\quad\hbox{is smooth near }s=1.
\end{equation}
When $\kappa>0$, the inner endpoint has the same power/logarithmic
model with order $\gamma_\kappa$.  In particular, with $D_u=u\partial_u$, for every
fixed $n\ge0$ and $\sigma\in\{-1,1\}$,
\begin{equation}\label{eq:A-inner-conormal-section2}
 \bigl|D_u^n\bigl(\Acal(\kappa+\sigma u)-\Acal(\kappa)\bigr)\bigr|
 \le C_nu\left(1+\one_{\{\gamma_\kappa=1\}}\log\frac1u\right),
 \qquad 0<u<u_0.
\end{equation}
Here $u_0$ is chosen sufficiently small that the chart does not
meet another velocity endpoint.
Finally, for every \(N\geq1\),
\begin{equation}
\label{eq:A-large-section2}
\Acal(s)
=\sum_{j=1}^{N}\frac{m_{2j}}{s^{2j}}
+O(|s|^{-2N-2}),\qquad |s|\to\infty.
\end{equation}
Under \eqref{H0}, in particular, $\Acal$ and $\Bcal$ are smooth
at $s=0$, and
\begin{equation}\label{eq:filled-ball-zero-sign-appendix}
 \Acal(0)=-\int_0^1\mu(u^2/2)\,\rd u<0,\qquad\Bcal(0)=0.
\end{equation}
\end{lemma}

\begin{proof}
We first establish smoothness away from the velocity endpoints.
The identity \eqref{eq:mu0-prime-section2} shows that $\mu_0'$ is
smooth near every such point, including zero under \eqref{H0}.
Fix one of these points $s_0$ and choose $a<s_0<b$ so that
$\mu_0'$ is smooth on $[a,b]$.  The Hilbert-transform representation
above gives, for $s\in(a,b)$,
\[
\begin{aligned}
 4\pi\Acal(s)
 ={}&\int_a^b\frac{\mu_0'(u)-\mu_0'(s)}{u-s}\,\rd u
   +\mu_0'(s)\log\frac{b-s}{s-a}
 +\int_{\R\setminus[a,b]}\frac{\mu_0'(u)}{u-s}\,\rd u.
\end{aligned}
\]
The first quotient extends smoothly across $u=s$, since
\[
 \frac{\mu_0'(u)-\mu_0'(s)}{u-s}
 =\int_0^1\mu_0''(s+\theta(u-s))\,\rd\theta.
\]
The logarithm is smooth near $s_0$.  In the last integral the
denominator stays away from zero there, and $\mu_0'$ is compactly
supported and integrable, so differentiation under the integral is
valid to every fixed order.  This proves the asserted smoothness
of $\Acal$; that of $\Bcal$ follows from \eqref{eq:A-B-section2}.
Under \eqref{H0}, evaluating the same definition at zero gives
\eqref{eq:filled-ball-zero-sign-appendix}.

Choose a smooth cutoff supported near \(u=1\) and equal to one on a
smaller neighborhood.  The complementary part of \(\Acal\) is smooth
for \(s\) near \(1\).  On the localized part write
\[
g(s,u):=\frac{u^2\mu_{\rm out}(u)}{s+u},\qquad
\frac{u^2\mu(u^2/2)}{s^2-u^2}
=g(s,u)\frac{(1-u)^{\gamma_1}}{s-u}.
\]
The identity
\[
g(s,u)=g(s,s)+(u-s)
\int_0^1\partial_2g\bigl(s,s+\theta(u-s)\bigr)\,\rd\theta
\]
shows that the term containing the difference quotient is smooth in
\(s\).  Since \(g(s,s)=s\mu_{\rm out}(s)/2\), the only nonsmooth
factor is the model integral
\begin{equation}
\label{eq:toy-Hilbert-section2}
T_{\gamma_1}(s)
:=\pv\int_0^1\frac{y^{\gamma_1}}{y+s-1}\,\rd y .
\end{equation}
Changing the upper limit of the $y$-integral to another fixed positive
value changes only a smooth function of $s$ near $s=1$.

Let \(a=s-1\).  If \({\gamma_1}\notin\mathbb N\), put
\(n=\lfloor{\gamma_1}\rfloor\) and repeatedly use
\[
\frac1{y+a}
=\sum_{j=0}^{n}\frac{(-a)^j}{y^{j+1}}
+\frac{(-a)^{n+1}}{y^{n+1}(y+a)}.
\]
The sum integrates to a polynomial in \(a\).  In the last term, the
change of variables $y=|a|r$ isolates $|a|^{\gamma_1}$.
With $\alpha:=n+1-\gamma_1\in(0,1)$, its two one-sided coefficients
are determined by
\[
\int_0^\infty\frac{r^{-\alpha}}{1+r}\,\rd r
=\frac{\pi}{\sin(\pi\alpha)},\qquad
\pv\int_0^\infty\frac{r^{-\alpha}}{r-1}\,\rd r
=\pi\cot(\pi\alpha),\qquad 0<\alpha<1.
\]
Here $\sin(\pi\alpha)=(-1)^n\sin(\pi\gamma_1)$ and
$\cot(\pi\alpha)=-\cot(\pi\gamma_1)$; extending the rescaled upper
limit to infinity changes only the smooth part.
Together with the factor $(-a)^{n+1}$, these identities prove
\eqref{eq:A-outer-noninteger-section2}.
For $\gamma_1=m\in\mathbb N$, the same division gives the singular
term $-(1-s)^m\log|s-1|$, yielding \eqref{eq:A-outer-integer-section2}.

When $\kappa>0$, translating the integration variable at the inner
endpoint by $\kappa$ gives the same half-line model with reflected argument,
which yields \eqref{eq:A-inner-conormal-section2} because its
conormal derivatives retain the endpoint order and
$u^\gamma(1+|\log u|)=O(u)$ for $\gamma>1$.
The same model proves the stated inner regularity, and evenness gives
the negative endpoints.

For \(s>1\), differentiation is classical and gives
\[
\Acal'(s)
=-2s\int_\kappa^1
\frac{u^2\mu(u^2/2)}{(s^2-u^2)^2}\,\rd u<0,
\]
which proves \eqref{eq:A-monotone-section2}.  Finally, for
\(|s|>1\),
\[
\frac1{s^2-u^2}
=\sum_{j=0}^{N-1}\frac{u^{2j}}{s^{2j+2}}
+\frac{u^{2N}}{s^{2N}(s^2-u^2)}.
\]
Integration proves \eqref{eq:A-large-section2}, with the same expansion
after any fixed number of \(s\)-derivatives.
\end{proof}

The other boundary function has no principal-value contribution.  In
particular, for \(s\) near \(1\),
\begin{equation}
\label{eq:B-outer-section2}
\Bcal(s)
=\pi s\mu_{\rm out}(s)(1-s)^{\gamma_1}\one_{\{s<1\}},
\end{equation}
and the analogous formula holds at \(s=\kappa\).  Equations
\eqref{eq:A-outer-noninteger-section2},
\eqref{eq:A-outer-integer-section2}, and
\eqref{eq:B-outer-section2} are the only velocity-endpoint expansions
used in the analysis near contact below.

\subsection{Stability and the pole branch}

We first exclude zeros of the dispersion denominator in the open
right half-plane.  This is the stability input for the contour
deformation in the spectral inversion below.

\begin{lemma}
\label{lem:right-half-plane-stability-section2}
For every \(k>0\), \(1-\Lcal(z,k)\) has no zero in
\(\operatorname{Re}z>0\).
\end{lemma}

\begin{proof}
Suppose that \(1-\Lcal(z,k)=0\), and put
\(a=\ri z/k=a_1+\ri a_2\), where \(a_2>0\).  By
\eqref{eq:analytic-cauchy-section2},
\begin{equation}
\label{eq:stability-complex-equation-section2}
k^2=\int_\R\frac{\mu_0'(u)}{u-a}\,\rd u .
\end{equation}
Since \(\mu_0\) is even and
\(\mu_0'(u)=-2\pi u\mu(u^2/2)<0\) on \((\kappa,1)\), the imaginary part
of the right-hand side equals
\[
a_2\int_0^\infty\mu_0'(u)
\left[
\frac1{(u-a_1)^2+a_2^2}
-\frac1{(u+a_1)^2+a_2^2}
\right]\rd u .
\]
It has sign \(-\operatorname{sgn}a_1\) when \(a_1\ne0\), and therefore
can vanish only if \(a_1=0\).  In that case the right-hand side of
\eqref{eq:stability-complex-equation-section2} is
\[
2\int_0^\infty\frac{u\mu_0'(u)}{u^2+a_2^2}\,\rd u<0,
\]
which cannot equal \(k^2\).
\end{proof}

The possible imaginary-axis zeros can now be read directly from
\eqref{eq:Omega-boundary-section2}.  Define
\begin{equation}
\label{eq:kstar-section2}
k_*:=\sqrt{4\pi\Acal(1)}>0.
\end{equation}

The monotonicity of $\Acal$ outside the velocity support now
determines the pole branch and its residues.  Its expansion at small
$k$ 
uses the radial moments $m_4$ and $m_6$ in \eqref{eq:radial-moments-section2}.

\begin{proposition}
\label{prop:pole-branch-section2}
For every \(0<k<k_*\), there is a unique \(s_{\rm p}(k)>1\) satisfying
\begin{equation}
\label{eq:pole-equation-section2}
k^2=4\pi\Acal(s_{\rm p}(k)).
\end{equation}
The only imaginary-axis zeros of \(1-\Lcal\) off the cut \eqref{eq:physical-cut-section2} are
\[
z=\pm\ri\tau_{\rm p}(k),\qquad
\tau_{\rm p}(k):=k\,s_{\rm p}(k)>k .
\]
They are simple.  At \(k=k_*\) the positive and negative zeros reach
the endpoints \(\pm\ri k_*\) of the cut, while for \(k>k_*\) there is
no imaginary-axis zero. For every \(0<k<k_*\), the Bromwich residues are
\begin{equation}
\label{eq:bromwich-residues-section2}
\varrho_\sigma(k)
:=\Res_{z=\sigma\ri\tau_{\rm p}(k)}R(z,k)
=\sigma\ri\,r_{\rm p}(k),\qquad
r_{\rm p}(k)
:=-\frac{k^3}{4\pi\Acal'(s_{\rm p}(k))}>0,
\quad \sigma=\pm .
\end{equation}
Both $\tau_{\rm p}(k)$ and $r_{\rm p}(k)$ extend as smooth functions
of $k^2$ near zero, with $\tau_{\rm p}(0)=1$ and $r_{\rm p}(0)=1/2$.
Moreover,
\begin{align}
\tau_{\rm p}(k)
&=1+2\pi m_4k^2+
\bigl(2\pi m_6-10\pi^2m_4^2\bigr)k^4+O(k^6),
\label{eq:taup-low-section2}\\
\varrho_+(k)
&=\ri\left[
\frac12-\pi m_4k^2+
\bigl(15\pi^2m_4^2-3\pi m_6\bigr)k^4+O(k^6)
\right],
\qquad
\varrho_-=\overline{\varrho_+}.
\label{eq:residue-low-section2}
\end{align}
In particular,
\begin{equation}
\label{eq:taup-hessian-section2}
D_\xi^2\bigl[\tau_{\rm p}(|\xi|)\bigr]_{\xi=0}
=4\pi m_4 I_3>0.
\end{equation}
\end{proposition}

\begin{proof}
For \(s>1\), \eqref{eq:A-monotone-section2} shows that
\(4\pi\Acal\) decreases bijectively from \(k_*^2\) to \(0\), which
proves existence and uniqueness in \eqref{eq:pole-equation-section2}.
Inside $\kappa<|s|<1$, the imaginary part $\Bcal(s)$ does not
vanish, while for $|s|\le\kappa$ in the shell case one has
$\Acal(s)<0$.
Under \eqref{H0}, the remaining point $s=0$ is covered by
\eqref{eq:filled-ball-zero-sign-appendix}.  Hence there are no other
imaginary-axis zeros.
  Along the positive branch,
\[
\partial_z\Lcal(\ri\tau_{\rm p}(k),k)
=-\frac{4\pi\ri}{k^3}\Acal'(s_{\rm p}(k)).
\]
Since \(R=\Lcal/(1-\Lcal)\) and \(\Lcal=1\) at the zero,
\[
\Res_{z=\ri\tau_{\rm p}}R
=-\frac1{\partial_z\Lcal(\ri\tau_{\rm p},k)}
=-\frac{\ri k^3}{4\pi\Acal'(s_{\rm p})}
=\ri r_{\rm p}(k).
\]
The symmetry $\Lcal(\bar z,k)=\overline{\Lcal(z,k)}$, which follows
from \eqref{eq:dispersion-decay-section2} and the evenness of $\mu_0$,
gives $\varrho_-=\overline{\varrho_+}$.  This proves
\eqref{eq:bromwich-residues-section2}; in particular, the factor
\(\ri\) comes from the Bromwich variable and agrees with
\eqref{eq:kzero-resolvent}.

For the smooth dependence at $k=0$, put $\lambda=k^2$ and define
\[
 \mathscr D(\tau,\lambda)
 :=1-4\pi\int_\kappa^1
       \frac{u^2F(u)}{\tau^2-\lambda u^2}\,\rd u.
\]
Its denominator is uniformly separated from zero near $(1,0)$,
so all fixed-order parameter derivatives may be taken under the
integral.  Mass normalization gives
$\mathscr D(1,0)=0$ and $\partial_\tau\mathscr D(1,0)=8\pi m_2=2$.
The smooth implicit function theorem therefore yields
$\tau_{\rm p}(k)=\mathsf T(k^2)$ for $k$ near zero, where
$\mathsf T(0)=1$.  Differentiating the integral for $\Acal$ gives
\[
 \Acal'(\tau/k)
 =-2\tau k^3\int_\kappa^1
          \frac{u^2F(u)}{(\tau^2-k^2u^2)^2}\,\rd u.
\]
Consequently \eqref{eq:bromwich-residues-section2} reads
\[
 r_{\rm p}(k)
 =\frac1{\partial_\tau\mathscr D(\tau_{\rm p}(k),k^2)}
 =:\mathsf R(k^2),\qquad \mathsf R(0)=\frac12,
\]
with $\mathsf R$ smooth near zero.

To compute the Taylor coefficients, insert
$s_{\rm p}=\tau_{\rm p}/k$ into \eqref{eq:A-large-section2}
and use \eqref{eq:pole-equation-section2} together with $4\pi m_2=1$.
This gives
\[
1=\frac1{\tau_{\rm p}^2}
+\frac{4\pi m_4k^2}{\tau_{\rm p}^4}
+\frac{4\pi m_6k^4}{\tau_{\rm p}^6}
+O(k^6).
\]
Solving this identity successively at orders \(k^2\) and \(k^4\)
gives \eqref{eq:taup-low-section2}.  Differentiating the large-\(s\)
series for \(\Acal\) and substituting it in
\eqref{eq:bromwich-residues-section2} gives
\eqref{eq:residue-low-section2}.  Since a smooth radial function
\(1+2\pi m_4|\xi|^2+O(|\xi|^4)\) has Hessian
\(4\pi m_4I_3\) at the origin,
\eqref{eq:taup-hessian-section2} follows.
\end{proof}

Away from this pole graph, the dispersion denominator stays
separated from zero.  The following bounds will control the regular
parts of the inverse transform.

\begin{lemma}
\label{lem:dispersion-separation-section2}
Fix \({\varepsilon_{\rm sep}}>0\).  There is \(c_{\varepsilon_{\rm sep}}>0\) such that
\begin{equation}
\label{eq:superterminal-separation-section2}
\inf_{\substack{k\geq k_*+{\varepsilon_{\rm sep}}\\ \tau\in\R}}
\bigl|1-\Lcal_+(\tau,k)\bigr|\geq c_{\varepsilon_{\rm sep}} .
\end{equation}
If \({\varepsilon_{\rm sep}}\leq k\leq k_*-{\varepsilon_{\rm sep}}\), the same lower bound holds after
removing the two sets
\[
|\tau-\tau_{\rm p}(k)|<{\varepsilon_{\rm sep}},\qquad
|\tau+\tau_{\rm p}(k)|<{\varepsilon_{\rm sep}},
\]
with \(c_{\varepsilon_{\rm sep}}\) reduced if necessary.  Moreover, on every compact
interval \(I\Subset(0,k_*)\), there are \(c_I,C_I,\varepsilon_I>0\)
such that, for \(k\in I\), \(\sigma=\pm\), and
\(\operatorname{Re}z\geq0\),
\begin{equation}
\label{eq:simple-pole-separation-section2}
c_I|z-\sigma\ri\tau_{\rm p}(k)|
\leq |1-\Lcal(z,k)|
\leq C_I|z-\sigma\ri\tau_{\rm p}(k)|
\end{equation}
whenever \(|z-\sigma\ri\tau_{\rm p}(k)|<\varepsilon_I\).
\end{lemma}

\begin{proof}
The classification in Proposition~\ref{prop:pole-branch-section2}
follows from \eqref{eq:Omega-boundary-section2}: on the open cut
away from zero $\Bcal\ne0$, the central gap has $\Acal<0$, and
under \eqref{H0} the center satisfies
\eqref{eq:filled-ball-zero-sign-appendix}.  For $s>1$, $\Acal$ is strictly decreasing;
the negative exterior interval is treated by evenness. The asserted bounds on compact
sets therefore follow by continuity.  They remain uniform as
\(|\tau|\to\infty\) because \(\Lcal_+(\tau,k)\to0\), and as
\(k\to\infty\) because \(\Acal\) and \(\Bcal\) are bounded and
\(\Lcal_+=k^{-2}\Om(\tau/k)\).  Finally,
\eqref{eq:simple-pole-separation-section2} is Taylor's theorem together
with \(\partial_z\Lcal(\sigma\ri\tau_{\rm p}(k),k)\ne0\).
\end{proof}

\subsection{Pole-branch geometry and outer-edge asymptotics}

To follow the pole as it approaches the outer edge of the continuum,
we use the ratio
\[
\rho:=\frac{k}{\tau}.
\]
This places the positive cut endpoint at \(\rho=1\) for every $k$,
so that pole and continuum boundary values can be described in the
same coordinate.  On the pole side $0<\rho<1$ and in a fixed
neighborhood of $\rho=1$, choose the same square-root branch
\begin{equation}
\label{eq:Omega0-definition-section2}
\Om_0(\rho)
:=\rho^{-1}\sqrt{\Om(\rho^{-1})},
\qquad \operatorname{Re}\Om_0(\rho)>0,
\end{equation}
continued from the positive real square root for \(0<\rho<1\), and
write
\[
\Om_0(\rho)=\omega(\rho)+\ri\nu(\rho).
\]
Then $\nu=0$ for $0<\rho\le1$ and $\nu>0$ for
\(\rho>1\) sufficiently close to \(1\).  On the pole side \(0<\rho<1\), define
\begin{equation}
\label{eq:q-V-section2}
q(\rho):=\rho\Om_0(\rho),\qquad
V(\rho):=
\frac{\omega'(\rho)}
{\omega(\rho)+\rho\omega'(\rho)}.
\end{equation}
We also define \(V\) for \(\rho>1\) sufficiently close to \(1\)
by the same real-valued quotient.
We parametrize the positive pole curve of
Proposition~\ref{prop:pole-branch-section2} by
$\rho=k/\tau_{\rm p}(k)\in(0,1)$.  The dispersion relation on this
curve gives
\begin{equation}
\label{eq:pole-parametrization-section2}
k=q(\rho),\qquad
\tau_{\rm p}(k)=\Om_0(\rho),\qquad
\tau_{\rm p}'(k)=V(\rho),\qquad 0<\rho<1.
\end{equation}
Proposition~\ref{prop:terminal-jet-section2}\textup{(i)} below gives
$\Om_0'>0$ on $(0,1)$, hence
$q'(\rho)=\Om_0(\rho)+\rho\Om_0'(\rho)>0$.
The chain rule along this curve yields
\begin{equation}\label{eq:pole-chain-rule-Green-revision}
\begin{aligned}
 \tau_{\rm p}'(q(\rho))&=\frac{\Om_0'(\rho)}{q'(\rho)}=V(\rho),\\
 \tau_{\rm p}''(q(\rho))&=\frac{V'(\rho)}{q'(\rho)}
 =\frac{\Om_0\Om_0''-2(\Om_0')^2}
        {(\Om_0+\rho\Om_0')^3}(\rho).
\end{aligned}
\end{equation}
The positivity of these derivatives is proved below.  The function
$q$ parametrizes the pole curve; the threshold spatial frequency is
still $k_*$.

The next proposition collects the geometry of the branch, its
endpoint expansions, and the group-velocity bounds used in the
threshold inversion.  Write \(h=\rho-1\),
\(h_+=\max\{h,0\}\), and \(h_-=\max\{-h,0\}\).

\begin{proposition}
\label{prop:terminal-jet-section2}
\textup{(i)}\enspace
The function \(\Om_0\) is positive, strictly increasing, and strictly
convex on \(0<\rho<1\), and \(V\) is positive and strictly increasing
there.  Consequently, \(\tau_{\rm p}\) is strictly increasing and
strictly convex on \(0<k<k_*\), while its phase velocity
\(\tau_{\rm p}(k)/k=\rho^{-1}\) decreases from \(+\infty\) to \(1\).
As \(\rho\to0\),
\begin{equation}
\label{eq:Omega0-origin-section2}
\Om_0(\rho)
=1+\varkappa_1\rho^2+\varkappa_2\rho^4+O(\rho^6),
\qquad
\varkappa_1=2\pi m_4,\qquad
\varkappa_2=2\pi m_6-2\pi^2m_4^2>0.
\end{equation}
Moreover,
\begin{equation}
\label{eq:terminal-frequency-section2}
\omega_*:=\Om_0(1)=\sqrt{4\pi\Acal(1)}=k_*.
\end{equation}

\par\medskip\noindent
\textup{(ii)}\enspace
If \({\gamma_1}>1\), then
\begin{equation}
\label{eq:Omega0-terminal-general-section2}
\Om_0(1+h)
=\omega_*+\omega_1h+\omega_2h^2+\Theta_{\gamma_1}(h)
+E_\epsilon(h),\qquad \omega_1>0,
\end{equation}
where, for each sufficiently small $\epsilon>0$,
\[
 |\partial_h^jE_\epsilon(h)|
 \le C_{\epsilon,j}|h|^{\min\{3,\gamma_1+1\}-\epsilon-j},
 \qquad j=0,1,2,
\]
separately on $h<0$ and $h>0$ near zero.  For
\({\gamma_1}\notin\mathbb N\),
\begin{equation}
\label{eq:Theta-noninteger-section2}
\Theta_{\gamma_1}(h)=c_-h_-^{\gamma_1}+c_+h_+^{\gamma_1},
\qquad
c_-=-\frac{A_{\gamma_1}}{\sin(\pi{\gamma_1})},\qquad
c_+=A_{\gamma_1}\bigl(-\cot(\pi{\gamma_1})+\ri\bigr)
=\re^{-\ri\pi{\gamma_1}}c_-,
\end{equation}
where
\begin{equation}
\label{eq:Agamma-section2}
A_{\gamma_1}
:=\frac{\pi^{3/2}\mu_{\rm out}(1)}
{2\sqrt{\Acal(1)}}>0.
\end{equation}
For \({\gamma_1}=m\in\mathbb N\), \(m\geq2\),
\begin{equation}
\label{eq:Theta-integer-section2}
\Theta_m(h)
=c_{\log}h^m\log|h|+c_-h_-^m+c_+h_+^m,
\qquad
c_{\log}<0,\quad c_-\in\R,\quad \operatorname{Im}c_+>0.
\end{equation}
In all cases \({\gamma_1}>1\),
\begin{equation}
\label{eq:nu-terminal-section2}
\nu(1+h)
=\bigl(\operatorname{Im}c_+\bigr)h^{\gamma_1}
+O_\epsilon(h^{{\gamma_1}+1-\epsilon}),
\qquad h\downarrow0.
\end{equation}

\par\medskip\noindent
\textup{(iii)}\enspace
For the borderline exponent $\gamma_1=1$, a real smooth function
$\Om_{0,\rm reg}$ and one-sided constants give
\begin{equation}\label{eq:gamma-one-terminal-section2}
 \Om_0(1+h)=\Om_{0,\rm reg}(1+h)
       +c_{\log}h\log|h|+c_+h_++c_-h_-+R_1(h),
\end{equation}
where $|\partial_h^jR_1(h)|\le C|h|^{2-j}(1+|\log|h||)^2$ for
$j=0,1,2$ on either open side.  In particular,
\begin{equation}
\label{eq:gamma-one-leading-section2}
c_{\log}
=-\frac{\sqrt\pi\,\mu_{\rm out}(1)}
{2\sqrt{\Acal(1)}}=:-a_0<0,
\qquad
\operatorname{Im}c_+
=\frac{\pi^{3/2}\mu_{\rm out}(1)}
{2\sqrt{\Acal(1)}}>0.
\end{equation}
More precisely, for \(y\downarrow0\),
\begin{align}
\Om_0'(1+\sigma y)
&=a_0\log(1/y)+d_\sigma
+O\!\left(y\log^2(1/y)\right),
\label{eq:gamma-one-first-derivative-section2}\\
\Om_0''(1+\sigma y)
&=-\sigma\frac{a_0}{y}
+O\!\left(\log^2(1/y)\right),
\qquad \sigma=\pm.
\label{eq:gamma-one-second-derivative-section2}
\end{align}

\par\medskip\noindent
\textup{(iv)}\enspace
The limiting pole group velocity at contact is therefore
\begin{equation}
\label{eq:vstar-section2}
v_*:=\lim_{k\uparrow k_*}\tau_{\rm p}'(k)
=
\begin{cases}
\displaystyle\frac{\omega_1}{\omega_*+\omega_1},
&{\gamma_1}>1,\\[2mm]
1,&{\gamma_1}=1.
\end{cases}
\end{equation}
For \({\gamma_1}>1\), \(q'(1)=\omega_*+\omega_1>0\), and the pole weight has
the useful form
\begin{equation}
\label{eq:residue-rho-section2}
r_{\rm p}(q(\rho))
=\frac{\Om_0(\rho)^2}
{2\bigl(\Om_0(\rho)+\rho\Om_0'(\rho)\bigr)}.
\end{equation}
It tends to
\(\omega_*^2/[2(\omega_*+\omega_1)]\) if \({\gamma_1}>1\), and to \(0\) if
\({\gamma_1}=1\).

For the group-velocity argument, set
\[
 I_{\rm osc}=
 \begin{cases}
 (1-{\varepsilon_{\rm tr}},1),&1\le{\gamma_1}\le2,\\
 (1-{\varepsilon_{\rm tr}},1)\cup(1,1+{\varepsilon_{\rm tr}}),&{\gamma_1}>2.
 \end{cases}
\]
On each component of $I_{\rm osc}$, $V$ is monotone and its derivative
is bounded away from zero in absolute value.  In particular,
\begin{equation}\label{eq:V-sublevel-section2}
 |\{\rho\in I_{\rm osc}:|V(\rho)-v|\le\eta\}|
 \le C\eta,\qquad 0<\eta\le1.
\end{equation}
When ${\gamma_1}=1$, more precisely,
\begin{align}
 V(1+\sigma y)&=1-\frac{\omega_*}{a_0\log(1/y)}
       +O(\log(1/y)^{-2}),\label{eq:V-gamma-one-section2}\\
 V'(1+\sigma y)&=-\sigma\frac{\omega_*}{a_0y\log^2(1/y)}
       \bigl(1+O(\log(1/y)^{-1})\bigr),
       \label{eq:Vprime-gamma-one-section2}
\end{align}
where $\sigma=-1$; only the pole side is used when $\gamma_1=1$.
\end{proposition}

\begin{proof}
\emph{Part (i): pole-side geometry.}
For \(0<\rho<1\), \(\Bcal(\rho^{-1})=0\), and a direct change of
variables gives the useful positive representation
\begin{equation}
\label{eq:Omega0-positive-representation-section2}
\Om_0(\rho)^2
=\rho^{-2}\Om(\rho^{-1})
=4\pi\rho^{-2}\Acal(\rho^{-1})
=4\pi\int_\kappa^1
\frac{u^2\mu(u^2/2)}{1-\rho^2u^2}\,\rd u .
\end{equation}
In terms of \(r>1\), differentiation under the integral gives
\begin{gather}
r\Acal'(r)+2\Acal(r)
=-2\int_\kappa^1
\frac{u^4\mu(u^2/2)}{(r^2-u^2)^2}\,\rd u<0,
\label{eq:A-rigid-first-section2}\\
\Acal(r)\Acal''(r)-\frac32\Acal'(r)^2>0.
\label{eq:A-rigid-second-section2}
\end{gather}
{
For completeness, put
\[
 I_j(r):=\int_\kappa^1
 \frac{u^2\mu(u^2/2)}{(r^2-u^2)^j}\,\rd u,\qquad j=1,2,3.
\]
Then $\Acal=I_1$, $\Acal'=-2rI_2$, and
$\Acal''=-2I_2+8r^2I_3$.  Cauchy--Schwarz gives
$I_2^2\le I_1I_3$, and consequently
\begin{equation}\label{eq:A-rigid-positive-lower-section2}
 \Acal\Acal''-\tfrac32(\Acal')^2
 =-2I_1I_2+8r^2I_1I_3-6r^2I_2^2
 \ge 2I_1\int_\kappa^1
       \frac{u^4\mu(u^2/2)}{(r^2-u^2)^3}\,\rd u>0.
\end{equation}
The strict lower bound follows by restricting the positive integral
to a fixed compact subinterval of $(\kappa,1)$, also with $\kappa=0$
under \eqref{H0}.
}  With \(r=\rho^{-1}\),
\eqref{eq:A-rigid-first-section2} implies \(\Om_0'(\rho)>0\).
Furthermore,
\begin{equation}
\label{eq:Vprime-identity-section2}
V'(\rho)
=\frac{\Om_0(\rho)\Om_0''(\rho)-2\Om_0'(\rho)^2}
{\bigl(\Om_0(\rho)+\rho\Om_0'(\rho)\bigr)^2},
\qquad 0<\rho<1,
\end{equation}
and the numerator is positive by
\eqref{eq:A-rigid-second-section2}.  More explicitly,
\[
2\Om_0^2\bigl(\Om_0\Om_0''-2(\Om_0')^2\bigr)
=\frac{(4\pi)^2}{\rho^8}
\left(\Acal(r)\Acal''(r)-\frac32\Acal'(r)^2\right)>0.
\]
This proves $\Om_0''>0$ and $V'>0$.  Since $q'>0$,
\eqref{eq:pole-chain-rule-Green-revision} now gives both
$\tau_{\rm p}'>0$ and $\tau_{\rm p}''>0$.

Expanding \eqref{eq:Omega0-positive-representation-section2} at
\(\rho=0\) gives
\[
\Om_0(\rho)
=\left(1+4\pi m_4\rho^2+4\pi m_6\rho^4+O(\rho^6)\right)^{1/2},
\]
which is \eqref{eq:Omega0-origin-section2}.  The positivity of
\(\varkappa_2\) follows from Cauchy--Schwarz and \(4\pi m_2=1\).

\emph{Parts (ii)--(iii): transfer of the endpoint expansion.}
Near \(s=1\), write
\[
\Acal=\Acal_{\rm reg}+\Acal_{\rm sig},
\qquad
\Om_{\rm reg}:=4\pi\Acal_{\rm reg},
\qquad
\Om_{\rm sig}:=4\pi\Acal_{\rm sig}+2\pi\ri\Bcal ,
\]
where \(\Acal_{\rm reg}\) is real and smooth and
\(\Acal_{\rm sig}\) is the explicit term in
\eqref{eq:A-outer-noninteger-section2} or
\eqref{eq:A-outer-integer-section2}.  Set
\[
Z_{\rm reg}(\rho):=\Om_{\rm reg}(\rho^{-1}),
\qquad
Z_{\rm sig}(\rho):=\Om_{\rm sig}(\rho^{-1}).
\]
After shrinking the neighborhood, \(Z_{\rm reg}\) remains positive
and
\[
\Om_0(\rho)
=\rho^{-1}\sqrt{Z_{\rm reg}(\rho)}
+\rho^{-1}
\left(\sqrt{Z_{\rm reg}(\rho)+Z_{\rm sig}(\rho)}
-\sqrt{Z_{\rm reg}(\rho)}\right).
\]
If \({\gamma_1}>1\), Taylor's formula for the square root yields
\begin{equation}
\label{eq:square-root-first-section2}
\Om_0(\rho)
=\Om_{0,\rm reg}(\rho)
+\frac{Z_{\rm sig}(\rho)}
{2\rho\sqrt{Z_{\rm reg}(\rho)}}
+O_\epsilon(|\rho-1|^{2{\gamma_1}-\epsilon}).
\end{equation}
Since \(2{\gamma_1}>{\gamma_1}+1\), the last term is absorbed in the
remainder of \eqref{eq:Omega0-terminal-general-section2}; freezing the
smooth coefficient in the middle term costs
\(O_\epsilon(|h|^{{\gamma_1}+1-\epsilon})\).
The conormal product and reciprocal rules preserve these bounds after
$j\le2$ derivatives, each derivative costing one power of $|h|$;
this gives the stated derivative bounds for $E_\epsilon$.

For noninteger \({\gamma_1}\), the singular term on the pole side
\(h<0\) is
\[
-\frac{2\pi^2\mu_{\rm out}(1)}{\sin(\pi{\gamma_1})}h_-^{\gamma_1},
\]
before division by \(2\sqrt{4\pi\Acal(1)}\).  On the continuum side
\(h>0\), the real Hilbert-transform term and
\eqref{eq:B-outer-section2} combine as
\[
2\pi^2\mu_{\rm out}(1)
\bigl(-\cot(\pi{\gamma_1})+\ri\bigr)h_+^{\gamma_1} .
\]
This proves \eqref{eq:Theta-noninteger-section2}, including
\(c_+=\re^{-\ri\pi{\gamma_1}}c_-\).  At an integer exponent, the same
calculation with the logarithmic boundary value gives
\eqref{eq:Theta-integer-section2}.  The imaginary part in
\eqref{eq:B-outer-section2} also proves
\eqref{eq:nu-terminal-section2}.

When $\gamma_1=1$, $Z_{\rm sig}=O_{\rm con}(h\log|h|)$.
Here $h=\rho-1$ and the conormal bounds are required through order
two in the conventions at the end of Section~\ref{sec:introduction}; explicitly,
$|(h\partial_h)^jZ_{\rm sig}|\le C_j|h|(1+|\log|h||)$ for $j\le2$
on each open side.
The first-order Taylor formula for the square root has a remainder
$O_{\rm con}(h^2(1+|\log|h||)^2)$ through two derivatives.
It gives \eqref{eq:gamma-one-terminal-section2}; differentiation yields
\eqref{eq:gamma-one-first-derivative-section2}--
\eqref{eq:gamma-one-second-derivative-section2}.

\emph{Part (iv): limiting velocity at contact and the four regimes.}
For \({\gamma_1}>1\), the pole-side derivative has the finite value
\[
\omega_1
=-\frac{\sqrt{4\pi}}{2\sqrt{\Acal(1)}}
\bigl(2\Acal(1)+\Acal'(1)\bigr)>0,
\]
where the sign follows from \eqref{eq:A-rigid-first-section2}.
Taking $\rho\uparrow1$ in
$\tau_{\rm p}'(q(\rho))=\Om_0'(\rho)/(\Om_0+\rho\Om_0')$,
as given by \eqref{eq:pole-chain-rule-Green-revision}, proves
\eqref{eq:vstar-section2}.
  Combining the same change of variables
with \eqref{eq:bromwich-residues-section2} gives
\eqref{eq:residue-rho-section2}.  When \({\gamma_1}=1\),
\(\Om_0'(\rho)\to+\infty\) on the pole side, which gives
\(v_*=1\) and the vanishing limiting pole residue.

For any twice differentiable real branch \(\omega\), direct
differentiation gives
\[
V'(\rho)
=\frac{\omega(\rho)\omega''(\rho)-2\omega'(\rho)^2}
{\bigl(\omega(\rho)+\rho\omega'(\rho)\bigr)^2}.
\]
We finally verify the one-sided monotonicity required in
\eqref{eq:V-sublevel-section2}.  For \({\gamma_1}>1\), the denominator tends to \((\omega_*+\omega_1)^2>0\) as \(\rho\to1\). Suppose first that \(1<{\gamma_1}<2\), and write
\(\rho=1-y\), \(y\downarrow0\), on the pole side.  From
\eqref{eq:Theta-noninteger-section2}, 	
	\[
	\omega''(1-y)
	=
	{\gamma_1}({\gamma_1}-1)c_-y^{{\gamma_1}-2}
	+o(y^{{\gamma_1}-2}),
	\qquad
	c_-=-\frac{A_{\gamma_1}}{\sin(\pi{\gamma_1})}>0.
	\]
	Since \(\omega'\) remains bounded,
	\[
	V'(1-y)
	=
	\frac{{\gamma_1}({\gamma_1}-1)\omega_*c_-}
	{(\omega_*+\omega_1)^2}
	y^{{\gamma_1}-2}(1+o(1))>0.
	\]
	For \({\gamma_1}=2\), the pole-side logarithmic term in
	\eqref{eq:Theta-integer-section2} gives
	\[
	\omega''(1-y)=2c_{\log}\log y+O(1),
	\]
	and therefore
	\[
	V'(1-y)
	=
	\frac{2\omega_*c_{\log}}
	{(\omega_*+\omega_1)^2}\log y+O(1)>0
	\]
	for \(y\) sufficiently small, since \(c_{\log}<0\).
	Thus, for \(1<{\gamma_1}\le2\), \(V\) is strictly increasing on a
	sufficiently small pole-side interval and \(|V'|\) is bounded below
	there.  No corresponding continuum-side estimate is needed in this
	range, since the continuum center is treated directly by the damping
	\(\nu(1+h)\asymp h^{\gamma_1}\).

If $\gamma_1>2$, the endpoint expansion is $C^2$ across $\rho=1$:
$\partial_h^2\Theta_{\gamma_1}(h)=o(1)$ in
\eqref{eq:Theta-noninteger-section2} or \eqref{eq:Theta-integer-section2},
and the explicit remainder bound after
\eqref{eq:Omega0-terminal-general-section2} gives
$\partial_h^2E_\epsilon(h)=o(1)$ for sufficiently small $\epsilon>0$.
{
 On the pole side, the integrands defining $I_1,I_2,I_3$ at $r=1$
 behave at worst like $(1-u)^{{\gamma_1}-3}$, which is integrable because
 ${\gamma_1}>2$.  Dominated convergence in
 \eqref{eq:A-rigid-positive-lower-section2} therefore gives the
 strictly positive endpoint lower bound
 \[
 \bigl(\Acal\Acal''-\tfrac32(\Acal')^2\bigr)(1)
 \ge 2\Acal(1)\int_\kappa^1
       \frac{u^4\mu(u^2/2)}{(1-u^2)^3}\,\rd u>0.
 \]
 Using \eqref{eq:Vprime-identity-section2} and the displayed identity
 relating its numerator to $\Acal$ now yields
 \[
 \lim_{\rho\uparrow1}V'(\rho)
 =\frac{2(\omega_*\omega_2-\omega_1^2)}
       {(\omega_*+\omega_1)^2}=:V_*'>0.
 \]
}
	By \eqref{eq:Omega0-terminal-general-section2}, the continuum limits
of $\omega$, $\omega'$, and $\omega''$ are likewise
$\omega_*$, $\omega_1$, and $2\omega_2$.
Hence $V'(\rho)\to V_*'>0$ also as $\rho\downarrow1$.
After shrinking the neighborhood of $\rho=1$, $V'\ge c>0$ on both sides.
	
	Finally, let \({\gamma_1}=1\), put \(L=\log(1/y)\), and use
	\eqref{eq:gamma-one-first-derivative-section2}--
	\eqref{eq:gamma-one-second-derivative-section2}.  For
	\(\sigma=-1\),
	\[
	\omega'(1+\sigma y)=a_0L+O(1),
	\qquad
	\omega''(1+\sigma y)
	=-\sigma\frac{a_0}{y}+O(L^2),
	\qquad
	\omega(1+\sigma y)=\omega_*+O(yL).
	\]
	Substitution in the formulas for \(V\) and \(V'\) gives
	\begin{gather*}
		V(1+\sigma y)=
		1-\frac{\omega_*}{a_0L}
		+O(L^{-2}),\\
		V'(1+\sigma y)=
		-\sigma\frac{\omega_*}{a_0yL^2}
		\bigl(1+O(L^{-1})\bigr),
	\end{gather*}
	which are \eqref{eq:V-gamma-one-section2}--
	\eqref{eq:Vprime-gamma-one-section2}.  Hence \(V\) is monotone on
	the pole side and \(|V'|\) is bounded below after
	shrinking the chart.
	
	We have therefore obtained, on every connected component $I_{\rm osc}^{(j)}$ of
	$I_{\rm osc}$, $|V'(\rho)|\ge c_0>0.$
	Since $V$ is monotone on $I_{\rm osc}^{(j)}$, the mean-value theorem gives
	\[
	\bigl|
	\{\rho\in I_{\rm osc}^{(j)}:|V(\rho)-v|\le\eta\}
	\bigr|
	\le \frac{2\eta}{c_0}.
	\]
	Summing over the at most two components of \(I_{\rm osc}\) proves
	\eqref{eq:V-sublevel-section2}.
\end{proof}

The boundary factorization through $\Om_0$ identifies the pole
prescription and the exact continuum damping coefficient.

\begin{corollary}
\label{cor:terminal-factorization-section2}
Let $\tau>0$, $k=\tau\rho$, and $\rho$ be sufficiently close to $1$.
The positive-frequency boundary value $Q_+(\tau,k)$ in
\eqref{eq:resolvent-boundary-section2} factorizes as
\begin{equation}
\label{eq:terminal-factorization-section2}
Q_+(\tau,\tau\rho)
=\frac{\tau^2}{\tau+\Om_0(\rho)}
\frac1{\tau-\omega(\rho)-\ri(\nu(\rho)+0)}.
\end{equation}
Here $+0$ is the boundary prescription inherited from $\Re z>0$;
it is relevant only when $\nu=0$.  On the continuum side,
\begin{equation}\label{eq:damping-exact-factor-section2}
 \nu_{\rm c}(h):=\nu(1+h)
 =\frac{\pi^2\mu_{\rm out}((1+h)^{-1})}
        {(1+h)^{\gamma_1+3}\omega(1+h)}\,h^{\gamma_1},
 \qquad 0<h\ll1.
\end{equation}
In particular, $\nu_{\rm c}(h)\asymp h^{\gamma_1}$.
\end{corollary}

\begin{proof}
By \eqref{eq:Omega0-definition-section2},
$\Om(\rho^{-1})=\rho^2\Om_0(\rho)^2$.  Hence
\[
Q_+(\tau,\tau\rho)
=\frac{\tau^2\rho^2}
{\tau^2\rho^2-\Om(\rho^{-1})}
=\frac{\tau^2}
{(\tau-\Om_0(\rho))(\tau+\Om_0(\rho))},
\]
which gives \eqref{eq:terminal-factorization-section2} with the
orientation fixed in \eqref{eq:bromwich-model-section2}.
Taking imaginary parts of
$\Om_0(\rho)^2=\rho^{-2}\Om(\rho^{-1})$ for $\rho=1+h>1$ gives
\[
 \nu_{\rm c}(h)=\frac{\pi\rho^{-2}\Bcal(\rho^{-1})}{\omega(\rho)}.
\]
Substitution of \eqref{eq:B-outer-section2} proves
\eqref{eq:damping-exact-factor-section2}; its coefficient tends to
the positive constant $\pi^2\mu_{\rm out}(1)/\omega_*$.
\end{proof}

\section{Spectral representation and the wave--damping decomposition}
\label{sec:spectral-decomposition}

This section gives the spectral representation of the response and
its wave--damping decomposition.  We first derive the pole--cut
formula away from the survival threshold.  We then construct the
spectral measure across contact and recover the resolvent and the time
response from it.  After fixing the spatial-frequency decomposition,
we collect the frequency and endpoint estimates for the regular
response.

\subsection{Spectral representation away from the threshold}

Put
\begin{equation}\label{def:I_kappa}
	\mathcal I_\kappa=[-1,-\kappa]\cup[\kappa,1].
\end{equation}
The jump in \eqref{eq:jump-density-section2} is understood to be zero
outside $k\mathcal I_\kappa$.  The following proposition represents the
causal response and its Laplace resolvent with the same pole weights
and cut density.  This common representation fixes the normalization
used to continue through contact and to invert in space.

\begin{proposition}
\label{prop:fixed-k-representation-section3}
Recall that $\Sigma_k=k\mathcal I_\kappa=[-k,-\kappa k]\cup[\kappa k,k]$,
which is $[-k,k]$ when $\kappa=0$.
For $k>0$, $k\ne k_*$, and $t>0$,
\begin{equation}
\label{eq:fixed-k-representation-section3}
\widehat G(t,k)
=\one_{\{k<k_*\}}
\sum_{\sigma=\pm}\varrho_\sigma(k)
\re^{\ri\sigma t\tau_{\rm p}(k)}
+\int_{\Sigma_k}\re^{\ri t\tau}J_k(\tau)\,\rd\tau .
\end{equation}
The same pole weights and cut density reconstruct the resolvent:
\begin{equation}
\label{eq:cauchy-slit-appendix}
R(z,k)
=\one_{\{k<k_*\}}\sum_{\sigma=\pm}
\frac{\varrho_\sigma(k)}{z-\ri\sigma\tau_{\rm p}(k)}
+\int_{\Sigma_k}\frac{J_k(\tau)}{z-\ri\tau}\,\rd\tau,
\qquad \Re z>0.
\end{equation}
At $k=0$ the continuous definition is
\begin{equation}
\label{eq:fixed-k-zero-section3}
\widehat G(t,0)=-\sin t,\qquad t>0.
\end{equation}
The response is extended by zero to $t<0$.
\end{proposition}

\begin{proof}
Fix $k>0$, $k\ne k_*$.  We deform the Bromwich contour in
\eqref{eq:Green-fun-G-t-k} for the original resolvent $R(\cdot,k)$.
The Cauchy formula \eqref{eq:analytic-cauchy-section2} extends
$\Lcal(\cdot,k)$ analytically to $\C\setminus\ri\Sigma_k$ and gives
\[
 \Lcal(-\overline z,k)=\overline{\Lcal(z,k)}.
\]
A zero of $1-\Lcal$ in the open left half-plane would reflect to
one in the open right half-plane, excluded by
Lemma~\ref{lem:right-half-plane-stability-section2}.
Proposition~\ref{prop:pole-branch-section2} classifies the remaining
zeros of $1-\Lcal$, which are the poles of $R$: for $0<k<k_*$
there are exactly the two simple poles
$z=\pm\ri\tau_{\rm p}(k)$ outside the slits, and for $k>k_*$ there
are no imaginary-axis poles.

At each normalized endpoint $a=\pm1$, and also $a=\pm\kappa$
when $\kappa>0$, the denominator $k^2-\Om(s)$ has a nonzero limit
as $s\to a$ along the cut.  This limit is $k^2-k_*^2$ at the outer
endpoints and $k^2-4\pi\Acal(\kappa)>0$ at the inner endpoints,
since $\Acal(\kappa)<0$ by \eqref{eq:A-B-section2}.
Lemma~\ref{lem:A-endpoint-section2} and the endpoint factor in
\eqref{H2} or \eqref{H0} give
\[
 J_k(ks)=O_k(|s-a|^{\gamma_{|a|}})
\]
through \eqref{eq:jump-density-section2}; the negative endpoints
follow by parity.  In particular, the jump is integrable.

For the endpoint circles in the contour deformation, we need
$|R(z,k)|\le C_k$ in complex neighborhoods of the slit endpoints.
The denominator limits just computed are nonzero for this fixed
$k>0$, $k\ne k_*$.  It therefore suffices to establish the corresponding
finite complex limits of $\Lcal$.  By the oddness of $\mu_0'$, the Cauchy representation
\eqref{eq:analytic-cauchy-section2} also reads
\[
 \Lcal(z,k)=\frac1{k^2}\int_{\R}
                   \frac{\mu_0'(s)}{s+\ri z/k}\,\rd s.
\]
Near $z=\ri ka$, choose a smooth cutoff $\chi_a$, equal to one
near $s=a$ and supported in a small neighborhood of that endpoint.
The complementary integral is
\[
 \frac1{k^2}\int_{\R}
 \frac{(1-\chi_a(s))\mu_0'(s)}{s+\ri z/k}\,\rd s.
\]
Its denominator is bounded away from zero on the numerator's support
for $z$ near $\ri ka$, so it is holomorphic in a full neighborhood
of that endpoint.  An affine change of variable reduces the localized
nonsmooth term to
\[
 I(w)=\int_0^\delta\frac{f(v)}{v-w}\,\rd v,
 \qquad f(v)=v^{\gamma_{|a|}}b(v),\qquad
 w\notin[0,\delta],\quad w\longrightarrow0,
\]
where $b$ is smooth and vanishes near $\delta$.
Write $w=x+\ri y$.  If $0<x<\delta/2$, then
\[
 I(w)=\int_0^\delta\frac{f(v)-f(x)}{v-w}\,\rd v
       +f(x)\bigl[\log(\delta-w)-\log(-w)\bigr].
\]
The first integrand is bounded by the Lipschitz constant of $f$,
since $|v-w|\ge|v-x|$; dominated convergence sends its integral to
$\int_0^\delta f(v)/v\,\rd v$.  The second term tends to zero,
being bounded by $C|w|^{\gamma_{|a|}}(1+|\log|w||)$.
For $x\le0$, the inequality $|v-w|\ge v$ gives the same limit by
direct dominated convergence.  Thus $\Lcal$ has a finite limit from
both sides of the slit at each endpoint.  The nonzero denominator
limits give $|R(z,k)|\le C_k$ in the required neighborhoods, with
constants allowed to depend on the fixed distance of $k$ from $k_*$.

On each open cut interval, the boundary formulas
\eqref{eq:Omega-boundary-section2} are locally uniform and smooth.
Their resolvent denominators are nonzero because $\Bcal\ne0$ away
from zero and, under \eqref{H0}, $\Acal(0)<0$ by
\eqref{eq:filled-ball-zero-sign-appendix}.  Thus the same local
uniformity and smoothness hold for $R_\pm$.
This local uniformity rules out any additional singular measure
on an open lip.  At infinity, the Cauchy representation on the slit
plane gives the same bound as \eqref{eq:dispersion-decay-section2},
also in the left half-plane:
\[
 |\Lcal(z,k)|\le (|z|-k)^{-2}\le4|z|^{-2},
 \qquad |R(z,k)|\le8|z|^{-2}
\]
for sufficiently large $|z|>2k$.

Orient the Bromwich line upwards, truncate it at
$\sigma\pm\ri T$, and close it by the left semicircle centered at
$\sigma$, as in Figure~\ref{fig:bromwich}.  For fixed $t>0$, the
outer arc contributes at most $C \re^{\sigma t}T^{-1}$.
\begin{figure}[H]
\centering
\begin{tikzpicture}[x=1cm,y=.8cm,>=Stealth,font=\small]
\draw[->,gray] (-3.7,0)--(2.15,0) node[below left] {$\Re z$};
\draw[->,gray] (0,-3.7)--(0,3.9) node[above] {$\Im z$};
\draw[thick,->] (1.5,-3.4)--(1.5,3.4) node[pos=.87,right] {$\sigma$};
\draw[->] (1.5,3.4) arc[start angle=90,end angle=180,radius=3.4];
\draw[->] (-1.9,0) arc[start angle=180,end angle=270,radius=3.4];
\node[left] at (-1.9,-.6) {$|z-\sigma|=T$};
\draw[line width=2.2pt,gray] (0,.8)--(0,2.0);
\draw[line width=2.2pt,gray] (0,-.8)--(0,-2.0);
\draw[->] (.22,1.92)--(.22,.86) node[pos=.45,right] {$R_+$};
\draw[->] (-.22,.86)--(-.22,1.92) node[pos=.45,left] {$R_-$};
\draw[->] (.22,-.86)--(.22,-1.92);
\draw[->] (-.22,-1.92)--(-.22,-.86);
\draw[dashed] (0,2.0) ellipse (.22 and .18);
\draw[dashed] (0,.8) ellipse (.22 and .18);
\draw[dashed] (0,-2.0) ellipse (.22 and .18);
\draw[dashed] (0,-.8) ellipse (.22 and .18);
\fill (0,2.65) circle (2pt) node[right=5pt] {$\ri\tau_{\rm p}$};
\fill (0,-2.65) circle (2pt) node[right=5pt] {$-\ri\tau_{\rm p}$};
\draw[->] (.29,2.65) arc[start angle=0,end angle=-300,radius=.29];
\draw[->] (.29,-2.65) arc[start angle=0,end angle=-300,radius=.29];
\node[left] at (-.28,2.0) {$\ri k$};
\node[left] at (-.28,.8) {$\ri\kappa k$};
\node[anchor=west,text width=5.5cm,align=left] at (2.7,1.15)
{Outer boundary: counterclockwise.\\[3pt]
Slit lips and pole circles: clockwise.\\[3pt]
The right lip runs downward as an inner boundary; moving it to the other side gives the jump $R_+-R_-$.};
\node[anchor=west,text width=5.5cm,align=left] at (2.7,-1.65)
{Fixed $0<k<k_*$ is shown.\\[3pt]
At $k=k_*$ the pole meets the cut: use the joint spectral measure, not a uniform small-circle estimate.};
\end{tikzpicture}
\caption{The Bromwich deformation for $R(\cdot,k)$ at a fixed $0<k<k_*$ for the shell case $\kappa>0$, schematically. For $k>k_*$ the pole circles are absent; contact $k=k_*$ is treated by the joint spectral measure below. The line $\Re z=\sigma$ is oriented upward. Clockwise inner boundaries give the residue terms and $(R_+-R_-)\,\rd\tau/(2\pi)$ after transposition. At fixed outer radius, the lip width first tends to zero, then the noncritical endpoint radii; the outer radius is sent to infinity last.}
\label{fig:bromwich}
\end{figure}
Delete disks of radius $d$ about the slit endpoints, thin
neighborhoods of width $\eta\ll d$ about the remaining slits, and
small disks about the two poles when $k<k_*$.
At fixed $d,T$, first let $\eta\downarrow0$ using the locally
uniform boundary traces.  Next let $d\downarrow0$; the endpoint
arcs and omitted lip pieces are $O_{k,t}(d)$ by the complex
endpoint bound.  The pole circles tend to their residues.  Finally
let $T\to\infty$ using the outer-arc bound above.

The inner boundaries are clockwise: their right lips run downwards
and their left lips upwards.  Transposing those integrals reverses
the orientations, so $\rd z=\ri\,\rd\tau$ gives
\begin{equation}
\label{eq:bromwich-slit-appendix}
\begin{split}
\frac1{2\pi\ri}\int_{\sigma-\ri\infty}^{\sigma+\ri\infty}
 \re^{zt}R(z,k)\,\rd z
={}&\one_{\{k<k_*\}}\sum_{\sigma'=\pm}
 \re^{\ri\sigma't\tau_{\rm p}(k)}
 \Res_{z=\ri\sigma'\tau_{\rm p}(k)}R(z,k)\\
&+\frac1{2\pi}\int_{\Sigma_k}\re^{\ri t\tau}
 (R_+(\tau,k)-R_-(\tau,k))\,\rd\tau .
\end{split}
\end{equation}
The clockwise pole circles likewise give positive residue terms
after transposition.  Substituting
\eqref{eq:jump-density-section2} and
\eqref{eq:bromwich-residues-section2} proves
\eqref{eq:fixed-k-representation-section3}.  The factor
$1/(2\pi\ri)$ already divides out the contour normalization, so
an atom has weight exactly $\varrho_\sigma$, with no extra factor.
For $t<0$ close to the right; there are no enclosed singularities
and the same outer-arc estimate applies.

For fixed $\Re z>0$, apply the same slit contour to
$R(w,k)/(w-z)$ on a disk containing $z$, the poles and the cuts.
Its outer integral is $O(T^{-2})$, and its endpoint integrals are
$O_{k,z}(d)$.  The residue at $w=z$, the pole fractions and the
oriented lips yield \eqref{eq:cauchy-slit-appendix}.
Finally, \eqref{eq:kzero-resolvent} gives
\eqref{eq:fixed-k-zero-section3}. The identities $\varrho_-=\overline{\varrho_+}$ and
$J_k(-\tau)=\overline{J_k(\tau)}$ show that the sum of the two pole
terms and the sum of the two cut contributions are each real.
\end{proof}
We collect the pole atoms and the cut density in the finite complex
spectral measure.  For $k>0$, $k\ne k_*$, define
\begin{equation}
	\label{eq:physical-measure-nonterminal-section3}
	\rd\mathfrak M_k(\tau)
	:=
	\one_{\{k<k_*\}}
	\sum_{\sigma=\pm}
	\varrho_\sigma(k)\,
	\delta_{\sigma\tau_{\rm p}(k)}(\rd\tau)
	+
	\one_{\Sigma_k}(\tau)J_k(\tau)\,\rd\tau .
\end{equation}
The next subsection extends this same measure through $k_*$.

\subsection{Spectral measure across the pole--continuum threshold}

To continue the spectral measure through $k_*$, we first isolate
its contributions near the two outer spectral edges.  For $\gamma_1>1$, a continuum
Lorentzian for $k>k_*$ converges to the same atom as
the pole for $k<k_*$.  Their joint limit therefore
determines the contribution at contact.

Near the positive outer endpoint, use $s=\tau/k$ and set
\begin{equation}
\label{eq:terminal-XY-section3}
a=k^2-k_*^2,\qquad y=1-s,\qquad
X_k(y)=k^2-4\pi\Acal(1-y),\qquad
Y(y)=2\pi\Bcal(1-y).
\end{equation}
Choose $\varepsilon_{\rm tr}>0$ so small that this chart meets
only the outer endpoint and $\partial_yX_k<0$ on
$0<|y|<{\varepsilon_{\rm tr}}$ for $k$ near $k_*$.
The continuous function $X_k$ is strictly decreasing across $y=0$
as well.  If $k<k_*$, put
\[
y_{\rm p}(k)=1-\frac{\tau_{\rm p}(k)}k<0.
\]
Fix a sufficiently small $k$-neighborhood
$|k-k_*|<\eta_{\rm tr}$, with $0<\eta_{\rm tr}\le\varepsilon_{\rm tr}$,
so that $k>0$, $X_k(\varepsilon_{\rm tr})<0$, and $y_{\rm p}(k)$
lies strictly inside the chart whenever $k<k_*$.
The sign at $\varepsilon_{\rm tr}$ also keeps the supercritical
zero inside the chart when $k>k_*$.  The local measures below are defined on
$[-\varepsilon_{\rm tr},\varepsilon_{\rm tr}]$ and extended by zero
outside it.  Their densities give no mass to the artificial
boundaries $y=\pm\varepsilon_{\rm tr}$, and no atoms are placed there.
This makes explicit the restriction to the outer-edge chart; the
possible critical atom at $y=0$ is retained.  Under
$T_k(y)=k(1-y)$, the local measure $\mathfrak m_k^+$ below is the
coordinate representation of the positive edge slice of
\eqref{eq:physical-measure-nonterminal-section3}.
Indeed, $(T_k)_\#(\ri r_{\rm p}\delta_{y_{\rm p}})
=\ri r_{\rm p}\delta_{\tau_{\rm p}}$, so the atom keeps its mass,
whereas the continuous density becomes
$kJ_k(k(1-y))\,\rd y$ because $|T_k'|=k$.
The negative edge slice and the rest of the cut must still be
added to recover the full spectral measure.

\begin{proposition}
\label{prop:terminal-measure-section3}
For $0<|k-k_*|<\eta_{\rm tr}$, define on
$[-\varepsilon_{\rm tr},\varepsilon_{\rm tr}]$
\begin{equation}
\label{eq:terminal-measure-section3}
\begin{split}
\rd\mathfrak m_k^+(y)
:={}\ri\one_{\{k<k_*\}}r_{\rm p}(k)
\delta_{y_{\rm p}(k)}(\rd y)+\one_{\{y>0\}}
\frac{\ri k^3Y(y)}
{\pi\bigl(X_k(y)^2+Y(y)^2\bigr)}\,\rd y .
\end{split}
\end{equation}
For $0<y\le\varepsilon_{\rm tr}$, the second term is the restriction
of the positive-cut measure from \eqref{eq:jump-density-section2}
to this chart:
\begin{equation}
\label{eq:terminal-cut-weight-section3}
kJ_k(k(1-y))\,\rd y
=\frac{\ri k^3Y(y)}
{\pi\bigl(X_k(y)^2+Y(y)^2\bigr)}\,\rd y .
\end{equation}
The total variations are bounded uniformly for $k$ in this
neighborhood.  As $k$ approaches $k_*$ from either side, the measures
converge weakly, against continuous test functions on
$[-\varepsilon_{\rm tr},\varepsilon_{\rm tr}]$, to the same limiting
measure $\mathfrak m_{k_*}^+$ described below.

If ${\gamma_1}>1$, set
\begin{equation}
\label{eq:terminal-constants-section3}
c_*=-4\pi\Acal'(1)>0,\qquad
d_*=2\pi^2\mu_{\rm out}(1)>0,\qquad
r_*:=\frac{k_*^3}{c_*}.
\end{equation}
Then
\begin{equation}
\label{eq:terminal-critical-measure-greater-section3}
\rd\mathfrak m_{k_*}^+(y)
=\ri r_*\delta_0(\rd y)
+\one_{\{y>0\}}
\frac{\ri k_*^3Y(y)}
{\pi\bigl(X_{k_*}(y)^2+Y(y)^2\bigr)}\,\rd y .
\end{equation}
The density is $O(y^{{\gamma_1}-2})$.  If ${\gamma_1}=1$, there is no
atom at contact and
\begin{equation}
\label{eq:terminal-critical-measure-one-section3}
\rd\mathfrak m_{k_*}^+(y)
=\one_{\{y>0\}}
\frac{\ri k_*^3Y(y)}
{\pi\bigl(X_{k_*}(y)^2+Y(y)^2\bigr)}\,\rd y,
\qquad
\left|\frac{k_*^3Y(y)}
{X_{k_*}(y)^2+Y(y)^2}\right|
\lesssim\frac1{y\log^2(1/y)}.
\end{equation}
In particular, $\mathfrak m_{k_*}^+$ is finite in all cases
${\gamma_1}\ge1$.
\end{proposition}

\begin{proof}
Identity \eqref{eq:terminal-cut-weight-section3} follows from
\eqref{eq:jump-density-section2} after the change $\tau=k(1-y)$, reversing integration limits so
that the measure Jacobian is $|\rd\tau|=k\,\rd y$.  Since
$|\partial_yX_k(y_{\rm p})|=-4\pi\Acal'(1-y_{\rm p})$,
\eqref{eq:bromwich-residues-section2} also gives the following
form of the pole mass, which puts it in the same denominator
variable $X_k$ as the Lorentzian density:
\begin{equation}
\label{eq:terminal-pole-delta-section3}
\ri r_{\rm p}(k)\delta_{y_{\rm p}(k)}
=\ri k^3\delta(X_k(y))\one_{\{y<0\}}.
\end{equation}

For the continuum variable $y\downarrow0$, substitute the outer
endpoint expansions
\eqref{eq:A-outer-noninteger-section2}--\eqref{eq:A-outer-integer-section2}
and \eqref{eq:B-outer-section2} into
\eqref{eq:terminal-XY-section3}.  With the constants in
\eqref{eq:terminal-constants-section3}, Lemma~\ref{lem:A-endpoint-section2}
gives, uniformly with one conormal derivative,
\begin{align}
X_k(y)&=a-c_*y
+O\bigl(y^2+y^{\gamma_1}(1+|\log y|)\bigr),
&Y(y)&=d_*y^{\gamma_1}(1+O(y)),
&&{\gamma_1}>1,
\label{eq:terminal-XY-greater-section3}\\
X_k(y)&=a-c_{1}y\log(1/y)+O(y),
&Y(y)&=\pi c_{1}y(1+O(y)),
&&{\gamma_1}=1,
\label{eq:terminal-XY-one-section3}
\end{align}
where $c_1=2\pi\mu_{\rm out}(1)>0$ and
$a=k^2-4\pi\Acal(1)$ as in \eqref{eq:terminal-XY-section3}.
The case $\gamma_1=1$ uses the integer outer-endpoint formula
\eqref{eq:A-outer-integer-section2} and the explicit imaginary part
\eqref{eq:B-outer-section2}.

First let $k\uparrow k_*$, so $a\le0$.
The continuum density, whose variable remains
$0<y<\varepsilon_{\rm tr}$, satisfies
\[
\frac{Y}{X_k^2+Y^2}\lesssim
\begin{cases}
y^{{\gamma_1}-2},&{\gamma_1}>1,\\
[y\log^2(1/y)]^{-1},&{\gamma_1}=1.
\end{cases}
\]
On the continuum side $0<y<\varepsilon_{\rm tr}$, the cut density
converges pointwise to its
\(k=k_*\) value and is dominated by the preceding \(k\)-independent
integrable majorant. 
Meanwhile, $y_{\rm p}(k)\longrightarrow0^-.$
If \({\gamma_1}>1\), then
\[
r_{\rm p}(k)
=-\frac{k^3}{4\pi\Acal'(s_{\rm p}(k))}
\longrightarrow
-\frac{k_*^3}{4\pi\Acal'(1)}
=r_*,
\]
and therefore $\ri r_{\rm p}(k)\delta_{y_{\rm p}(k)}\rightharpoonup\ri r_*\delta_0$. 
If \({\gamma_1}=1\), the logarithmic endpoint expansion gives $\ri r_{\rm p}(k)\delta_{y_{\rm p}(k)}\rightharpoonup0$. 
Thus, as $k\uparrow k_*$, the continuum part converges in
\(L^1\), while the pole atom converges to
\(\ri r_*\delta_0\) for \({\gamma_1}>1\) and disappears in total variation
for \({\gamma_1}=1\).

For $k\downarrow k_*$, so $a>0$,
the real dispersion zero enters the continuum: $X_k$ has a unique
zero $y_a\in(0,\varepsilon_{\rm tr})$, with $y_a\to0$.
Since $Y(y_a)>0$, this gives a narrow resonance in the cut density
rather than an isolated pole.  Set
\[
 c_a=|\partial_yX_k(y_a)|,\qquad
 w_a=\frac{Y(y_a)}{c_a}.
\]
Linearizing $X_k$ and freezing $Y$ at $y_a$ gives the profile
\[
 \frac{Y(y_a)}{c_a^2(y-y_a)^2+Y(y_a)^2}
 =\frac1{c_a}\frac{w_a}{(y-y_a)^2+w_a^2}.
\]
Thus $w_a$ is the natural Lorentzian width.  We take a larger
resonant core $|y-y_a|\le L_aw_a$, with $L_a\to\infty$ to capture
its full limiting mass and $L_aw_a=o(y_a)$ to freeze the
coefficients.  The choices
$L_a=y_a^{-({\gamma_1}-1)/2}$ for ${\gamma_1}>1$ and
$L_a=\log(1/y_a)^{1/2}$ for ${\gamma_1}=1$ have both properties:
$w_a\asymp y_a^{\gamma_1}$ in the first case and
$w_a\asymp y_a/\log(1/y_a)$ in the second.
For $k$ sufficiently close to $k_*$, this core lies in $(y_a/2,2y_a)$,
and $y/y_a\to1$ uniformly on it.  The differentiated endpoint
expansions give $X'_k(y)/X'_k(y_a)\to1$ and
$Y(y)/Y(y_a)\to1$ uniformly there.  Integrating the first ratio
from $y_a$ to $y$ shows that, on $|y-y_a|\le L_aw_a$, we have
\[
X_k(y)=-c_a(y-y_a)(1+o(1)),\qquad
Y(y)=Y(y_a)(1+o(1)).
\]
Thus, for $\phi\in C_c^1(-{\varepsilon_{\rm tr}},{\varepsilon_{\rm tr}})$ and
$z=(y-y_a)/w_a$,
\begin{equation}
\label{eq:terminal-lorentzian-section3}
\int_{|y-y_a|\le L_aw_a}\phi(y)\,\rd\mathfrak m_k^+(y)
=\frac{\ri k^3}{\pi c_a}\phi(y_a)
\int_{-L_a}^{L_a}\frac{\rd z}{1+z^2}+o(1)
=\frac{\ri k^3}{c_a}\phi(y_a)+o(1).
\end{equation}

To control the remaining continuum and the Lorentzian tails, split
outside the core into
\[
(0,y_a/2),\qquad
(y_a/2,2y_a)\setminus
\{|y-y_a|\le L_aw_a\},\qquad
(2y_a,{\varepsilon_{\rm tr}}).
\]
For \({\gamma_1}>1\), monotonicity of \(X_k\), $|X'_k(y)|\gtrsim 1$ for $0<y<{\varepsilon_{\rm tr}}$ and
\(X_k(y_a)=0\) give
\[
|X_k(y)|\gtrsim
\begin{cases}
	y_a,&0<y<y_a/2,\\
	c_a|y-y_a|,&y_a/2<y<2y_a,\\
	y,&2y_a<y<{\varepsilon_{\rm tr}},
\end{cases}
\]
while \(Y(y)\asymp y^{\gamma_1}\).  Hence, on the first and third
regions,
\[
\frac{Y(y)}{X_k(y)^2+Y(y)^2}
\lesssim y^{{\gamma_1}-2}.
\]
For ${\gamma_1}=1$, integration of $|X'_k(y)|\asymp\log(1/y)$
from its zero gives
\[
 |X_k(y)|\gtrsim
 \begin{cases}
 y_a\log(1/y_a),&0<y<y_a/2,\\
 y\log(1/y),&2y_a<y<{\varepsilon_{\rm tr}}.
 \end{cases}
\]
For the second line one may restrict the derivative integral to
$[y/2,y]$; for the first, to $[y_a/2,y_a]$.  As
$y\log(1/y)$ is increasing near zero, both lines imply
\[
 \frac{Y(y)}{X_k(y)^2+Y(y)^2}
 \lesssim\frac1{y\log^2(1/y)}
\]
on their respective regions.
Both majorants are integrable at \(y=0\).  Since
\((0,y_a/2)\) shrinks to the endpoint, its contribution tends to
zero.  On \(2y_a<y<{\varepsilon_{\rm tr}}\), the density converges pointwise to its
critical value, and dominated convergence gives
\[
\int_{2y_a}^{{\varepsilon_{\rm tr}}}\phi(y)\,\rd\mathfrak m_k^+(y)
\longrightarrow
\int_0^{\varepsilon_{\rm tr}}
\phi(y)
\frac{\ri k_*^3Y(y)}
{\pi(X_{k_*}(y)^2+Y(y)^2)}
\,\rd y .
\]

On the middle annulus \(y\asymp y_a\),
\(Y(y)\asymp Y(y_a)\) and
\(|X_k(y)|\gtrsim c_a|y-y_a|\).
Thus, with \(w_a=Y(y_a)/c_a\),
\[
\int_{\substack{y_a/2<y<2y_a\\
		|y-y_a|>L_aw_a}}
\frac{Y(y)}
{X_k(y)^2+Y(y)^2}\,\rd y
\lesssim
\frac1{c_a}
\int_{|z|>L_a}\frac{\rd z}{1+z^2}
\lesssim L_a^{-1}.
\]
Hence the Lorentzian tail vanishes.

Together with
\eqref{eq:terminal-lorentzian-section3}, these estimates show that the
only concentrating mass is the central Lorentzian core; the remaining
continuum converges to the absolutely continuous part of the critical
measure. Indeed, if ${\gamma_1}>1$, then $c_a\to c_*$ and
\eqref{eq:terminal-lorentzian-section3} supplies the atom
$\ri r_*\delta_0$.  If ${\gamma_1}=1$, then
$c_a\sim c_1\log(1/y_a)$, so the concentrating mass tends to zero.
The core and complementary-region estimates also give a uniform
bound on the total variation.
The densities converge uniformly near the fixed chart boundaries,
so testing first with $C_c^1(-{\varepsilon_{\rm tr}},{\varepsilon_{\rm tr}})$ and then approximating
continuous test functions loses no boundary mass.  This proves both
one-sided limits and
\eqref{eq:terminal-critical-measure-greater-section3}--
\eqref{eq:terminal-critical-measure-one-section3}.
\end{proof}

For $|k-k_*|<\eta_{\rm tr}$, the change
$\tau=k(1-y)$ maps the local $y$-interval onto the positive edge
slice.  Define this slice and its reflection by
\[
\mathcal U_{+,k}
:=
\left\{
\tau\in\R:
\left|\frac{\tau}{k}-1\right|\le\varepsilon_{\rm tr}
\right\},
\qquad
\mathcal U_{-,k}
:=
\left\{
\tau\in\R:
\left|\frac{\tau}{k}+1\right|\le\varepsilon_{\rm tr}
\right\}.
\]
The choice of $\eta_{\rm tr}$ ensures that the pole frequencies
$\pm\tau_{\rm p}(k)$ lie strictly inside their respective slices
when $k<k_*$.

Push the local measure $\mathfrak m_k^+$ forward to the
$\tau$-coordinate by
\begin{equation}
\label{eq:terminal-pushforward-section3}
\int_\R\phi(\tau)\,\rd\mathfrak M_k^+(\tau)
:=\int_{-\varepsilon_{\rm tr}}^{\varepsilon_{\rm tr}}
\phi(k(1-y))\,\rd\mathfrak m_k^+(y).
\end{equation}
Thus $\mathfrak M_k^+$ is supported on $\mathcal U_{+,k}$.
The boundary symmetry following \eqref{eq:Omega-boundary-section2}
gives $R_\pm(-\tau,k)=\overline{R_\pm(\tau,k)}$, and hence
$J_k(-\tau)=\overline{J_k(\tau)}$; the pole residues likewise obey
$\varrho_-=\overline{\varrho_+}$.
The negative local measure, supported on $\mathcal U_{-,k}$, is
therefore fixed by Hermitian reflection:
\begin{equation}
\label{eq:terminal-negative-section3}
\int_\R\phi(\tau)\,\rd\mathfrak M_k^-(\tau)
:=\overline{\int_\R\overline{\phi(-\tau)}
\,\rd\mathfrak M_k^+(\tau)}.
\end{equation}
The measures $\mathfrak M_k^\pm$ contain the complete contribution
within the two edge slices, rather than the whole positive or
negative spectrum.  For $k<k_*$ each contains the pole atom and the
adjacent cut density; for $k>k_*$ it contains the local continuum.
At $k=k_*$ it contains the critical continuum and, when
$\gamma_1>1$, the corresponding atom from
Proposition~\ref{prop:terminal-measure-section3}.

On the remaining open cut intervals away from $\tau=0$,
$\Bcal(\tau/k)\ne0$.  Under \eqref{H0}, the center is regular by
\eqref{eq:filled-ball-zero-sign-appendix}.  At the ordinary inner
endpoints when $\kappa>0$, the limiting denominator is controlled
by $k^2-4\pi\Acal(\kappa)>0$; the outer endpoints are already
inside the edge slices.
These bounds are uniform for $k$ near $k_*$.  Consequently, the
ordinary density $J_k$, restricted to
$\Sigma_k\setminus(\mathcal U_{+,k}\cup\mathcal U_{-,k})$
and extended by zero outside this set, depends continuously on $k$
in $L^1_\tau$ near $k_*$.  We therefore define the complete physical
spectral measure near and at $k_*$ by
\begin{equation}
	\label{eq:full-physical-measure-section3}
		\rd\mathfrak M_k(\tau)
		=
		\rd\mathfrak M_k^+(\tau)
		+\rd\mathfrak M_k^-(\tau)
		+\one_{\Sigma_k\setminus
			(\mathcal U_{+,k}\cup\mathcal U_{-,k})}
		J_k(\tau)\,\rd\tau .
\end{equation}
For $k\ne k_*$ this is the decomposition of
\eqref{eq:physical-measure-nonterminal-section3} into its restrictions
to the two edge slices and their complement, so it recovers the
entire measure for $k\ne k_*$.  At $k=k_*$ it defines the full extension
using the common local limits.  Each critical atom for
$\gamma_1>1$ occurs once, in its corresponding edge piece; the complementary
term contains neither atom.  The three contributions are supported
on disjoint pieces apart from artificial boundary points of zero
mass.  Proposition~\ref{prop:terminal-measure-section3}
gives a uniform bound on the total variation of the two edge
pieces, while the ordinary cut contribution is uniformly \(L^1\).
Hence
\[
\sup_{|k-k_*|<\eta_{\rm tr}}
\|\mathfrak M_k\|_{\mathrm{TV}}<\infty .
\]
For every continuous test function on a common compact set containing
the nearby supports, $\phi\circ T_k\to\phi\circ T_{k_*}$ uniformly
on the fixed chart.  Together with the uniform total-variation bound,
this transfers weak convergence to the edge pieces in the
$\tau$-coordinate.  Away from $k_*$ their continuity follows directly
from the pole and density formulas; the ordinary part is continuous
in total variation.  Thus $k\mapsto\mathfrak M_k$ is weakly continuous
in a neighborhood of $k_*$.

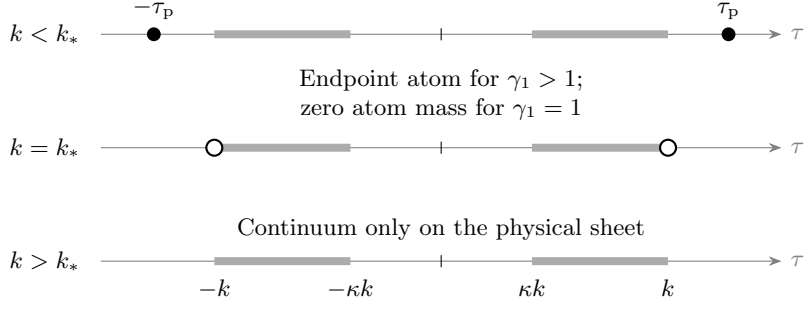
\begin{figure}[t]
\centering
\begin{tikzpicture}[x=1cm,y=1cm,>=Stealth,font=\small]
\foreach \y in {0,-1.5,-3.0}{
 \draw[->,gray] (-4.5,\y)--(4.5,\y) node[right] {$\tau$};
 \draw[line width=3pt,gray!65] (-3,\y)--(-1.2,\y);
 \draw[line width=3pt,gray!65] (1.2,\y)--(3,\y);
 \draw (0,\y-.08)--(0,\y+.08);
}
\fill (-3.8,0) circle (2.6pt);\fill (3.8,0) circle (2.6pt);
\node[above] at (-3.8,.08) {$-\tau_{\rm p}$};\node[above] at (3.8,.08) {$\tau_{\rm p}$};
\node[left,align=right] at (-4.65,0) {$k<k_*$};
\node[left,align=right] at (-4.65,-1.5) {$k=k_*$};
\node[left,align=right] at (-4.65,-3) {$k>k_*$};
\draw[thick,fill=white] (-3,-1.5) circle (3pt);\draw[thick,fill=white] (3,-1.5) circle (3pt);
\node[above,align=center] at (0,-1.25) {Endpoint atom for ${\gamma_1}>1$;\\zero atom mass for ${\gamma_1}=1$};
\node[below] at (-3,-3.12) {$-k$};\node[below] at (-1.2,-3.12) {$-\kappa k$};
\node[below] at (1.2,-3.12) {$\kappa k$};\node[below] at (3,-3.12) {$k$};
\node[above] at (0,-2.84) {Continuum only on the physical sheet};
\end{tikzpicture}
\caption{Support of the response spectral measure $\mathfrak M_k$
in the shell case, shown at three spatial frequencies. Gray segments denote the continuum; filled points denote isolated poles. The open edge markers indicate the location of contact, not a fixed atom mass. All scales are schematic.}
\label{fig:spectrum}
\end{figure}

Figure~\ref{fig:spectrum} summarizes the three physical spectral
regimes.  The critical atom for ${\gamma_1}>1$ and its absence for
${\gamma_1}=1$ are distinguished explicitly.

\subsection{Reconstruction of the resolvent and the response}

We now show that the Cauchy transform of the complete measure
\eqref{eq:full-physical-measure-section3} recovers the original
resolvent $R$, and hence identifies the Volterra response $G$.
\begin{proposition}
\label{prop:terminal-reconstruction-section3}
The family $k\mapsto\mathfrak M_k$ is weakly continuous through
$k_*$.  For every $\operatorname{Re}z>0$,
\begin{equation}
\label{eq:terminal-reconstruction-section3}
R(z,k)=\int_\R\frac{\rd\mathfrak M_k(\tau)}{z-\ri\tau},
\end{equation}
including at $k=k_*$.  Consequently,
\begin{equation}
\label{eq:terminal-time-reconstruction-section3}
\widehat G(t,k)=\one_{\{t>0\}}
\int_\R\re^{\ri t\tau}\,\rd\mathfrak M_k(\tau)
\end{equation}
is the Bromwich inverse of $R$ throughout a fixed neighborhood of
$k_*$; no distribution supported only at $k=k_*$ may be added.
\end{proposition}
\begin{proof}
For \(k\neq k_*\), the measure
\(\mathfrak M_k\) is the pole--cut measure
\eqref{eq:physical-measure-nonterminal-section3}.
The Cauchy identity \eqref{eq:cauchy-slit-appendix} in
Proposition~\ref{prop:fixed-k-representation-section3} therefore gives
\[
R(z,k)
=
\int_{\R}\frac{\rd\mathfrak M_k(\tau)}{z-\ri\tau},
\qquad \Re z>0.
\]
To pass through $k_*$, use the uniform total-variation bound and
weak continuity supplied by
Proposition~\ref{prop:terminal-measure-section3} and the assembly
\eqref{eq:full-physical-measure-section3}.  For fixed $\Re z>0$,
the Cauchy kernel is bounded and continuous on the real axis.
All nearby spectral supports lie in one compact set, so the kernel
may be multiplied by a fixed compactly supported test function equal
to one on that set.  Weak convergence thus allows passage to the
limit $k\to k_*$.  On the other hand,
$R(z,k)\to R(z,k_*)$ by \eqref{eq:dispersion-Laplace},
\eqref{eq:rho-resolvent-section2}, and the nonvanishing denominator in
Lemma~\ref{lem:right-half-plane-stability-section2}.  This proves the Cauchy identity \eqref{eq:terminal-reconstruction-section3} at
contact $k=k_*$.

To identify the time inverse in
\eqref{eq:terminal-time-reconstruction-section3}, define
\[
g_k(t)
:=
\one_{\{t>0\}}
\int_{\R}\re^{\ri t\tau}\,\rd\mathfrak M_k(\tau).
\]
Since \(\mathfrak M_k\) is finite, \(g_k\) is a bounded function,
extended by zero for negative times.  For \(\operatorname{Re}z>0\), Fubini's theorem gives
\[
\begin{aligned}
	\int_0^\infty \re^{-zt}g_k(t)\,\rd t
	=
	\int_{\R}
	\left(
	\int_0^\infty
	\re^{-(z-\ri\tau)t}\,\rd t
	\right)
	\rd\mathfrak M_k(\tau)=
	\int_{\R}
	\frac{\rd\mathfrak M_k(\tau)}{z-\ri\tau}
	=R(z,k).
\end{aligned}
\]
The use of Fubini is justified by $
\int_0^\infty\int_{\R}
\re^{-(\operatorname{Re}z)t}
\,\rd|\mathfrak M_k|(\tau)\,\rd t
=
\frac{\|\mathfrak M_k\|_{\rm TV}}
{\operatorname{Re}z}<\infty$. 

Thus $g_k$ and the original Volterra response
$\widehat G(\,\cdot\,,k)$ have the same Laplace transform on their
common half-plane $\Re z>1$, where the Volterra transform was
established in \eqref{eq:far-right-response-v3}.  Uniqueness of the Laplace transform for
exponential-order functions, and continuity for $t>0$, give
\eqref{eq:terminal-time-reconstruction-section3}.  This identifies
the spectral inverse with the original response also at $k=k_*$.
\end{proof}
The reconstructed response is bounded at each fixed spatial
frequency.  This improves the half-plane of convergence of the
actual Laplace integrals; the density estimate below records the
growth bound needed for that purpose.

\begin{corollary}
	\label{cor:fixed-k-boundedness-section3}
	For every fixed \(k\ge0\),
	\begin{equation}
		\label{eq:fixed-k-boundedness-section3}
		\sup_{t\ge0}|\widehat G(t,k)|<\infty.
	\end{equation}
	Moreover, for \(k\) in a fixed neighborhood of \(k_*\),
	\begin{equation}
		\label{eq:uniform-terminal-k-boundedness-section3}
		\sup_{|k-k_*|<\eta_{\rm tr}}\sup_{t\ge0}
		|\widehat G(t,k)|
		\le
		\sup_{|k-k_*|<\eta_{\rm tr}}\|\mathfrak M_k\|_{\mathrm{TV}}
		<\infty.
	\end{equation}
	Consequently, for every fixed \(k\ge0\), the actual Laplace transform
	$
	\int_0^\infty \re^{-zt}\widehat G(t,k)\,\rd t
	$
	converges for every \(\Re z>0\) and equals \(R(z,k)\).
Moreover, for every fixed $\xi$, the density represented by
\eqref{eq:rho-G-free-section2} satisfies
\begin{equation}\label{eq:density-growth-after-reconstruction-v3}
 |\widehat\rho(t,\xi)|\le C_\xi(1+t)\|h_{\rm in}\|_{L^1_{x,v}},
 \qquad
 \widetilde\rho(z,\xi)=Q(z,|\xi|)\widetilde{\rho^{\rm free}}(z,\xi),
 \quad\Re z>0,
\end{equation}
where both transforms in the last identity are actual Laplace integrals.
\end{corollary}
\begin{proof}
	For \(k>0\), \(k\neq k_*\), Proposition~\ref{prop:fixed-k-representation-section3} gives
	\[
	|\widehat G(t,k)|
	\le
	\one_{\{k<k_*\}}
	\sum_{\sigma=\pm}|\varrho_\sigma(k)|
	+
	\int_{\Sigma_k}|J_k(\tau)|\,\rd\tau,
	\]
	and the right-hand side is finite.  At \(k=k_*\),
	\eqref{eq:terminal-time-reconstruction-section3} gives
	$
	|\widehat G(t,k_*)|
	\le \|\mathfrak M_{k_*}\|_{\mathrm{TV}},
	$
	while the uniform estimate near \(k_*\) follows from the total-variation
	bound following \eqref{eq:full-physical-measure-section3}.  Finally,
	\(\widehat G(t,0)=-\sin t\) by
	\eqref{eq:fixed-k-zero-section3}.  This proves
	\eqref{eq:fixed-k-boundedness-section3}.
Boundedness makes the response Laplace integral absolutely convergent
and holomorphic on $\Re z>0$.  It agrees with the holomorphic
function $R$ on their common convergence half-plane $\Re z>1$ by
\eqref{eq:far-right-response-v3}.  The identity theorem for
holomorphic functions extends this equality to the connected
half-plane $\Re z>0$.
Finally, the bound
$|\widehat{\rho^{\rm free}}(t,\xi)|\le\|h_{\rm in}\|_{L^1_{x,v}}$
in \eqref{eq:rho-G-free-section2} gives the stated density growth.
Its actual transform therefore also converges on $\Re z>0$.
Fubini in that formula is justified by
\[
 \int_0^\infty\int_0^t\re^{-(\Re z)t}
 |\widehat G(t-s,k)\widehat{\rho^{\rm free}}(s,\xi)|\,\rd s\,\rd t
 \le C_k(\Re z)^{-2}\|h_{\rm in}\|_{L^1_{x,v}},
\]
and gives the last identity in
\eqref{eq:density-growth-after-reconstruction-v3}.
\end{proof}

\subsection{The wave--damping decomposition}

We now realize \eqref{eq:green-wave-LD-intro} at the resolvent level:
the smooth pole pair and the complete threshold block form the wave
component, and the remaining continuum forms the damped component.
Choose $\varepsilon>0$ with $3\varepsilon<\eta_{\rm tr}$ and the
partition in \eqref{eq:k-partition-intro} so that
\begin{align*}
&\supp\chi_\ell\subset[0,k_*-2\varepsilon),
\qquad\chi_\ell=1\quad\hbox{on }[0,k_*-3\varepsilon],\\
&\supp\vartheta\subset(k_*-3\varepsilon,k_*+3\varepsilon),
\qquad \vartheta=1\quad\hbox{if }|k-k_*|\le2\varepsilon,\\
&\supp\chi_h\subset(k_*+2\varepsilon,\infty),
\qquad\chi_h=1\quad\hbox{on }[k_*+3\varepsilon,\infty).
\end{align*}
For $k<k_*$ let
\begin{equation}
\label{eq:pole-rational-section3}
P(z,k):=\sum_{\sigma=\pm}
\frac{\varrho_\sigma(k)}{z-\ri\sigma\tau_{\rm p}(k)},
\end{equation}
using the continuous values $\tau_{\rm p}(0)=1$ and
$\varrho_\pm(0)=\pm\ri/2$ at $k=0$, and set $P=0$ for $k\ge k_*$.
The product $\chi_\ell P$ is unaffected by
this extension because $\chi_\ell$ vanishes on a full neighborhood of $k_*$.
Define
\begin{equation}
\label{eq:three-blocks-section3}
R_{\rm KG}=\chi_\ell P,\qquad
R_{\rm th}=\vartheta R,\qquad
R_{\rm LD}=\chi_\ell(R-P)+\chi_h R.
\end{equation}
Then, before taking any boundary value,
\begin{equation}
\label{eq:three-block-sum-section3}
R=R_{\rm KG}+R_{\rm th}+R_{\rm LD}.
\end{equation}
Indeed, below contact this is
$\chi_\ell P+\vartheta R+\chi_\ell(R-P)=(\chi_\ell+\vartheta)R$, while above contact it is
$(\vartheta+\chi_h)R$.

Let \(\widehat G_{\rm KG}\), \(\widehat G_{\rm th}\), and
\(\widehat G_{\rm LD}\) denote the Bromwich inverses of the three
terms in \eqref{eq:three-blocks-section3}, and let
\(G_{\rm KG}\), \(G_{\rm th}\), and \(G_{\rm LD}\) be their
inverse spatial Fourier transforms.  Away from \(k_*\), their
fixed-frequency representations follow from
Proposition~\ref{prop:fixed-k-representation-section3}.  At
\(k=k_*\), Proposition~\ref{prop:terminal-reconstruction-section3}
defines the threshold inverse by
\[
\widehat G_{\rm th}(t,k_*)
=
\one_{\{t>0\}}
\int_{\R}\re^{\ri t\tau}\,\rd\mathfrak M_{k_*}(\tau),
\]
since \(\vartheta(k_*)=1\). Thus
\begin{equation}
\label{eq:wave-components-section3}
G_{\rm wave}=G_{\rm KG}+G_{\rm th},\qquad
G=G_{\rm wave}+G_{\rm LD},
\end{equation}
and this is exactly the definition in
\eqref{eq:green-wave-LD-intro}. 
Each pole or spectral inverse here is understood for $t>0$ and
extended by zero to $t<0$.  Spatial multipliers preserve this forward
time support, so all three components are causal response kernels.

\subsection{Frequency and endpoint structure of the regular response}
\label{subsec:low-frequency-regular-section3}

We record the low- and high-frequency endpoint bounds for the
continuum jump and the quadratic low-frequency factor obtained after
subtraction of the pole pair.
The last assertion gives the corresponding Cauchy representation
on the slit plane.
Recall that $
\mathcal I_\kappa=[-1,-\kappa]\cup[\kappa,1],\,
\Sigma_k=k\mathcal I_\kappa.$
For $a\in\{\kappa,1\}$ with $a>0$ (only $a=1$ under \eqref{H0}), let
\[
s_a(u):=a+\epsilon_a u,
\qquad
\epsilon_\kappa=1,\qquad
\epsilon_1=-1,
\qquad
D_u:=u\partial_u,
\]
for $0<u<u_0$, where $u_0$ lies within the endpoint neighborhoods
in \eqref{H2} or \eqref{H0} and is small enough that each chart
meets no other velocity endpoint.  The negative endpoints are obtained
from the evenness of $\Acal$ and the oddness of $\Bcal$.
Each endpoint estimate uses its own order $\gamma_a$; the other
endpoint contributes a smooth function on the local chart.
The measure and geometry near contact use the outer order $\gamma_1$.

\begin{proposition}
	\label{prop:low-frequency-factor-section3}
	There exist constants $0<k_0<k_*,\, k_\infty>k_*,$
	such that the following statements hold.
{
Every estimate below is valid for arbitrary fixed nonnegative
integers $m,n$, with constants depending on that pair, the fixed
endpoint orders, the equilibrium and the cutoffs.
}

\par\medskip\noindent
\textup{(i)}\enspace
For \(0<k\le k_0\) and \(s\in\mathcal I_\kappa\),
	\begin{equation}
		\label{eq:jump-low-frequency-section3}
		J_k(ks)
		=
		k^2j_0(s)+k^4j_1(k,s),
		\qquad
		j_0(s)
		=
		\frac{2\ri\Bcal(s)}{|\Om(s)|^2}.
	\end{equation}
	For $a\in\{\kappa,1\}$ with $a>0$,
	\begin{equation}
		\label{eq:jump-conormal-section3}
		\left|
		(k\partial_k)^mD_u^n j_1(k,s_a(u))
		\right|
		+
		\left|
		D_u^n j_0(s_a(u))
		\right|
		\le C_{m,n}u^{\gamma_a} ,
	\end{equation}
	and the same estimates hold at \(-s_a(u)\).  If
	\(K\Subset\mathcal I_\kappa^\circ\), then
	\begin{equation}
		\label{eq:jump-interior-low-section3}
		\sup_{\substack{0<k\le k_0\\s\in K}}
		\left(
		|(k\partial_k)^m\partial_s^n j_1(k,s)|
		+
		|\partial_s^n j_0(s)|
		\right)
		\le C_{K,m,n}.
	\end{equation}
	
\par\medskip\noindent
\textup{(ii)}\enspace
At sufficiently high frequency,
	\begin{equation}
		\label{eq:jump-high-frequency-section3}
		\left|
		\partial_k^mD_u^n
		J_k\bigl(k s_a(u)\bigr)
		\right|
		\le
		C_{m,n}k^{-2-m}u^{\gamma_a},
		\qquad
		k\ge k_\infty,
	\end{equation}
	with the same estimate at the negative endpoints.  Moreover, for
	\(K\Subset\mathcal I_\kappa^\circ\),
	\begin{equation}
		\label{eq:jump-interior-high-section3}
		\sup_{s\in K}
		\left|
		\partial_k^m\partial_s^n J_k(ks)
		\right|
		\le
		C_{K,m,n}k^{-2-m},
		\qquad
		k\ge k_\infty.
	\end{equation}
	
\par\medskip\noindent
\textup{(iii)}\enspace
After subtraction of the pole pair, the low-frequency response
gains two powers of \(k\).  For the boundary side $\varsigma\in\{+,-\}$, set
$(R-P)_{\varsigma}(\tau,k):=R_{\varsigma}(\tau,k)-P(\ri\tau,k)$.
On the physical cut,
	\begin{equation}
		\label{eq:R-minus-P-factor-section3}
		(R-P)_{\varsigma}(ks,k)
		=
		k^2\mathcal R_{\varsigma}(k,s),
		\qquad
		0<k\le k_0,
	\end{equation}
	where, for either endpoint sign $\sigma\in\{-1,1\}$ independently,
	\begin{align}
		\left|
		(k\partial_k)^mD_u^n
		\mathcal R_{\varsigma}(k,\sigma s_a(u))
		\right|
		&\le C_{m,n},
		\label{eq:Rcal-conormal-section3}\\
		\left|
		(k\partial_k)^mD_u^n
		\bigl(\mathcal R_+-\mathcal R_-\bigr)
		(k,\sigma s_a(u))
		\right|
		&\le C_{m,n}u^{\gamma_a} .
		\label{eq:Rcal-jump-conormal-section3}
	\end{align}
For every $K\Subset\mathcal I_\kappa^\circ$, one also has
\begin{equation}
\label{eq:Rcal-interior-section3}
 \sup_{\substack{0<k\le k_0\\s\in K}}
 \left(
 |(k\partial_k)^m\partial_s^n\mathcal R_{\varsigma}(k,s)|
 +|(k\partial_k)^m\partial_s^n
       (\mathcal R_+-\mathcal R_-)(k,s)|
 \right)\le C_{K,m,n}.
\end{equation}
	
\par\medskip\noindent
\textup{(iv)}\enspace
For the same range $0<k\le k_0$, set
$\mathscr S_\kappa:=\C\setminus\ri\mathcal I_\kappa$.  Then
	\begin{equation}
		\label{eq:R-minus-P-slit-factor-section3}
		R(k\zeta,k)-P(k\zeta,k)
		=
		k^2\widetilde{\mathcal R}(\zeta,k),
		\qquad
		\zeta\in\mathscr S_\kappa,
	\end{equation}
	where
	\begin{equation}
		\label{eq:Rtilde-cauchy-section3}
		\widetilde{\mathcal R}(\zeta,k)
		=
		\int_{\mathcal I_\kappa}
		\frac{j_0(s)+k^2j_1(k,s)}
		{\zeta-\ri s}\,\rd s .
	\end{equation}
	In particular, \(R-P\) extends holomorphically across the two pole
	points, and for every compact \(K\Subset\mathscr S_\kappa\),
	\begin{equation}
		\label{eq:Rtilde-symbol-section3}
		\sup_{\substack{0<k\le k_0\\\zeta\in K}}
		\left|
		(k\partial_k)^m\partial_\zeta^n
		\widetilde{\mathcal R}(\zeta,k)
		\right|
		\le C_{K,m,n},
		\qquad
		m,n\ge0.
	\end{equation}

\end{proposition}

\begin{proof}
\emph{Part (i): the low-frequency jump.}
On the cut, \eqref{eq:jump-density-section2} gives
	\begin{equation}
		\label{eq:J-explicit-proof-section3}
		J_k(ks)
		=
		\frac{2\ri k^2\Bcal(s)}
		{d(k,s)},
		\qquad
		d(k,s)
		:=
		\bigl(k^2-4\pi\Acal(s)\bigr)^2+
		\bigl(2\pi\Bcal(s)\bigr)^2.
	\end{equation}
	At \(k=0\), $d(0,s)=|\Om(s)|^2$.
This is bounded away from zero on the closed set
$\mathcal I_\kappa$: in its interior away from zero $\Bcal\ne0$,
at $s=\pm1$ one has $4\pi\Acal(1)=k_*^2$, and at $s=\pm\kappa$
when $\kappa>0$ one has $\Acal(\kappa)<0$.
For \eqref{H0}, \eqref{eq:filled-ball-zero-sign-appendix} controls
the remaining point $s=0$.  Thus
	\begin{equation}
		\label{eq:D0-lower-proof-section3}
		d(0,s)\ge c_0>0,
		\qquad
		d(k,s)\ge c_0/2
	\end{equation}
	for \(s\in\mathcal I_\kappa\) and \(0<k\le k_0\), after decreasing
	\(k_0\) if necessary.
	
Since $d(k,s)=d(0,s)-8\pi k^2\Acal(s)+k^4$, the identity
\[
 \frac1{d(k,s)}-\frac1{d(0,s)}
 =\frac{k^2(8\pi\Acal(s)-k^2)}{d(0,s)d(k,s)}
\]
gives the exact expansion
\[
 J_k(ks)
 =k^2\frac{2\ri\Bcal(s)}{d(0,s)}
 +k^4\frac{2\ri\Bcal(s)(8\pi\Acal(s)-k^2)}{d(0,s)d(k,s)}.
\]
	This proves \eqref{eq:jump-low-frequency-section3}, with
	\begin{equation}
		\label{eq:j1-explicit-section3}
		j_1(k,s)
		=
		\frac{2\ri\Bcal(s)(8\pi\Acal(s)-k^2)}
		{d(0,s)d(k,s)}.
	\end{equation}
The endpoint formulas in Lemma~\ref{lem:A-endpoint-section2},
together with the definition of $\Bcal$ in \eqref{eq:A-B-section2}
and the endpoint factors in \eqref{H2} or \eqref{H0}, imply,
for every fixed $n$,
	\begin{equation}\label{eq:conormal-est-A-B}
		|D_u^n\Bcal(\pm s_a(u))|
		\lesssim u^{\gamma_a},
		\qquad
		|D_u^n\Acal(\pm s_a(u))|
		\lesssim1.
	\end{equation}
	Together with \eqref{eq:D0-lower-proof-section3} and the conormal
	product rule, the explicit formulas for \(j_0\) and \(j_1\) therefore
	give \eqref{eq:jump-conormal-section3}.  On a compact
	\(K\Subset\mathcal I_\kappa^\circ\), the functions \(\Acal\) and
	\(\Bcal\) are smooth and the same denominators remain uniformly
	nonzero; ordinary differentiation of
	\eqref{eq:j1-explicit-section3} proves
	\eqref{eq:jump-interior-low-section3}.
	
\emph{Part (ii): high frequencies.}
Normalize the denominator by $k^4$.  Choose \(k_\infty\) so large that
\[
	\sup_{s\in\mathcal I_\kappa}
	\left(
	\frac{|4\pi\Acal(s)|}{k_\infty^2}
	+
	\frac{|2\pi\Bcal(s)|}{k_\infty^2}
	\right)
	\le\frac14.
\]
	Then
	\[
	J_k(ks)
	=
	2\ri k^{-2}\Bcal(s)
	\left[
	\left(1-\frac{4\pi\Acal(s)}{k^2}\right)^2
	+
	\left(\frac{2\pi\Bcal(s)}{k^2}\right)^2
	\right]^{-1},
	\]
	and the bracket is uniformly bounded above and below for
	\(k\ge k_\infty\).  Differentiating this expression in \(k\), and
	using \eqref{eq:conormal-est-A-B}, proves
	\eqref{eq:jump-high-frequency-section3}.  Ordinary differentiation on
	\(K\Subset\mathcal I_\kappa^\circ\) proves
\eqref{eq:jump-interior-high-section3}.

\emph{Part (iii): subtraction on the cut.}
We compare the common order-zero terms of $R$ and $P$.
Set $\Om_+(s)=\Om(s)$ and $\Om_-(s)=\overline{\Om(s)}$.
From \eqref{eq:Omega-boundary-section2},
	\[
	R_{\varsigma}(ks,k)
	=
	-1+\frac{k^2}{k^2-\Om_{\varsigma}(s)}
	=
	-1-\frac{k^2}{\Om_{\varsigma}(s)}+O_{\rm con}(k^4).
	\]
	Here the $O_{\rm con}(k^4)$ bound is preserved under each fixed
	$(k\partial_k)^mD_u^n$ on the endpoint charts, and under
	$(k\partial_k)^m\partial_s^n$ on compact subsets of the cut interior,
	with constants depending on $m,n$ and, in the latter case, the compact set.
	On the other hand,
	\[
	P(\ri ks,k)
	=
	-\frac{2r_{\rm p}(k)\tau_{\rm p}(k)}
	{\tau_{\rm p}(k)^2-k^2s^2}.
	\]
	The low-frequency expansions
	\eqref{eq:taup-low-section2}--\eqref{eq:residue-low-section2} give $
	2r_{\rm p}(k)\tau_{\rm p}(k)=1+O(k^4),\, 
	\tau_{\rm p}(k)^2
	=
	1+4\pi m_4k^2+O(k^4),$
	and hence, uniformly for \(s\in\mathcal I_\kappa\),
	\[
	P(\ri ks,k)
	=
	-1+\bigl(4\pi m_4-s^2\bigr)k^2+O(k^4).
	\]
Since $\tau_{\rm p}$ and $r_{\rm p}$ are smooth in $k^2$ by
Proposition~\ref{prop:pole-branch-section2} and the displayed
denominator stays bounded away from zero for small $k$, this
$O(k^4)$ remainder retains its bound under every fixed
$(k\partial_k)^m\partial_s^n$, hence also under
$(k\partial_k)^mD_u^n$ on the endpoint charts.
	Subtracting yields \eqref{eq:R-minus-P-factor-section3}; in particular, $
	\mathcal R_{\varsigma}(0,s)
	=
	s^2-4\pi m_4-\frac1{\Om_{\varsigma}(s)}$. 
The conormal bounds \eqref{eq:conormal-est-A-B} and the lower
bound \eqref{eq:D0-lower-proof-section3} for $|\Om_{\varsigma}|^2$ give
	\eqref{eq:Rcal-conormal-section3}.  Since \(P\) has no cut jump, $
	k^2(\mathcal R_+-\mathcal R_-)
	=
	R_+-R_-
	=
	2\pi J_k$, 
	so \eqref{eq:jump-low-frequency-section3} and
	\eqref{eq:jump-conormal-section3} give
\eqref{eq:Rcal-jump-conormal-section3}.  On a compact cut-interior
set, ordinary $s$ differentiation of the same expressions gives
\eqref{eq:Rcal-interior-section3}.
	
\emph{Part (iv): the slit-plane representation.}
Since $k_0<k_*$, subtracting the pole fractions in
\eqref{eq:cauchy-slit-appendix} of
Proposition~\ref{prop:fixed-k-representation-section3} gives
\[
 R(z,k)-P(z,k)
 =\int_{\Sigma_k}\frac{J_k(\tau)}{z-\ri\tau}\,\rd\tau,
 \qquad \Re z>0.
\]
The principal parts of $R$ and $P$ agree at the two simple poles,
so their difference extends holomorphically across those points.
The right-hand side is holomorphic on the connected slit plane
$\C\setminus\ri\Sigma_k$.  The identity theorem for holomorphic
functions therefore extends the equality from $\Re z>0$ to this
entire domain.  Now put $\tau=ks$ and $z=k\zeta$:
\[
 R(k\zeta,k)-P(k\zeta,k)
 =\int_{\mathcal I_\kappa}
   \frac{kJ_k(ks)}{k(\zeta-\ri s)}\,\rd s
 =k^2\int_{\mathcal I_\kappa}
   \frac{j_0(s)+k^2j_1(k,s)}{\zeta-\ri s}\,\rd s.
\]
Here the Jacobian $k$ cancels the denominator's factor $k$, and
\eqref{eq:jump-low-frequency-section3} supplies the remaining
$k^2$.  This proves \eqref{eq:R-minus-P-slit-factor-section3}--
\eqref{eq:Rtilde-cauchy-section3} on $\mathscr S_\kappa$.
  On every compact
	\(K\Subset\mathscr S_\kappa\), the Cauchy kernel and its fixed finite
	\(\zeta\)-derivatives are uniformly bounded; the estimates already
	proved for \(j_0,j_1\) then give
	\eqref{eq:Rtilde-symbol-section3}.
\end{proof}

\section{Pointwise bounds for the Green kernel}
\label{sec:green}

We prove the pointwise estimates for the Green kernel.  The
physical-space Volterra series first gives finite propagation and
the short-time bounds.  We then estimate the smooth pole wave, the
Landau-damped continuum and its fronts, and the threshold block,
before combining these bounds to prove Theorem~\ref{thm:main-green}.

\subsection{Volterra bounds and finite propagation}

The physical-space Volterra series gives finite propagation and uniform bounds on bounded time intervals.

\begin{lemma}
\label{lem:finite-propagation-green-section4}
\label{lem:short-time-green-section4}
For \(t>0\), the physical-space Volterra kernel is
\begin{equation}\label{eq:first-iterate-section4}
 K(t,x)=-t^{-2}f^0(x/t).
\end{equation}
The Volterra iterates $\mathcal K_1=K$,
$\mathcal K_{n+1}=K*_{t,x}\mathcal K_n$ satisfy
\begin{equation}\label{eq:volterra-iterates-v4}
 \|\mathcal K_n(t)\|_1\le C^n\frac{t^{2n-1}}{(2n-1)!},\qquad
 \|\mathcal K_n(t)\|_\infty
 \le C^n\frac{n^3t^{2n-4}}{(2n-2)!}\quad(n\ge2).
\end{equation}
The response defined by \eqref{eq:green-volterra-definition}
has the physical-space representation
\begin{equation}\label{eq:physical-volterra-series-v3}
 G(t)=\sum_{n=1}^\infty\mathcal K_n(t)\quad\hbox{in }L^1_x,
\end{equation}
absolutely and uniformly in \(L^1_x\) on every bounded time interval.
Writing $\overline B(0,t)=\{x\in\R^3:|x|\le t\}$ for the closed ball,
we have
\begin{align}
 \supp G(t)&\subset\overline B(0,t),\label{eq:finite-propagation-green-section4}\\
 \sup_{0<t\le2}\|\ang x^mG(t)\|_1&<\infty
                \quad\text{for each fixed }m\ge0,\label{eq:short-weighted-G-v4}\\
 \sup_{0<t\le1}\|G(t)-K(t)\|_\infty&<\infty,\label{eq:short-time-response-section4}\\
 |G_{\rm LD}(t,x)|&\le C(1+t^{-2}\one_{\{\kappa t\le|x|\le t\}}),
                       \quad0<t\le1.\label{eq:short-time-LD-section4}
\end{align}
\end{lemma}
\begin{proof}
The change of variables $x=tv$ gives
\[
 \mathcal F_x[-t^{-2}f^0(x/t)](\xi)
 =-t\widehat f^0(t\xi)=\widehat K(t,\xi),
\]
which proves \eqref{eq:first-iterate-section4} by \eqref{eq:kernel-K}.
Write $\Delta_{n-1}(t)=\{t_j>0:\sum_{j=1}^nt_j=t\}$, with measure
$\rd t_1\cdots\rd t_{n-1}$ and $t_n=t-\sum_{j<n}t_j$.  Iteration gives
\[
 \mathcal K_n(t)=\int_{\Delta_{n-1}(t)}
 K(t_1)*_x\cdots *_xK(t_n)\,\rd t_1\cdots\rd t_{n-1}.
\]
The simplex identity
\[
 \int_{\Delta_{n-1}(t)}\prod_{j=1}^nt_j^{a_j-1}\,\rd t_1\cdots\rd t_{n-1}
 =t^{\sum a_j-1}\frac{\prod_j\Gamma(a_j)}{\Gamma(\sum_j a_j)}
 \qquad(a_j>0)
\]
follows by scaling $t_j=t u_j$ and repeated beta integration;
$\Gamma$ denotes Euler's gamma function.  Since
$\|K(t_j)\|_1\le Ct_j$, its case $a_j=2$ proves the first estimate.
For the second, partition the simplex into the $n$ sets on which a
specified $t_j$ is maximal, assigning ties arbitrarily.  There
$t_j\ge t/n$ and
\[
 \|K(t_1)*_x\cdots *_xK(t_n)\|_\infty
 \le C^n n^2t^{-2}\prod_{\ell\ne j}t_\ell .
\]
Enlarge each of these sets to the full simplex and apply the identity
with $a_j=1$ and $a_\ell=2$ for $\ell\ne j$.  Its integral is
$t^{2n-2}/(2n-2)!$, and summation over the $n$ choices gives the factor
$n^3$ in \eqref{eq:volterra-iterates-v4}.

The first bound makes the series in
\eqref{eq:physical-volterra-series-v3} absolutely and uniformly
convergent in \(L^1_x\) on bounded time intervals.  Since
$\widehat{\mathcal K_n}(t,\xi)=(\widehat K)^{*n}(t,\xi)$,
its spatial Fourier transform is the series \eqref{eq:volterra-series-section2};
by the uniqueness argument immediately following that series, it is exactly the
previously defined response \(\widehat G\).  This identifies the
physical-space sum with \(G\).

For each point of the simplex,
$\overline B(0,t_1)+\cdots+\overline B(0,t_n)\subset\overline B(0,t)$.
Every iterate is therefore supported in this common closed ball,
as is its $L^1$ limit.  Consequently, for every $T<\infty$,
\[
 \sup_{0<t\le T}\|\ang x^m G(t)\|_1
 \le \ang T^m\sum_{n\ge1}C^n\frac{T^{2n-1}}{(2n-1)!}<\infty.
\]
Finally $\sum_{n\ge2}C^nn^3t^{2n-4}/(2n-2)!$ is uniformly bounded
for $0<t\le1$, proving \eqref{eq:short-time-response-section4}.
For the pole block, \eqref{eq:pole-rational-section3}--
\eqref{eq:three-blocks-section3} and the pole inverse
\eqref{eq:bromwich-model-section2} give
\begin{equation}\label{eq:KG-fourier-representation-section4}
 \widehat G_{\rm KG}(t,\xi)
 =\chi_\ell(k)\sum_{\sigma=\pm}\varrho_\sigma(k)
        \re^{\ri\sigma t\tau_{\rm p}(k)},\qquad k=|\xi|,\quad t>0.
\end{equation}
Proposition~\ref{prop:pole-branch-section2} shows that the pole
frequency and residues are smooth functions of $k^2$ near zero.
Since $\chi_\ell=1$ there, the compactly supported amplitudes
$\chi_\ell(|\xi|)\varrho_\sigma(|\xi|)$ are smooth in Cartesian coordinates,
including at $\xi=0$.  Direct Fourier inversion therefore gives
\[
 \|G_{\rm KG}(t)\|_\infty
 \le\frac1{(2\pi)^3}\sum_{\sigma=\pm}
       \int_{\R^3}|\chi_\ell(|\xi|)\varrho_\sigma(|\xi|)|\,\rd\xi
 \le C.
\]
For the threshold block, the $L^1$ bound already proved above yields
\[
 \|G_{\rm th}(t)\|_\infty
 =\|\check\vartheta*_xG(t)\|_\infty
 \le\|\check\vartheta\|_\infty\|G(t)\|_1\le C,
 \qquad 0<t\le1.
\]
Using $G_{\rm LD}=(G-K)+K-G_{\rm KG}-G_{\rm th}$ now proves
\eqref{eq:short-time-LD-section4}, since the remaining singular term
is $K=-t^{-2}f^0(x/t)$, supported on $\kappa t\le|x|\le t$.
\end{proof}

\subsection{The smooth pole component}

The two amplitude classes needed here and for the electric field are
treated together in Lemma~\ref{lem:radial-dispersive-section4} in
Appendix~\ref{app:endpoints}.  The scalar pole amplitude is Cartesian
smooth, so the first conclusion of that lemma supplies any fixed
spatial weight.
\begin{proposition}
\label{prop:pole-dispersive-outline}
For all $N\ge0$ and $t>0$, 
\begin{equation}
\label{eq:KG-weighted-section4}
|G_{\rm KG}(t,x)|
\le C_N\ang t^{-3/2}
\frac{\ang t^N}{\ang{t,x}^N}.
\end{equation}
\end{proposition}

\begin{proof}
Applying spatial Fourier inversion
\eqref{eq:fourier-transform-intro} to
\eqref{eq:KG-fourier-representation-section4} gives
\[
G_{\rm KG}(t,x)=\frac1{(2\pi)^3}
\sum_{\sigma=\pm}\int_{\R^3}
\re^{\ri x\cdot\xi+\ri\sigma t\tau_{\rm p}(|\xi|)}
\chi_\ell(|\xi|)\varrho_\sigma(|\xi|)\,\rd\xi .
\]
To apply Lemma~\ref{lem:radial-dispersive-section4}, we check its
Hessian condition on a neighborhood of the amplitude support.
On the pole branch, the chain rule
\eqref{eq:pole-chain-rule-Green-revision} gives
\[
 q'(\rho)=\Om_0(\rho)+\rho\Om_0'(\rho)>0,\qquad
 \tau_{\rm p}''(q(\rho))=\frac{V'(\rho)}{q'(\rho)}>0,
 \qquad 0<\rho<1.
\]
Here $V'>0$ follows from \eqref{eq:Vprime-identity-section2} and
the strict lower bound \eqref{eq:A-rigid-positive-lower-section2}.
Also $\tau_{\rm p}'(k)/k>0$ for $k>0$, since $V>0$.
At the origin, \eqref{eq:taup-low-section2}--\eqref{eq:taup-hessian-section2} give
\[
 \tau_{\rm p}''(0)
 =\lim_{k\downarrow0}\frac{\tau_{\rm p}'(k)}k
 =4\pi m_4>0.
\]
Since $\supp\chi_\ell$ lies a fixed distance below $k_*$, continuity gives
a common positive lower bound for these two eigenvalues on a
slightly larger compact frequency ball.  The phase $-\tau_{\rm p}$
has the corresponding negative-definite Hessian; the lemma treats
it by conjugation.  The amplitudes are compactly supported and
Cartesian smooth, as verified after
\eqref{eq:KG-fourier-representation-section4}.  Applying the Cartesian
conclusion of Lemma~\ref{lem:radial-dispersive-section4} to both phases proves
\eqref{eq:KG-weighted-section4} for integer $N\ge0$; the remaining
real orders follow by increasing $N$ to the next integer.
\end{proof}

\subsection{The Landau-damped continuum, fronts, and far field}
The ordinary fronts are the rays $|x|=at$, $a\in\{\kappa,1\}$, $a>0$.
We estimate conic neighborhoods of these fronts separately from the
interior and exterior regions.

\begin{proposition}
\label{prop:green-LD-outline}
For $t>0$ and $t+|x|\ge1$, the Landau-damped component satisfies
\begin{equation}
\label{eq:green-LD-outline}
 |G_{\rm LD}(t,x)|
 \le \frac C{\ang{t,x}^{6}}
 \left[1+\sum_{\substack{a\in\{\kappa,1\}\\a>0}}
 \left(\frac{\ang{t,x}}{\ang{|x|-at}}\right)^{(4-\gamma_a)_+}\right].
\end{equation}
This establishes the large-space--time part of
\eqref{eq:main-green-LD-envelope}.
\end{proposition}

\begin{proof}
Put $r:=|x|$ and $X:=\ang{t,x}$. The comparison
\[
 X\le1+t+r\le\sqrt3X
\]
allows the same fixed-scale decompositions with $X$; the phases and
front distances remain $tk$ and $r-at$.
Set
	\[
	A(k,s):=(1-\vartheta(k))J_k(ks).
	\]
	Since
	\[
	R_{\rm LD}=\chi_\ell(R-P)+\chi_hR,
	\qquad
	\chi_\ell+\chi_h=1-\vartheta,
	\]
	the pole part has been removed on the support of \(\chi_\ell\), while no pole
	is present on the support of \(\chi_h\).  Thus
\eqref{eq:fixed-k-representation-section3}, together with
\eqref{eq:three-blocks-section3}, gives
	\begin{equation}
		\label{eq:LD-joint-integral-section4}
		G_{\rm LD}(t,x)
		=
		\frac1{(2\pi)^3}
		\int_{\mathbb S^2}\int_0^\infty\int_{\mathcal I_\kappa}
		\re^{\ri k(x\cdot\theta+ts)}
		A(k,s)k^3\,\rd s\,\rd k\,\rd\theta .
	\end{equation}
	The factor \(k^3\) is the product of the radial spatial measure
	\(k^2\rd k\) and the Jacobian
	\(\rd\tau=k\,\rd s\).
	{
	When the high-frequency integral is not absolutely convergent, its
	meaning is the limit obtained by inserting $\chi(k/H)$, where
	$\chi\in C_c^\infty([0,\infty))$ equals one near zero and
	$H\to\infty$.  The smooth cutoffs at $k\asymp X^{-1}$
	used below are retained throughout this limit.  On $k\asymp H$,
	$|\partial_k^j\chi(k/H)|\le C_j k^{-j}$, so every integration by
	parts below includes these terms with the same symbol order.
	The resulting differentiated tails are integrable and their
	$H$-annuli tend to zero.  This proves existence and independence of
	the truncation pointwise in the regions under consideration, as well
	as agreement with the distributional Fourier inverse in
	\eqref{eq:LD-joint-integral-section4}.
	}
	
	We may assume that \(X\) is sufficiently large.
The bounded range is treated directly.  The iterate bounds
\eqref{eq:volterra-iterates-v4}, together with $K=-t^{-2}f^0(x/t)$,
give $\sup_{1/4\le t\le T}\|G(t)\|_\infty<\infty$ for each fixed $T$.
The two frequency-localized components are bounded on this interval
as in Lemma~\ref{lem:short-time-green-section4}.  If $t<1/4$ and
$t+r\ge1$, then $r>t$, so finite propagation gives
$G_{\rm LD}=-G_{\rm KG}-G_{\rm th}$, which is again bounded.
	
	\medskip
	\noindent
	\emph{Low frequencies.}
	By Proposition~\ref{prop:low-frequency-factor-section3}(i),
	for \(0<k\le k_0\),
	\begin{equation}
		\label{eq:LD-low-k-symbol-section4}
		A(k,s)
		=
		k^2j_0(s)+k^4j_1(k,s),
	\end{equation}
	with uniform finite conormal bounds.  In particular,
	\[
	|A(k,s)|\lesssim k^2.
	\]
Choose a fixed sufficiently large constant $c>0$ and a smooth function
$\psi:[0,\infty)\to[0,1]$ equal to zero on $[0,c]$ and one on
$[2c,\infty)$.  We split with $1-\psi(X k)$ and
$\psi(X k)$; the same lower cutoff is used below.
The first piece satisfies
	\[
	\left|\int_{\mathbb S^2}\int_0^\infty\int_{\mathcal I_\kappa}
	\re^{\ri k(x\cdot\theta+ts)}
	A(k,s)k^3[1-\psi(X k)]\,\rd s\,\rd k\,\rd\theta
	\right|
	\lesssim
	\int_0^{CX^{-1}}k^5\,\rd k
	\lesssim
	X^{-6}.
	\]
	Thus it remains to consider
	\begin{equation}
		\label{eq:LD-k-away-zero-section4}
		k\gtrsim X^{-1}.
	\end{equation}
	
	\medskip
	\noindent
	\emph{Fixed-gap pieces.}
Take a fixed partition of $\mathcal I_\kappa$ into small endpoint
charts and a compact interior set, using one chart across zero
under \eqref{H0}.  Each endpoint chart lies inside the corresponding
chart of Section~\ref{subsec:low-frequency-regular-section3}; its
radius $u_0$ may be decreased further below.  Lemma~\ref{lem:A-endpoint-section2} gives the
smoothness there, so integrations by parts cross zero without an
artificial boundary term.  On the interior set the amplitude is a
smooth fixed-gap symbol.
	
	If \(r\le2t\), then \(t\asymp X\), and repeated integration
	by parts in \(s\) gives, with $M=7$,
	\[
	\left|
	\int_{\mathcal I_\kappa^{\rm int}}
	\re^{\ri tks}A(k,s)\,\rd s
	\right|
	\lesssim_M
	(tk)^{-M}\mathcal A_M(k),
	\]
Here the fixed interior cutoff is included in $A$, and we set
\[
 \mathcal A_M(k):=\max_{0\le j\le M}
 \sup_{s\in\mathcal I_\kappa^{\rm int}}|\partial_s^jA(k,s)|.
\]
For every fixed pair of nonnegative integers $j_1,j_2$,
\begin{equation}\label{eq:fixed-gap-amplitude-symbol-c2}
 \mathcal A_M(k)+
 \sup_{s\in\mathcal I_\kappa^{\rm int}}
 |D_k^{j_1}\partial_s^{j_2}A(k,s)|
 \le C_{M,j_1,j_2}
 \begin{cases}
 k^2,&0<k\le k_\infty,\\
 k^{-2},&k\ge k_\infty.
 \end{cases}
\end{equation}
Here $D_k=k\partial_k$.  Since
$k^{j_1}\partial_k^{j_1}=\prod_{j=0}^{j_1-1}(D_k-j)$, the same bounds give
\[
 \sup_{s\in\mathcal I_\kappa^{\rm int}}
 |\partial_k^{j_1}\partial_s^{j_2} A(k,s)|
 \le C_{j_1,j_2}
 \begin{cases}
 k^{2-j_1},&0<k\le k_\infty,\\
 k^{-2-j_1},&k\ge k_\infty.
 \end{cases}
\]
For $0<k\le k_0$ and $k\ge k_\infty$, respectively, these bounds follow
from Proposition~\ref{prop:low-frequency-factor-section3}(i)--(ii).
On $k_0\le k\le k_\infty$, the open cut away from zero has
$\Bcal\ne0$; under \eqref{H0}, the center is separated by
\eqref{eq:filled-ball-zero-sign-appendix}.  At the endpoints,
$1-\vartheta$ removes the outer contact, whereas $\Acal(\kappa)<0$
separates the inner endpoint denominator when $\kappa>0$.
Lemma~\ref{lem:dispersion-separation-section2} and compactness thus
bound the normalized denominator away from zero on the relevant
closed support.  Its derivatives are bounded in the interior and
satisfy the endpoint conormal estimates on the endpoint charts.
The compact-frequency constant is at most $Ck^2$ after increasing
$C$ by a fixed factor depending on $k_0$.
Using
	\(k\gtrsim X^{-1}\), one obtains
	\[
	\int_{X^{-1}}^{k_{\infty}}
	k^5(tk)^{-M}\,\rd k
	+
	\int_{k_{\infty}}^\infty
	k(tk)^{-M}\,\rd k
	\lesssim
	X^{-6}
	\]
	once \(M>6\).
	
For $r>0$, writing $e=x/r$, the identity
\[
 \int_{\mathbb S^2}\re^{\ri kr\theta\cdot e}\,\rd\theta
 =4\pi\frac{\sin(kr)}{kr}
\]
gives the exact radial formula
\begin{equation}\label{eq:LD-radial-form-section4}
 G_{\rm LD}(t,r)=\frac1{2\pi^2r}\int_0^\infty
 k^2\sin(kr)\int_{\mathcal I_\kappa}\re^{\ri tks}A(k,s)\,\rd s\,\rd k.
\end{equation}
In particular,
\[
 \re^{\ri tks}\sin(kr)
 =\frac{\re^{\ri k(ts+r)}-\re^{\ri k(ts-r)}}{2\ri}.
\]
At $r=0$ we keep the full-sphere formula
\eqref{eq:LD-joint-integral-section4}.  If $r>2t$, then
$r\asymp X$, and the two radial phases satisfy
	\[
	|ts\pm r|\ge r-t\ge r/2
	\]
	for \(s\in\mathcal I_\kappa\).  Repeated integration by parts in
	\(k\), using \eqref{eq:fixed-gap-amplitude-symbol-c2},
	again gives \(O(X^{-6})\).  More explicitly, in the exact
	radial formula the differentiated amplitude is $k^2A(k,s)$.
	Including the lower cutoff, six $k$ derivatives have integral at most
	$C(X+1)$, since their low-frequency integral is bounded by
	$C\int_{c/X}^1 k^{-2}\,\rd k$ and their high-frequency
	integral by $C\int_1^\infty k^{-6}\,\rd k$.  Thus the bound is
	$Cr^{-1}X^{-6}(X+1)\lesssim X^{-6}$.
	The upper truncation annulus contributes $O(H^{-5})$ after these
	derivatives and hence vanishes.
	Consequently,
	\begin{equation}
		\label{eq:LD-fixed-gap-section4}
		|G_{\rm LD}^{\rm int}(t,x)|
		\lesssim
		X^{-6}.
	\end{equation}
	
	In an endpoint chart the amplitude is only conormal in the
	\(s\)-variable, so the preceding \(s\)-integration-by-parts argument
	is no longer available.
	{
	Put $a_{\min}=\kappa$ for the annulus and $a_{\min}=1$ for the
	filled-ball variant, which has only the outer endpoint charts.
	Shrink the fixed endpoint charts so that $|s|\ge a_{\min}/2$ on
	them.  In the interior $r\le a_{\min}t/4$, keep the spherical
	integral intact: its phase satisfies
	\[
	 |x\cdot\theta+ts|\ge t|s|-r\ge a_{\min}t/4
	 \qquad(\theta\in\mathbb S^2).
	\]
	Here $t\asymp X$.  The conormal $s$ bounds and the three frequency
regimes used in \eqref{eq:fixed-gap-amplitude-symbol-c2} give
	\[
	 \int_{\mathcal I_\kappa}\int_0^\infty
	 \left|\partial_k^7\bigl[\psi(X k)k^3A(k,s)\bigr]\right|
	 \,\rd k\,\rd s
	 \le C\left(\int_{c/X}^1 k^{-2}\,\rd k
	             +\int_1^\infty k^{-6}\,\rd k\right)
	 \le CX.
	\]
	Derivatives of the lower cutoff are supported at $k\asymp
	X^{-1}$ and have exactly the displayed bound.  The upper
	truncation terms are $O(H^{-5})$.  Seven integrations in $k$
	therefore give $Ct^{-7}X\lesssim X^{-6}$,
	uniformly down to $r=0$.
	}
	For the remaining points $r\ge a_{\min}t/4$ one has
	$r\asymp X$ when $X$ is large.  Away from the
	corresponding front $r=at$, the radial $k$-phase is uniformly
	nonstationary.  Indeed, writing
	\[
	s=a+\epsilon_a u,\qquad 0<u<u_0,
	\]
	the two phases arising from \(\sin(kr)\) are
	\[
	k(at+r+\epsilon_a tu),
	\qquad
	k(at-r+\epsilon_a tu).
	\]
	The first is always nonstationary, while, if
	\[
	|r-at|\ge c_0X,
	\]
	then, after choosing \(u_0>0\) sufficiently small,
	\[
	|at-r+\epsilon_a tu|
	\ge cX
	\]
	uniformly on the endpoint chart.  Repeated integration by parts in
	\(k\), using the same endpoint symbol bounds in the three
frequency regimes of \eqref{eq:fixed-gap-amplitude-symbol-c2}, gives
	\(O(X^{-6})\): the six-derivative estimate for $k^2A$
	just proved applies after integration of its $u^{\gamma_a}$ factor.
	The negative endpoint charts give the same bounds after $s\mapsto-s$.
	Thus only a fixed conic neighborhood of
	\(r=at\) requires the endpoint transform.
	
	\medskip
	\noindent
\emph{Endpoint reduction and the fronts.}
The endpoint estimate consists of two Fourier transfers.  Vanishing
of order $\beta$ in the inward variable $u$ first produces the decay
$(tk)^{-\beta-1}$.  Radial inversion then leaves the power
$k^{3-\beta}=k^{(4-\beta)-1}$, whose Fourier transform measures the
distance to the front.  The two steps are supplied by
Lemmas~\ref{lem:finite-endpoint-v4} and~\ref{lem:truncated-radial-v4}.
The value $\beta=4$ leaves the borderline kernel $1/k$; its treatment
will use the pairing of the reflected endpoints before taking
absolute values.

Let $a\in\{\kappa,1\}$ be a positive endpoint, with only $a=1$ in
the ball case, and write
\[
 s=a+\epsilon_a u,\qquad u>0,\qquad
 \epsilon_\kappa=1,\qquad\epsilon_1=-1.
\]
Choose a smooth inward cutoff $\chi_{\rm ep}$ supported in
$[0,u_0)$ and equal to one near zero.  On this chart set
$\chi_a(s):=\chi_{\rm ep}(\epsilon_a(s-a))$ and
\[
 w(k):=k^2\ang k^{-4},\qquad
 f_a(k,u):=w(k)^{-1}A(k,a+\epsilon_a u)\chi_{\rm ep}(u),
 \qquad
 \beta:=\min\{\gamma_a,9/2\}.
\]
Since $1-\vartheta$ removes the contact band,
Proposition~\ref{prop:low-frequency-factor-section3}(i)--(ii),
together with the compact-frequency denominator bound used in
\eqref{eq:fixed-gap-amplitude-symbol-c2}, gives
\begin{equation}\label{eq:green-whole-conormal-v3}
 |D_k^{j_1}\partial_u^{j_2} f_a(k,u)|\le C_{j_1,j_2}u^{\beta-j_2},
 \qquad j_1+j_2\le10.
\end{equation}
The weight $w$ incorporates both limiting frequency orders, and the
extra endpoint power in the power--log germ absorbs its logarithms.
The order ten has a specific role: the endpoint transform uses
$L=6>\beta+1$ derivatives, and four derivatives along $tk$ may all
transfer to $u$, giving $L+4=10$.  These are derivatives of the
fixed equilibrium; this count does not change the data index $M_*$.

With $T_f$ defined in Lemma~\ref{lem:finite-endpoint-v4}, set
\begin{equation}\label{eq:green-whole-transform-v3}
 q_{a,t}(k):=T_{f_a}(k,tk)
 =(tk)^{\beta+1}\int_0^{u_0}
        \re^{\ri\epsilon_a tku}f_a(k,u)\,\rd u.
\end{equation}
The ray estimate \eqref{eq:endpoint-ray-symbol-v196} gives directly
\[
 |D_k^j q_{a,t}(k)|\le C_j\min\{1,(tk)^{\beta+1}\},
 \qquad 0\le j\le4.
\]
Thus the whole endpoint contribution, including its remainder, is
\[
 I_a(t,k):=\int \re^{\ri tks}A(k,s)\chi_a(s)\,\rd s
 =\re^{\ri atk}w(k)(tk)^{-\beta-1}q_{a,t}(k).
\]
The factor $\re^{\ri atk}$ identifies the endpoint speed, while
$q_{a,t}$ remains a controlled radial symbol even on the small-$tk$
interval added back in Lemma~\ref{lem:truncated-radial-v4}.

We now work in the remaining cone $|r-at|\le c_0X$, where
$r\asymp t\asymp X$ after fixing $c_0>0$ sufficiently small in terms
of the endpoint speeds.  Thus $t/X$ has a fixed positive lower bound
there.  Choose $c$ in $\psi(Xk)$ so that $(t/X)c\ge1$ throughout this
cone; increasing $c$ only adds an annulus already bounded by $CX^{-6}$.
Substituting $I_a$ into \eqref{eq:LD-radial-form-section4} gives
\[
 \frac1{2\pi^2rt^{\beta+1}}\int_0^\infty
 \psi(Xk)\sin(kr)\re^{\ri atk}
 k^{3-\beta}\ang k^{-4}q_{a,t}(k)\,\rd k,
 \qquad
 k^2\,k^2\,k^{-\beta-1}=k^{3-\beta}.
\]
After expansion of the sine, each carrier is therefore bounded by
\begin{equation}\label{eq:LD-front-reduction-section4}
 CX^{-\beta-2}\left|
 \int_0^\infty \re^{-\ri k\ell}\psi(X k)
       k^{3-\beta}\ang k^{-4}q_{a,t}(k)\,\rd k\right|,
 \qquad \ell:=r-at\ \hbox{or }-(r+at).
\end{equation}
Here Lemma~\ref{lem:truncated-radial-v4} applies to
$\ang k^{-4}q_{a,t}$ with $\alpha=4-\beta$, $J=4$, and $B=4$.
If $1\le\gamma_a<4$, its positive-$\alpha$ branch gives
\[
 CX^{-\gamma_a-2}\ang\ell^{-(4-\gamma_a)}
 +CX^{-6}.
\]
If $\gamma_a>4$, then $\beta>4$, and the negative-$\alpha$ branch
gives
$CX^{-\beta-2}X^{\beta-4}
=CX^{-6}$ by absolute integration.
The carrier with $r+at\asymp X$ is always
$O(X^{-6})$, including the critical branch $\beta=4$.
The high-frequency tails have integrable order $k^{-1-\beta}$,
with the same symbol bounds through the four required derivatives.
Changing the endpoint cutoff adds only fixed-gap terms already
covered by \eqref{eq:LD-fixed-gap-section4}.

\medskip
\noindent
\emph{The fourth-order endpoint.}
For $\gamma_a=4$, a separate estimate of either endpoint in
\eqref{eq:LD-front-reduction-section4} would leave
$X^{-6}\log(2+X/\ang{r-at})$.
We remove this logarithm by first pairing the two physical endpoint
contributions.  Choose mirrored cutoffs at $a$ and $-a$.
The identities $A(k,-s)=-A(k,s)$ and $A(k,s)\in\ri\R$ will
pair their leading coefficients.  We first obtain the required
pure-power remainder directly from the jump formula.
Write $\Bcal(a+\epsilon_a u)=u^4\beta_a(u)$ with $\beta_a$ smooth.
Lemma~\ref{lem:A-endpoint-section2} gives
$\Acal(a+\epsilon_a u)-\Acal(a)=O_{\rm con}(u)$ through the required
finite orders: its fourth-order term $u^4\log u$ satisfies this
first-increment bound.  On the support of $1-\vartheta$, put
\[
 \widetilde d(k,u):=\ang k^{-4}
 \left[(k^2-4\pi\Acal(a+\epsilon_a u))^2
             +(2\pi\Bcal(a+\epsilon_a u))^2\right].
\]
The low-frequency bound \eqref{eq:D0-lower-proof-section3}, the
compact-frequency separation used in
\eqref{eq:fixed-gap-amplitude-symbol-c2}, and the high-frequency
normalization give $\widetilde d\ge c>0$ on a fixed endpoint chart.
Moreover, the product and reciprocal rules give
\[
 D_k^{j_1}D_u^{j_2}[\widetilde d(k,u)^{-1}
                 -\widetilde d(k,0)^{-1}]=O(u),\qquad j_1+j_2\le7.
\]
Indeed, the same bound first holds for the denominator increment;
all its mixed derivatives are bounded, and the reciprocal stays
uniformly bounded.  The explicit jump \eqref{eq:J-explicit-proof-section3}
therefore yields
\[
 \frac{A(k,a+\epsilon_a u)}{w(k)u^4}
 =\frac{2\ri(1-\vartheta(k))\beta_a(0)}{\widetilde d(k,0)}
       +O_{\rm con}(u)
\]
with the same mixed orders.  Derivatives of $\vartheta$ lie in a
fixed separated frequency band; derivatives of the endpoint cutoff
are supported away from zero.  Including both cutoffs and converting
conormal derivatives by the conventions at the end of
Section~\ref{sec:introduction}, we obtain
\begin{equation}\label{eq:green-fourth-germ-v4}
 \begin{aligned}
 A(k,a+\epsilon_a u)\chi_{\rm ep}(u)
   &=\ri c_a(k)u^4\chi_{\rm ep}(u)+R_a(k,u),\\
 |D_k^{j_1}\partial_u^jR_a(k,u)|
   &\le C_{j_1,j}w(k)u^{5-j},\qquad j_1+j\le7.
 \end{aligned}
\end{equation}
The coefficient
\[
 b_a(k):=\frac{c_a(k)}{k^2}
 =\frac{2\beta_a(0)(1-\vartheta(k))}{[k^2-4\pi\Acal(a)]^2}
\]
is real, smooth and even at zero, and has order $-4$ at infinity;
$\Acal(a)$ is the real boundary function in \eqref{eq:A-B-section2}.
At the outer endpoint this expression is extended by zero on the
fixed neighborhood where $\vartheta=1$.  At the inner endpoint
$\Acal(\kappa)<0$, so its denominator never vanishes.

Five endpoint integrations by parts give
\[
 \int_0^\infty \re^{\ri\epsilon_a\zeta u}u^4\chi_{\rm ep}(u)\,\rd u
 =24\ri\epsilon_a\zeta^{-5}+O_N(\zeta^{-N}),\qquad \zeta\ge1.
\]
The leading transforms at $a$ and $-a$ are consequently
$-24\epsilon_a c_a(k)(tk)^{-5}\re^{\ri atk}$ and
$-24\epsilon_a c_a(k)(tk)^{-5}\re^{-\ri atk}$.
Apply Lemma~\ref{lem:finite-endpoint-v4} to $R_a/w$, with exponent
five and $L=7$.  Including the rapidly decaying
cutoff remainder above, the paired transform is
\[
 -48\epsilon_ac_a(k)(tk)^{-5}\cos(atk)
       +O\bigl(w(k)(tk)^{-6}\bigr).
\]
We retain the same lower cutoff, so that $tk\ge1$.
Multiplication by the radial sine now uses the identity
\[
 \sin(kr)\cos(atk)
 =\tfrac12\bigl[\sin(k(r-at))+\sin(k(r+at))\bigr].
\]
Hence \eqref{eq:LD-radial-form-section4} makes the leading term a
fixed constant multiple of
\begin{equation}\label{eq:green-fourth-sine-v4}
 \frac1{rt^5}\int_0^\infty\psi(Xk)b_a(k)
 \frac{\sin(k(r-at))+\sin(k(r+at))}{k}\,\rd k.
\end{equation}
Each integral is uniformly bounded.  Indeed, for a fixed smooth
compact cutoff $\chi_*=1$ near zero,
$(b_a-b_a(0)\chi_*)/k\in L^1$.  Its contribution is bounded
absolutely, and the remaining sine integral is bounded by the
Dirichlet integral.  Removing $k\lesssim X^{-1}$ costs at
most $C|r\pm at|/X\le C$.
Since $r\asymp t\asymp X$, the leading term is
$CX^{-6}$.  This bound comes from the paired sine-over-$k$
integral, before any separate absolute estimate of the two endpoints.

The remainder has one additional endpoint power and can be
integrated absolutely:
\[
 \frac C{rt^6}\int_{c/X}^\infty
 k^{-2}\ang k^{-4}\,\rd k
 \le \frac{CX}{rt^6}\le CX^{-6}.
\]
The low-frequency, nonfront, interior and exterior regions were already
bounded by $CX^{-6}$, so no logarithm remains for this
fourth-order Green endpoint.  Only fixed-equilibrium derivatives
have entered this calculation.

For $1\le\gamma_a<4$, the front contribution just obtained is
$X^{-\gamma_a-2}\ang{r-at}^{\gamma_a-4}
=X^{-6}(X/\ang{r-at})^{4-\gamma_a}$; for $\gamma_a\ge4$ it is
$O(X^{-6})$.  Summing over the positive endpoints and their
Hermitian reflections proves \eqref{eq:green-LD-outline}.
\end{proof}

\subsection{Uniform estimates near the survival threshold}
\label{subsec:threshold-core-v4}
We start from the physical spectral measure constructed in
Section~\ref{sec:spectral-decomposition}.  Proposition~\ref{prop:terminal-reconstruction-section3}
gives
\begin{equation}\label{eq:threshold-measure-inverse-section4}
 \widehat G_{\rm th}(t,k)
 =\vartheta(k)\int_\R \re^{\ri t\tau}\,\rd\mathfrak M_k(\tau).
\end{equation}
For $r=|x|>0$, its radial inverse is therefore
\begin{equation}\label{eq:threshold-radial-measure-v18}
 G_{\rm th}(t,r)=\frac1{2\pi^2r}\int_0^\infty k\sin(kr)\vartheta(k)
       \int_\R \re^{\ri t\tau}\,\rd\mathfrak M_k(\tau)\,\rd k.
\end{equation}

\begin{proposition}
\label{prop:threshold-bound-outline}
For every integer $N\ge0$,
\begin{equation}\label{eq:threshold-global-outline}
 |G_{\rm th}(t,x)|\le C_N\ang t^{-3/2}
       \frac{\ang t^N}{\ang{t,x}^N},\qquad t>0.
\end{equation}
\end{proposition}

The proof, completed at the end of this subsection, first establishes
the unweighted bound in $r\le2t$ for the contact core and its ordinary
complement.  Finite propagation of the complete response then
recovers any fixed spatial weight.

\paragraph{The spectral core and its complement.}
Keep the physical cutoff $\vartheta$ fixed.  Choose a real-valued
$\zeta\in C_c^\infty((0,\infty))$, supported in $\{\vartheta=1\}$
and equal to one near
$k_*=\sqrt{4\pi\Acal(1)}$.  For a real-valued cutoff $\chi_y$, equal to
one near zero and with its small support fixed below, put
\[
 \chi_+(k,\tau)=\chi_y(1-\tau/k),\qquad
 \chi_-(k,\tau)=\chi_+(k,-\tau),\qquad
 \chi_{\rm tr}=\chi_++\chi_-.
\]
For each $k\in\supp\zeta$, the supports of $\chi_+(k,\cdot)$ and $\chi_-(k,\cdot)$ are disjoint. Define
\begin{equation}\label{eq:threshold-core-split-v18}
\begin{split}
 \widehat G_{\rm core}(t,k)
 &=\zeta(k)\int_\R \re^{\ri t\tau}\chi_{\rm tr}(k,\tau)
                         \,\rd\mathfrak M_k(\tau),\\
 \widehat G_{\rm rem}(t,k)
 &=(\vartheta-\zeta)(k)\int_\R \re^{\ri t\tau}\,\rd\mathfrak M_k(\tau)
   +\zeta(k)\int_\R \re^{\ri t\tau}(1-\chi_{\rm tr}(k,\tau))
                         \,\rd\mathfrak M_k(\tau).
\end{split}
\end{equation}
Thus $G_{\rm th}=G_{\rm core}+G_{\rm rem}$.
The core contains the pole--continuum contact; the complement can
still contain ordinary poles, whose dispersion is $t^{-3/2}$.
Fourier inversion, the uniform total-variation bound in
Proposition~\ref{prop:terminal-measure-section3}, and
\eqref{eq:full-physical-measure-section3} give
\begin{equation}\label{eq:threshold-core-bounded-time-v196}
 \sup_{x\in\R^3}|G_{\rm core}(t,x)|
 \le C\int_0^\infty k^2|\zeta(k)|
                \|\mathfrak M_k\|_{\mathrm{TV}}\,\rd k
 \le C,\qquad t\ge0.
\end{equation}
For any fixed $T_0\ge2$, this gives the required
$C_{T_0}t^{-3/2}$ bound on $2\le t\le T_0$.
For large times we estimate the same core in the conic region
$v_0t\le r\le2t$ by radial oscillation and damping, and in the interior
region $0\le r\le v_0t$ by its common holomorphic boundary value,
retaining the angular integral; $v_0>0$ is fixed below.
The common branch and transverse inversion provide both representations.
The next two preparations isolate their oscillatory and holomorphic
parts, to be estimated in Lemma~\ref{lem:terminal-inversion-appendix}.

\paragraph{The common holomorphic branch.}
We first construct the common holomorphic branch.  Its boundary
trace determines the transverse roots, and its positive imaginary
part controls the interior-region contour estimate.  The variable $z$
is an auxiliary ratio coordinate, distinct from the Laplace variable
in Sections~\ref{sec:dispersion}--\ref{sec:spectral-decomposition}.

{
Put $F(u)=\mu(u^2/2)$.  On a small upper chart with positive real part,
\[
 \mathcal U_{\varepsilon_{\rm tr}}
 :=\{z=\xi+\ri\eta:|\xi-1|<2\varepsilon_{\rm tr},\quad
                         0<\eta<2\varepsilon_{\rm tr}\},
\]
define
\begin{equation}\label{eq:Omega-holomorphic-v4}
 H(z):=4\pi\int_\kappa^1\frac{u^2F(u)}{1-z^2u^2}\,\rd u.
\end{equation}
On compact subsets of the chart the denominator stays away from
zero, so differentiation under the integral proves holomorphicity.
Near $z=1$, $H$ satisfies the following continuity estimate:
\begin{equation}\label{eq:H-endpoint-modulus-c2}
\begin{gathered}
 |H(z)-H(1)|\le C|z-1|
 \left(1+\one_{\{\gamma_1=1\}}\log\frac1{|z-1|}\right),\quad
 0<|z-1|<d_0<1,\quad
 z\in\overline{\mathcal U_{\varepsilon_{\rm tr}}},
\end{gathered}
\end{equation}
where $H(1)=\omega_*^2>0$, with $\omega_*$ defined in
\eqref{eq:terminal-frequency-section2}.  To prove the modulus, change
variables $u^{-1}=1+s$ near the outer endpoint and use partial
fractions.  Up to a smooth holomorphic term, $H$ is the Cauchy
transform of $s^{\gamma_1}b(s)$ with $b$ smooth.  Writing $w=z-1$,
its difference from the endpoint value is
\[
 \int_0^a\frac{s^{\gamma_1}b(s)}{s-w}\,\rd s
 -\int_0^a s^{\gamma_1-1}b(s)\,\rd s
 =w\int_0^a\frac{s^{\gamma_1-1}b(s)}{s-w}\,\rd s.
\]
For $\gamma_1=1$, subtract $b(\Re w)$ in the last numerator,
using a smooth extension across zero.  The difference is uniformly
integrable, and the constant part is
$b(\Re w)[\log(a-w)-\log(-w)]$.
For $\gamma_1>1$, the same subtraction uses the positive H\"older
regularity of $s^{\gamma_1-1}b(s)$ and bounds the last integral
uniformly.  This proves \eqref{eq:H-endpoint-modulus-c2}, including
the continuous boundary trace at~$1$.

Positivity on a fixed compact subinterval of $(\kappa,1)$ gives
\begin{equation}\label{eq:holomorphic-sector-bound-v4}
 \Im H(\xi+\ri\eta)=8\pi\xi\eta\int_\kappa^1
       \frac{u^4F(u)}{|1-(\xi+\ri\eta)^2u^2|^2}\,\rd u\ge c\eta.
\end{equation}
By \eqref{eq:H-endpoint-modulus-c2}, choose
$\mathcal U_{\varepsilon_{\rm tr}}$ small enough that
$|H-\omega_*^2|<\omega_*^2/2$.  Then
\[
 \Omega^\sharp(z):=\sqrt{H(z)},\qquad \Re\Omega^\sharp>0,
\]
is single valued and holomorphic, with
\[
 c\le\Re\Omega^\sharp\le C_1,\qquad
 \Im\Omega^\sharp
 =\frac{\Im H}{2\Re\Omega^\sharp}\ge c_1\eta>0.
\]
The Plemelj formula and \eqref{eq:A-B-section2} identify its upper
boundary trace:
\begin{equation}\label{eq:holomorphic-trace-v196}
 \begin{aligned}
 H(\rho+\ri0)
 &=\rho^{-2}\bigl[4\pi\Acal(\rho^{-1})
                   +2\pi\ri\Bcal(\rho^{-1})\bigr]
   =\rho^{-2}\Om(\rho^{-1}),\\
 \Omega^\sharp(\rho+\ri0)&=\Om_0(\rho).
 \end{aligned}
\end{equation}
Indeed, taking the reciprocal sends approach from the upper half-plane
to approach from below to $\rho^{-1}$.  Since $\rho>0$ and both square
roots have positive real part, the second identity follows from
\eqref{eq:Omega0-definition-section2}; continuity extends it through
$\rho=1$, also when $\gamma_1=1$.

We now fix the local supports and constants.  Choose
$\chi_y\in C_c^\infty((-\varepsilon_{\rm tr},\varepsilon_{\rm tr}))$
and an upper bound for $\epsilon>0$ so that
\[
 \rho_\epsilon(y):=(1-y-\ri\epsilon)^{-1}
       \in\mathcal U_{\varepsilon_{\rm tr}},\qquad
 q_\epsilon(y):=\rho_\epsilon(y)
                      \Omega^\sharp(\rho_\epsilon(y))
\]
for $y\in\supp\chi_y$, with $\Re q_\epsilon$ in a fixed interval
compactly contained in $\{\zeta=1\}$.  The chart is also chosen so
that $\Re(z\Omega^\sharp(z))>0$.  The one-sided geometric and damping
bounds used below hold on these supports.  In particular, if
$1<\gamma_1\le2$, choose a fixed
$0<\epsilon_v<(1-v_*)/4$, with $v_*$ from
\eqref{eq:vstar-section2}, and require on the continuum support that
\[
 |V(\rho)-v_*|<\epsilon_v/2,\qquad
 1-v\rho\ge c>0\quad\hbox{when }|v-v_*|<\epsilon_v.
\]
These choices are possible by the one-sided limits in
Proposition~\ref{prop:terminal-jet-section2}.
Set
\[
 X_0:=\sup_{z\in\mathcal U_{\varepsilon_{\rm tr}}}|\Re z|
       =1+2\varepsilon_{\rm tr},
\]
and choose $v_0,T_0$ so that
\[
 0<v_0\le\min\{(2X_0)^{-1},c_1/(4C_1)\},\qquad
 T_0\ge\max\{2,2X_0\}.
\]
All these choices, including the support of $\chi_y$, are fixed
independently of $t,x$; no further shrinking is used below.

\paragraph{Transverse inversion of the joint spectral measure.}
We derive an exact radial representation of $G_{\rm core}$ by
integrating in $k$ with positive regularization parameter and then
taking the boundary limit.  The Cauchy representation
\eqref{eq:analytic-cauchy-section2} gives
\begin{equation}\label{eq:causal-transverse-root-v196}
 \begin{aligned}
 q_\epsilon(y)^2
 =\rho_\epsilon(y)^2H(\rho_\epsilon(y))
   =4\pi\int_\kappa^1
       \frac{u^2F(u)}{(1-y-\ri\epsilon)^2-u^2}\,\rd u
 =k^2\Lcal(k\epsilon+\ri k(1-y),k).
 \end{aligned}
\end{equation}
Thus the roots come from the original causal Laplace half-plane
$\Re z_{\rm Lap}=k\epsilon>0$.  With $s=1-y>0$, the imaginary part is
\[
 \Im q_\epsilon(y)^2
 =8\pi s\epsilon\int_\kappa^1
       \frac{u^2F(u)}{|(s-\ri\epsilon)^2-u^2|^2}\,\rd u>0.
\]
The positive real-part branch therefore has $\Im q_\epsilon>0$;
its conjugate lies in the lower half-plane.  Using
$Y$ from \eqref{eq:terminal-XY-section3}, define the boundary root
\begin{equation}\label{eq:transverse-roots-v18}
 q_+(y):=\lim_{\epsilon\downarrow0}q_\epsilon(y)
       =\sqrt{4\pi\Acal(1-y)+\ri Y(y)},\qquad \Re q_+>0.
\end{equation}
The local joint measure is represented in the limiting-absorption sense by
\begin{equation}\label{eq:terminal-joint-measure-v18}
 \rd k\,\rd\mathfrak m_k^+(y)
 =\lim_{\epsilon\downarrow0}\frac{k^3}{2\pi}
 \left(\frac1{k^2-q_\epsilon(y)^2}
       -\frac1{k^2-\overline{q_\epsilon(y)}^{\,2}}\right)
       \,\rd k\,\rd y,
\end{equation}
locally in distributions on the fixed positive $k$-band and the
$y$-interval.  For $y>0$, the Lorentzian identity gives the continuum
density in \eqref{eq:terminal-cut-weight-section3}.  For $y<0$,
Plemelj gives $2\pi\ri\delta(k^2-4\pi\Acal(1-y))$, hence the pole
measure \eqref{eq:terminal-pole-delta-section3} after multiplication
by $k^3/(2\pi)$.  To pass to the joint limit, write
$A_\epsilon=\Re(q_\epsilon^2)$ and $B_\epsilon=\Im(q_\epsilon^2)>0$.
The absolute density is
$k^3B_\epsilon/[\pi((k^2-A_\epsilon)^2+B_\epsilon^2)]$.
On the fixed positive $k$-band, $X=k^2-A_\epsilon$ and
$\rd X=2k\,\rd k$ bound its $k$-integral uniformly in $\epsilon,y$.
Thus the one-sided weak limits in $k$ may be integrated in $y$ by
dominated convergence after the $k$-integration.  In particular,
a strip $|y|<\eta$ has joint mass at most $C\eta$ and creates no
additional mass on $y=0$.  The opposite pole-side boundary
prescriptions remain encoded by $q_\epsilon$ and
$\overline{q_\epsilon}$ until the transverse integration is performed.
The fixed critical fiber remains the weak-limit measure of
Proposition~\ref{prop:terminal-measure-section3}, including its atom
when $\gamma_1>1$; its single $k$-value does not alter the subsequent
spatial Fourier integral.

Recall that $y=1-\tau/k$, so $\tau=k(1-y)$.
For fixed $t$ and $r>0$, apply
\eqref{eq:terminal-joint-measure-v18} in the radial inversion
\eqref{eq:threshold-radial-measure-v18}.  Before taking the limit,
Fubini's theorem applies by the uniform $k$-integral bound above.
The measure's factor $k^3$ and the radial factor $k$ give $k^4$;
the Jacobian of $\tau\mapsto y$ is already included.
Expanding $\sin(kr)$ and integrating in $k$ with the real cutoff
$\zeta$ gives
\begin{equation}\label{eq:terminal-transverse-regularized-v199}
 \begin{split}
 H_+(t,r)=\lim_{\epsilon\downarrow0}\frac1{4\pi^2\ri r}
 \int\chi_y(y)\bigl[\mathcal T_4(t(1-y)+r,q_\epsilon(y))-\mathcal T_4(t(1-y)-r,q_\epsilon(y))\bigr]\,\rd y.
 \end{split}
\end{equation}
For each fixed Fourier parameter, the continuous boundary value
\eqref{eq:transverse-real-boundary-v199} gives convergence of the
transverse transform, including the real roots on the pole side.
The bound \eqref{eq:transverse-small-parameter-v18} is uniform for
all real parameters, so $2C|\chi_y(y)|$ dominates this integrand.
Dominated convergence in $y$ therefore yields
\begin{equation}\label{eq:terminal-transverse-radial-v18}
 H_+(t,r)=\frac1{4\pi^2\ri r}\int\chi_y(y)
       [\mathcal T_4(t(1-y)+r,q_+(y))
        -\mathcal T_4(t(1-y)-r,q_+(y))]\,\rd y,
 \qquad G_{\rm core}=2\Re H_+.
\end{equation}
Here the real mirror cutoffs and the Hermitian reflection
\eqref{eq:terminal-negative-section3} give
$H_-=\overline{H_+}$.  Put $\rho:=(1-y)^{-1}$,
$\chi_\rho(\rho):=\chi_y(1-\rho^{-1})$, and
$L_\pm(\rho):=t/\rho\pm r$.  Since
$\rd y=\rho^{-2}\rd\rho$ and $q_+=\rho\Om_0(\rho)$,
\begin{equation}\label{eq:terminal-transverse-rho-v196}
 H_+(t,r)=\frac1{4\pi^2\ri r}\int\chi_\rho(\rho)\rho^{-2}
 \bigl[\mathcal T_4(L_+,\rho\Om_0(\rho))
       -\mathcal T_4(L_-,\rho\Om_0(\rho))\bigr]\,\rd\rho.
\end{equation}

\paragraph{The small-parameter region and the conic phases.}
In $v_0t\le r\le2t$, we first remove the small-$L_-$ region and the
transverse errors; both contribute $O(t^{-2})$.  The remaining
terms are the incoming and outgoing oscillatory integrals below.
Set $v=r/t$ and fix a smooth real function $\eta_0$ equal to one
on $[-1,1]$ and zero outside $[-2,2]$.
For $t\ge T_0$ and $v_0t\le r\le2t$, the part of
$H_+$ in \eqref{eq:terminal-transverse-rho-v196} with $\eta_0(L_-)$
is bounded by \eqref{eq:transverse-small-parameter-v18}.
Since $\partial_\rho L_-=-t/\rho^2$ and $\rho\asymp1$ on
$\supp\chi_\rho$, this region has length $O(t^{-1})$.  Its contribution
to $H_+$, and hence to $G_{\rm core}$, is $O((rt)^{-1})=O(t^{-2})$.
The large-parameter errors from Lemma~\ref{lem:transverse-cauchy-v18}
are bounded by
\begin{equation}\label{eq:transverse-layer-error-v18}
 \frac Cr\int_{\substack{\rho\in\supp\chi_\rho\\|L_-|\ge1}}
       |L_-|^{-2}\,\rd\rho
 +\frac Cr\int_{\supp\chi_\rho}\ang{L_+}^{-2}\,\rd\rho
 \le Ct^{-2}.
\end{equation}
For the first integral use the change of variables $L_-$;
for the second, $L_+\ge t/X_0\gtrsim t$ on the support and
$t\ge T_0$ puts it in the large-parameter range.

For the leading transverse term, the factor $\ri$ in
\eqref{eq:transverse-residues-v18} cancels the radial factor $1/\ri$.
Including the Jacobian, define
\begin{equation}\label{eq:exact-terminal-residues-v3}
 \mathcal R_4(\rho):=\rho^{-2}\frac{(\rho\Om_0(\rho))^3}{2}
                   =\frac{\rho\Om_0(\rho)^3}{2}.
\end{equation}
This amplitude retains the complete complex value of $\Om_0$.
Write $\Om_+=\Om_0$, $\Om_-=\overline{\Om_0}$ and
$\mathcal R_4^+=\mathcal R_4$, $\mathcal R_4^-=\overline{\mathcal R_4}$.
With $\sigma=\operatorname{sgn}L_-$, the conic formula becomes
\begin{equation}\label{eq:terminal-conic-residues-v18}
\begin{split}
 G_{\rm core}(t,r)=\frac1{2\pi^2r}\Re\int\chi_\rho(\rho)
 \big[\mathcal R_4\re^{\ri(t+r\rho)\Om_0}
  -(1-\eta_0(L_-))\mathcal R_4^\sigma
                 \re^{\ri(t-r\rho)\Om_\sigma}\big]\,\rd\rho
 +O(t^{-2}).
\end{split}
\end{equation}
The two terms arise, respectively, from $\re^{\ri kr}$ and
$\re^{-\ri kr}$ in the radial Fourier kernel; we call them incoming
and outgoing.  Recall from \eqref{eq:Omega0-definition-section2}
and \eqref{eq:q-V-section2} that
$\Om_0=\omega+\ri\nu$ and
$V=\omega'/(\omega+\rho\omega')$ on either side of contact.
The outgoing real phase and absolute damping are
\begin{equation}\label{eq:terminal-real-phase-appendix}
 \phi_v(\rho):=(1-v\rho)\omega(\rho),\qquad
 \left|\re^{\ri(t-r\rho)\Om_\sigma(\rho)}\right|
       =\re^{-t|1-v\rho|\nu(\rho)}.
\end{equation}
Indeed, $\sigma=\operatorname{sgn}(1-v\rho)$ on the support of
the large-parameter cutoff.  Differentiating gives
\begin{equation}\label{eq:terminal-phase-derivative-appendix}
 \phi_v'=g(V-v),\qquad g:=\omega+\rho\omega'.
\end{equation}
The incoming real phase $(1+v\rho)\omega$ has derivative
$g(V+v)$.  Thus the group-velocity geometry controls the real
oscillation, while the root selected by $\operatorname{sgn}L_-$
gives nonnegative absolute damping.

\paragraph{Holomorphic phase estimates in the interior region.}
For $0\le r\le v_0t$, we retain the spherical representation of the
same inverse transform.  Its leading term uses the common boundary
trace \eqref{eq:holomorphic-trace-v196} across the contact point.
The following contour estimate will apply before the angular
integration.  For $|\lambda|\le v_0$, define
\[
 \Psi_\lambda^\sharp(z):=(1+\lambda z)\Omega^\sharp(z),
 \qquad z\in\mathcal U_{\varepsilon_{\rm tr}}.
\]
In the interior-region application $\lambda=x\cdot\theta/t$.
For both signs of $\lambda$, the fixed choices of $X_0,v_0$ give
\[
 \Im\Psi_\lambda^\sharp(\xi+\ri\eta)
 =(1+\lambda\xi)\Im\Omega^\sharp
       +\lambda\eta\Re\Omega^\sharp
 \ge(c_1/2-v_0C_1)\eta\ge c_1\eta/4.
\]
Hence
\begin{equation}\label{eq:deep-timelike-phase-lower-appendix}
 \Im\Psi_\lambda^\sharp(z)\ge c\Im z,
 \qquad z\in\mathcal U_{\varepsilon_{\rm tr}},\quad |\lambda|\le v_0.
\end{equation}
The same choice yields
$1/\rho+\lambda\ge(2X_0)^{-1}$ on $\supp\chi_\rho$.
The chart and constants are independent of $t,x$ and the spherical
direction, including when $\gamma_1=1$.

The next estimate uses a finite-order almost-analytic extension,
following the standard construction in \cite{Zwo12}.

\begin{lemma}
\label{lem:holomorphic-terminal-v4}
Let $I_{\rm tr}:=(1-2\varepsilon_{\rm tr},1+2\varepsilon_{\rm tr})$
be the real boundary interval of the upper chart
$\mathcal U_{\varepsilon_{\rm tr}}$, and let
$\Psi_\lambda^\sharp$ be the phase in
\eqref{eq:deep-timelike-phase-lower-appendix}.
Suppose that $A^\sharp$ is bounded and holomorphic on this chart and
has a continuous, possibly complex-valued boundary trace on $I_{\rm tr}$.
For $\chi\in C_c^\infty(I_{\rm tr})$ and every fixed integer $J\ge0$,
\begin{equation}\label{eq:holomorphic-cancellation-v4}
 \left|\int_{I_{\rm tr}}
 \re^{\ri t\Psi_\lambda^\sharp(\rho+\ri0)}
 A^\sharp(\rho+\ri0)\chi(\rho)\,\rd\rho\right|
 \le C_Jt^{-J-1},\qquad t\ge1,\quad |\lambda|\le v_0.
\end{equation}
The constant depends on the fixed chart, the uniform lower bound
in \eqref{eq:deep-timelike-phase-lower-appendix},
$\|A^\sharp\|_\infty$, and the derivative bounds for $\chi$ through order $J+1$.
\end{lemma}

\begin{proof}
Write $z=\xi+\ri\eta$ and set
$\bar\partial:=\tfrac12(\partial_\xi+\ri\partial_\eta)$.
Let $\chi_{\rm v}$ be a smooth vertical cutoff, equal to one for
$0\le\eta\le\varepsilon_{\rm tr}/2$ and vanishing for $\eta$ near
and above $\varepsilon_{\rm tr}$.  Define
\[
 P_J(\xi,\eta):=\sum_{j=0}^{J}\frac{(\ri\eta)^j}{j!}
                    \chi^{(j)}(\xi),
 \qquad
 \widetilde\chi_J(\xi+\ri\eta):=\chi_{\rm v}(\eta)P_J(\xi,\eta).
\]
The adjacent Taylor terms cancel:
\[
\begin{aligned}
 2\bar\partial P_J
 =\sum_{j=0}^{J}\frac{(\ri\eta)^j}{j!}\chi^{(j+1)}(\xi)
   -\sum_{j=0}^{J-1}\frac{(\ri\eta)^j}{j!}\chi^{(j+1)}(\xi)
 =\frac{(\ri\eta)^J}{J!}\chi^{(J+1)}(\xi).
\end{aligned}
\]
The second sum is zero when $J=0$.  Consequently,
\[
 2\bar\partial\widetilde\chi_J
 =\chi_{\rm v}(\eta)\frac{(\ri\eta)^J}{J!}\chi^{(J+1)}(\xi)
   +\ri\chi_{\rm v}'(\eta)P_J(\xi,\eta).
\]
Thus $|\bar\partial\widetilde\chi_J|\le C_J\eta^J$ for
$0\le\eta\le\varepsilon_{\rm tr}/2$, whereas
$|\bar\partial\widetilde\chi_J|\le C_J$ on the transition strip
$\varepsilon_{\rm tr}/2\le\eta\le\varepsilon_{\rm tr}$.
Only the derivatives of $\chi$ through order $J+1$ enter these bounds.

Choose $a<1<b$ so that $[a,b]\subset I_{\rm tr}$ and
$\supp\chi\Subset(a,b)$.  For $0<\epsilon<\varepsilon_{\rm tr}/2$,
consider the interior rectangle
\[
 \mathcal R_\epsilon
 :=\{\xi+\ri\eta:a<\xi<b,\quad
                   \epsilon<\eta<\varepsilon_{\rm tr}\}.
\]
Its closure is contained in the open upper chart.  The contour and the
two height regions are shown in Figure~\ref{fig:holomorphic}.
\begin{figure}[htpb]
\centering
\begin{tikzpicture}[x=1cm,y=1cm,>=Stealth]
\draw[dashed,gray] (-4.1,0) rectangle (4.1,4.4);
\fill[gray!9] (-3,.30) rectangle (3,1.10);
\fill[gray!20] (-3,1.10) rectangle (3,2.20);
\draw[thick] (-3,.30) rectangle (3,2.20);
\draw[dashed] (-3,1.10)--(3,1.10);
\draw[->] (-4.7,0)--(4.9,0) node[right] {$\xi=\Re z$};
\draw[->] (-4.7,-.12)--(-4.7,4.8) node[above] {$\eta=\Im z$};
\draw (-4.77,.30)--(-4.63,.30) node[left=5pt] {$\epsilon$};
\draw (-4.77,1.10)--(-4.63,1.10) node[left=5pt] {$\varepsilon_{\rm tr}/2$};
\draw (-4.77,2.20)--(-4.63,2.20) node[left=5pt] {$\varepsilon_{\rm tr}$};
\draw (-4.77,4.40)--(-4.63,4.40) node[left=5pt] {$2\varepsilon_{\rm tr}$};
\draw[->,thick] (-1.1,.30)--(1.1,.30);
\draw[->,thick] (3,.65)--(3,.95);
\draw[->,thick] (1.1,2.20)--(-1.1,2.20);
\draw[->,thick] (-3,.95)--(-3,.65);
\node at (0,3.75) {$\mathcal U_{\varepsilon_{\rm tr}}$};
\node at (0,2.65) {$\widetilde\chi_J=0$ near the top edge};
\node[rotate=90] at (-3.40,1.25) {$\widetilde\chi_J=0$};
\node[rotate=90] at (3.40,1.25) {$\widetilde\chi_J=0$};
\node at (0,.72) {$\bar\partial\widetilde\chi_J=O(\eta^J)$};
\node[align=center] at (0,1.63)
 {$|\bar\partial\widetilde\chi_J|\le C_J$\\
  $\re^{-ct\eta}\le \re^{-ct\varepsilon_{\rm tr}/2}$};
\draw[very thick] (-2.1,0)--(2.1,0);
\draw (-3,-.08)--(-3,.08) node[below=5pt] {$a$};
\draw (3,-.08)--(3,.08) node[below=5pt] {$b$};
\fill (0,0) circle (1.5pt);
\node[below=5pt] at (0,0) {$\rho=1$};
\draw (-2.1,-.48)--(-2.1,-.62)--(2.1,-.62)--(2.1,-.48);
\node[below] at (0,-.62) {$\supp\chi\Subset(a,b)$};
\end{tikzpicture}
\caption{The interior rectangle $\mathcal R_\epsilon$ in the upper chart.
Only the bottom edge contributes to Cauchy--Green.  The two shaded height
regions use Taylor cancellation and exponential damping, respectively.}
\label{fig:holomorphic}
\end{figure}
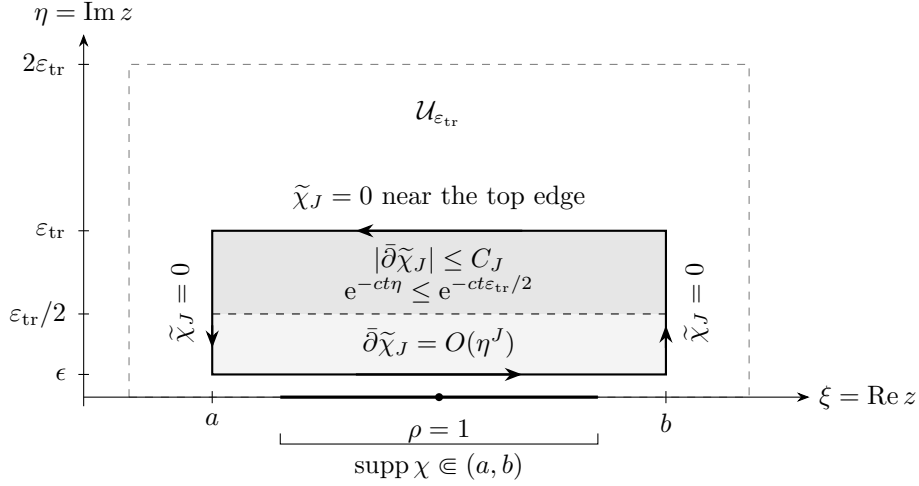
The extension vanishes near the two vertical sides, since $\chi$ and
its derivatives vanish near $a$ and $b$, and near the top edge by the
choice of $\chi_{\rm v}$.  Holomorphicity of the phase and amplitude gives
\[
 \bar\partial\!\left(\re^{\ri t\Psi_\lambda^\sharp}
                      A^\sharp\widetilde\chi_J\right)
 =\re^{\ri t\Psi_\lambda^\sharp}A^\sharp
                    \bar\partial\widetilde\chi_J.
\]
For the positively oriented boundary of $\mathcal R_\epsilon$, the
bottom edge runs from left to right.  The Gauss--Green formula, applied with
$\mathrm d(F\,\mathrm d z)=2\ri\bar\partial F\,\rd\xi\wedge\rd\eta$, gives
\[
\begin{aligned}
 \int_a^b
 \re^{\ri t\Psi_\lambda^\sharp(\xi+\ri\epsilon)}
 A^\sharp(\xi+\ri\epsilon)
 \widetilde\chi_J(\xi+\ri\epsilon)\,\rd\xi
 =2\ri\iint_{\mathcal R_\epsilon}
 \re^{\ri t\Psi_\lambda^\sharp(z)}A^\sharp(z)
           \bar\partial\widetilde\chi_J(z)\,\rd\xi\,\rd\eta.
\end{aligned}
\]

The boundary traces of $A^\sharp$ and $\Psi_\lambda^\sharp$ are continuous
on $[a,b]$, including across the single contact point $\rho=1$.
Their uniform bounds and the convergence
$\widetilde\chi_J(\xi+\ri\epsilon)\to\chi(\xi)$ justify passage to
the limit $\epsilon\downarrow0$ on the bottom edge by dominated
convergence.  By \eqref{eq:deep-timelike-phase-lower-appendix}, the
absolute value of the area integral is bounded uniformly in $\lambda$ by
\[
 C_J\int_0^{\varepsilon_{\rm tr}/2}\re^{-ct\eta}\eta^J\,\rd\eta
       +C_J\re^{-ct\varepsilon_{\rm tr}/2}
 \le C_Jt^{-J-1},\qquad t\ge1.
\]
This proves \eqref{eq:holomorphic-cancellation-v4}.
\end{proof}
}

We now prove the core estimate, using the radial representation
in the conic region and Lemma~\ref{lem:holomorphic-terminal-v4}
in the interior region.

\begin{lemma}
\label{lem:terminal-inversion-appendix}
For $t\ge2$ and $0\le r\le2t$,
\begin{equation}\label{eq:terminal-global-appendix}
 |G_{\rm core}(t,r)|\le Ct^{-3/2}.
\end{equation}
\end{lemma}

\begin{proof}
The bounded-time estimate \eqref{eq:threshold-core-bounded-time-v196}
reduces the proof to $t\ge T_0$.

\emph{The conic region $v_0t\le r\le2t$.}
The small-$L_-$ region of \eqref{eq:terminal-transverse-rho-v196}
and the transverse errors in \eqref{eq:transverse-layer-error-v18}
contribute $Ct^{-2}$.  Fix $h_0>0$ so that the continuum support of $\chi_\rho$ lies
in $(1,1+h_0)$.  For $v=r/t$, write the two remaining
integrals in \eqref{eq:terminal-conic-residues-v18} as
\[
 I_{\rm out}:=\int \re^{\ri t\phi_v(\rho)}B^{\rm out}_{t,v}(\rho)\,\rd\rho,
 \qquad
 I_{\rm in}:=\int \re^{\ri t\psi_v(\rho)}B^{\rm in}_{t,v}(\rho)\,\rd\rho,
\]
where $\phi_v$ is defined in
\eqref{eq:terminal-real-phase-appendix},
$\psi_v(\rho):=(1+v\rho)\omega(\rho)$, and the non-oscillatory
amplitudes are
\[
 \begin{aligned}
 B^{\rm out}_{t,v}(\rho)
 &:=\chi_\rho(\rho)(1-\eta_0(L_-))\mathcal R_4^\sigma(\rho)
       \re^{-t|1-v\rho|\nu(\rho)},\\
 B^{\rm in}_{t,v}(\rho)
 &:=\chi_\rho(\rho)\mathcal R_4(\rho)
       \re^{-t(1+v\rho)\nu(\rho)}.
 \end{aligned}
\]
Thus $G_{\rm core}=(2\pi^2r)^{-1}\Re(I_{\rm in}-I_{\rm out})
+O(t^{-2})$ in this region.

We first check the amplitude bounds.  By
Proposition~\ref{prop:terminal-jet-section2}, the derivatives of
$\mathcal R_4$ and its conjugate are bounded for $\gamma_1>1$
and at most logarithmic for $\gamma_1=1$; hence their variations
are finite.  The function $L_-$ is monotone in $\rho$, so
$1-\eta_0(L_-)$ has fixed variation, and the choice between
$\mathcal R_4$ and its conjugate adds at most one jump.
Lemma~\ref{lem:terminal-damping-bv-v2} bounds the variation of
both damping factors on the continuum side; on the pole side
both factors equal one.  The product inequality
\[
 \operatorname{Var}(fg)\le\|f\|_\infty\operatorname{Var}(g)
                      +\|g\|_\infty\operatorname{Var}(f)
\]
therefore gives, on each side of $\rho=1$,
\[
 \|B^{\rm out}_{t,v}\|_\infty+\operatorname{Var}(B^{\rm out}_{t,v})
 +\|B^{\rm in}_{t,v}\|_\infty+\operatorname{Var}(B^{\rm in}_{t,v})
 \le C,\qquad t\ge T_0,\quad v_0\le v\le2.
\]

For the phase factors in
\eqref{eq:terminal-phase-derivative-appendix}, direct differentiation
shows that
\[
 (g^{-1})'=-\frac{2\omega'+\rho\omega''}{g^2},\qquad
 V'=\frac{\omega\omega''-2(\omega')^2}{g^2}.
\]
When $\gamma_1>1$, the one-sided expansions
\eqref{eq:Omega0-terminal-general-section2}--
\eqref{eq:Theta-integer-section2} give bounded $\omega,\omega'$,
$g\to\omega_*+\omega_1>0$, and
\[
 |\omega''(1+h)|\lesssim
 1+\one_{\{1<\gamma_1<2\}}|h|^{\gamma_1-2}
   +\one_{\{\gamma_1=2\}}|\log|h||.
\]
This integrable bound proves that $g^{-1}$ and $V$ have bounded
variation on both sides.  For $\gamma_1=1$ on the pole side,
\eqref{eq:gamma-one-first-derivative-section2}--
\eqref{eq:gamma-one-second-derivative-section2} and
\eqref{eq:Vprime-gamma-one-section2} instead give, with $y=1-\rho>0$,
\[
 g(1-y)\asymp\log(1/y),\qquad
 |\partial_y g(1-y)^{-1}|+|\partial_y V(1-y)|
       \lesssim\frac1{y\log^2(1/y)}.
\]
Thus the same bounded-variation conclusion holds on the pole side
for every $\gamma_1\ge1$.

We estimate the outgoing integral first.  On the pole side,
$\nu=0$ and $\phi_v'=g(V-v)$.  The amplitude bound above,
\eqref{eq:V-sublevel-section2}, and
Lemma~\ref{lem:sublevel-vdc-appendix} give $Ct^{-1/2}$ before
$r^{-1}$, uniformly for $v\in[v_0,2]$ and every $\gamma_1\ge1$.
On the continuum side when $\gamma_1>2$, the same sublevel estimate
holds on the second component of $I_{\rm osc}$ defined before
\eqref{eq:V-sublevel-section2}, so the same argument applies.  These estimates allow ordinary stationary
points and give $Ct^{-3/2}$ after multiplication by $r^{-1}$.

For the continuum with $1<\gamma_1\le2$, recall the contact velocity
$v_*:=\lim_{\rho\to1^-}V(\rho)$ from \eqref{eq:vstar-section2}.
The cutoff was chosen so that, for a fixed
$0<\epsilon_v<(1-v_*)/4$,
\[
 |V-v_*|<\epsilon_v/2,\qquad
 1-v\rho\ge c>0\quad\hbox{if }|v-v_*|<\epsilon_v
\]
on its continuum support.  In this range of $v$,
$\nu(1+h)\asymp h^{\gamma_1}$ gives directly
\begin{equation}\label{eq:terminal-direct-damping-v4}
 \frac1r\left|\int_{\rho>1}\re^{\ri t\phi_v(\rho)}
                 B^{\rm out}_{t,v}(\rho)\,\rd\rho\right|
 \le\frac C r\int_0^{h_0}\re^{-cth^{\gamma_1}}\,\rd h
 \lesssim t^{-1-1/\gamma_1}\le t^{-3/2}.
\end{equation}
At $\gamma_1=3/2$, $\Re c_+=0$ but $\Im c_+>0$ in
\eqref{eq:Theta-noninteger-section2}; the estimate uses this
imaginary coefficient.  It also includes $\gamma_1=2$.
If $|v-v_*|\ge\epsilon_v$, then $|V-v|\ge\epsilon_v/2$ and
\[
 \frac1{\phi_v'}=\frac{g^{-1}}{V-v}
\]
has uniformly bounded supremum and variation.  One integration by
parts in the $BV$ sense therefore gives $Ct^{-1}$ before $r^{-1}$,
and hence $Ct^{-2}$, on this part of the continuum.

For $\gamma_1=1$ on the continuum side, put $h=\rho-1>0$ and
$h_v=(1-v)/v$.  Since $\nu(1+h)\asymp h$ and $v\ge v_0$, the
outgoing integral is bounded, before $r^{-1}$, by
\begin{equation}\label{eq:terminal-two-zero-damping-v18}
 C\int_0^{h_0}\re^{-cth|h_v-h|}\,\rd h\le Ct^{-1/2},
 \qquad h_v\in\R.
\end{equation}
Indeed, for $h_v\le0$, use $h|h_v-h|\ge h^2$.  For $h_v>0$, split
$(0,h_0)$ into its intersections with
\[
 (0,h_v/2),\qquad(h_v/2,3h_v/2),\qquad(3h_v/2,\infty).
\]
On the first two intervals, respectively,
$h|h_v-h|\ge h_vh/2$ and
$h|h_v-h|\ge h_v|h-h_v|/2$, so each integral is bounded by
$C\min\{h_v,(th_v)^{-1}\}\le Ct^{-1/2}$.
On the third interval $h\ge3h_v/2$ gives
$h|h-h_v|\ge h^2/3$, and the Gaussian integral gives the same
bound.  Empty intervals contribute zero. Thus the outgoing estimate remains uniform as $v\to1$, even as the two zeros of $h|h_v-h|$ merge.

For the incoming integral, $\psi_v'=g(V+v)$.
When $\gamma_1>1$ on both sides, and when $\gamma_1=1$ on the pole
side, $V+v\ge v_0$ and
$(\psi_v')^{-1}=g^{-1}/(V+v)$ has uniformly bounded supremum and
variation.  The amplitude estimate and one integration by parts
give $Ct^{-1}$ before $r^{-1}$.  On the remaining continuum side
when $\gamma_1=1$, direct damping gives
$Cr^{-1}\int_0^{h_0}\re^{-cth}\,\rd h\le Ct^{-2}$.
This proves the conic estimate.

\emph{The interior region $0\le r\le v_0t$.}
We return to the Fourier inversion underlying
\eqref{eq:terminal-transverse-radial-v18}, retaining the angular
integration.  The spherical measure contributes $k^2$ instead of
the radial factor $k$, so the transverse transform is $\mathcal T_5$
rather than $\mathcal T_4$.  The resulting formula has no $r^{-1}$
factor and applies also at $r=0$.  With the same change of variable
$\rho=(1-y)^{-1}$ and Hermitian reflection, it is
\begin{equation}\label{eq:terminal-transverse-sphere-v18}
 G_{\rm core}(t,x)=\frac2{(2\pi)^3}\Re
 \int_{\mathbb S^2}\int_0^\infty\chi_\rho(\rho)\rho^{-2}
 \mathcal T_5(t/\rho+x\cdot\theta,\rho\Om_0(\rho))
 \,\rd\rho\,\rd\theta.
\end{equation}
The $\chi_\rho$-weighted integrand is defined on the fixed contact
chart and extended by zero outside it; no value of $\Om_0$ outside
that chart is used.  For each $\theta\in\mathbb S^2$, set
$\lambda=x\cdot\theta/t$.
The fixed choices of $X_0,v_0,T_0$ give $|\lambda|\le v_0$ and
\[
 L=t(1/\rho+\lambda)\ge t/(2X_0)\ge1,\qquad t\ge T_0.
\]
Thus only the positive-$L$ branch of
\eqref{eq:transverse-residues-v18} occurs.
The leading term, including the Jacobian, has amplitude
$\ri\rho^{-2}(\rho\Om_0)^4/2=\ri\rho^2\Om_0^4/2$.

The spherical amplitude is the boundary trace of
\[
 A^\sharp(z):=\frac{z^2\Omega^\sharp(z)^4}{2},
 \qquad
 A^\sharp(\rho+\ri0)=\frac{\rho^2\Om_0(\rho)^4}{2}.
\]
It is bounded and holomorphic on the same upper chart, with a
continuous boundary trace, including at $\gamma_1=1$.
The fixed cutoff satisfies $\chi_\rho\in C_c^\infty(I_{\rm tr})$,
and the phase is precisely $\Psi_\lambda^\sharp$ in
\eqref{eq:deep-timelike-phase-lower-appendix}.  Hence the relevant
boundary integral can be written as
\begin{equation}\label{eq:deep-timelike-terminal-appendix}
 \begin{aligned}
 I_\lambda^{\rm hol}(t)
 :=\int \re^{\ri t(1+\lambda\rho)\Om_0(\rho)}
       \frac{\rho^2\Om_0(\rho)^4}{2}\chi_\rho(\rho)\,\rd\rho
 =\int \re^{\ri t\Psi_\lambda^\sharp(\rho+\ri0)}
       A^\sharp(\rho+\ri0)\chi_\rho(\rho)\,\rd\rho.
 \end{aligned}
\end{equation}
Lemma~\ref{lem:holomorphic-terminal-v4}, applied with $\chi=\chi_\rho$ and $J=1$, gives
\begin{equation}\label{eq:deep-timelike-paired-bound-appendix}
 |I_\lambda^{\rm hol}(t)|\le Ct^{-2},\qquad |\lambda|\le v_0.
\end{equation}
Retaining the factor $\ri$ and the spherical coefficient,
\eqref{eq:terminal-transverse-sphere-v18} becomes
\[
 G_{\rm core}(t,x)=\frac2{(2\pi)^3}\Re
       \left\{\ri\int_{\mathbb S^2}I_\lambda^{\rm hol}(t)\,\rd\theta\right\}
       +O_N(t^{-N}),
\]
where the error is uniform in $\theta$ by
Lemma~\ref{lem:transverse-cauchy-v18}.
The estimate is applied before the sphere integration, so it also
covers $r=0$.  Together with the conic estimate it proves the lemma.
\end{proof}

\begin{proof}[Proof of Proposition~\ref{prop:threshold-bound-outline}]
It remains to estimate the ordinary complement in
\eqref{eq:threshold-core-split-v18} and then recover the spatial
weight.  Combining the two remainder terms gives
\[
 \widehat G_{\rm rem}(t,k)
 =\int_\R \re^{\ri t\tau}
 [\vartheta(k)-\zeta(k)\chi_{\rm tr}(k,\tau)]\,
 \rd\mathfrak M_k(\tau).
\]
For $s>0$, put
$b_{\rm rem}(k,s):=\vartheta(k)-\zeta(k)\chi_y(1-s)$.
The mirror cutoff vanishes on the positive spectral axis, so
$\chi_{\rm tr}(k,\pm ks)=\chi_y(1-s)$ on the cut and at the poles.
For $k\ne k_*$, \eqref{eq:physical-measure-nonterminal-section3}
therefore gives
\[
 \begin{aligned}
 \widehat G_{\rm rem}(t,k)
 ={}\one_{\{k<k_*\}}
 b_{\rm rem}\bigl(k,\tau_{\rm p}(k)/k\bigr)
 \sum_{\sigma=\pm}\varrho_\sigma(k)\re^{\ri\sigma t\tau_{\rm p}(k)}
 +k\sum_{\sigma=\pm}\int_\kappa^1
 \re^{\ri\sigma tks}b_{\rm rem}(k,s)J_k(\sigma ks)\,\rd s.
 \end{aligned}
\]
The single value $k=k_*$ does not affect spatial Fourier inversion.
Since $\vartheta=\zeta=1$ near $k_*$, $\chi_y=1$ near zero,
and $\tau_{\rm p}(k)/k\to1$ as $k\uparrow k_*$, the pole
coefficient vanishes on a fixed left neighborhood of $k_*$.
The remaining pole amplitudes are smooth and compactly supported
in $(0,k_*)$.  On this support, the positive quantities
$\tau_{\rm p}''(k)$ and $\tau_{\rm p}'(k)/k$ from
Proposition~\ref{prop:terminal-jet-section2} have uniform lower
bounds.  Lemma~\ref{lem:radial-dispersive-section4}, and its
complex conjugate for the negative phase, give $Ct^{-3/2}$
for their inverse transforms, allowing ordinary stationary points.

For the continuum, $b_{\rm rem}$ vanishes in a neighborhood of
$(k_*,1)$.  Recall from \eqref{eq:J-explicit-proof-section3} that
\[
 J_k(ks)=\frac{2\ri k^2\Bcal(s)}{d(k,s)}.
\]
On the open positive cut $\kappa<s<1$, $\Bcal(s)>0$.
At $s=\kappa$, $\Acal(\kappa)<0$; at $s=1$, the denominator is
$(k^2-k_*^2)^2$, whose only zero has been removed by the cutoff.
Compactness of the remaining positive $k$-annulus therefore gives
$d(k,s)\ge c>0$ on the support of the localized density and its
derivatives.  Near a positive endpoint $a\in\{\kappa,1\}$, let
$s_a(u)=a+\epsilon_a u$, where $u>0$ is the inward coordinate.
The factorization $\Bcal(s_a(u))=u^{\gamma_a}\beta_a(u)$ with
smooth $\beta_a$, and the conormal bounds of
Lemma~\ref{lem:A-endpoint-section2}, imply
\[
 \left|\partial_u^j\left[
 b_{\rm rem}(k,s_a(u))J_k(\sigma k s_a(u))\right]\right|
 \le C_j u^{\gamma_a-j},\qquad \sigma=\pm,
\]
through the fixed orders required below.  Indeed, the separated
denominator and its reciprocal have bounded conormal derivatives,
and the localizing cutoffs preserve these bounds.
Lemma~\ref{lem:finite-endpoint-v4}, with Fourier parameter $tk$,
then gives $C(tk)^{-\gamma_a-1}\le Ct^{-2}$, since $k$ stays
bounded above and below.  Under \eqref{H0}, combine the two signs
across $s=0$ in a single chart, smooth by
Lemma~\ref{lem:A-endpoint-section2}.  On the remaining cut
interior, repeated integration by parts gives $O_N(t^{-N})$.
Integrating these
absolute bounds over the compact spatial-frequency annulus
preserves them uniformly in $x$.

Combining the ordinary poles and continuum with
Lemma~\ref{lem:terminal-inversion-appendix} yields
\begin{equation}\label{eq:threshold-cone-section4}
 |G_{\rm th}(t,x)|\lesssim t^{-3/2},\qquad
       t\ge2,\quad |x|\le2t.
\end{equation}
Propositions~\ref{prop:pole-dispersive-outline} and
\ref{prop:green-LD-outline}, \eqref{eq:threshold-cone-section4},
and Lemma~\ref{lem:finite-propagation-green-section4} then give
\[
 \|G(t)\|_{L^\infty_x}\lesssim t^{-3/2},\qquad
 \supp G(t)\subset\overline B(0,t).
\]
Consequently
\begin{equation}\label{eq:green-L1-from-cone-section4}
 \|G(t)\|_{L^1_x}\lesssim t^{3/2},\qquad t\ge2.
\end{equation}
Since $G_{\rm th}=\vartheta(D_x)G=\check\vartheta*_xG$ and
$\check\vartheta$ is Schwartz, for $r\ge2t$ we obtain
\[
 |G_{\rm th}(t,x)|\le C_Mr^{-M}\|G(t)\|_{L^1_x}
       \le C_Mt^{3/2}r^{-M}.
\]
Taking $M=N+3$ proves
$|G_{\rm th}(t,x)|\le C_Nt^{-3/2}(t/r)^N$ in this exterior
region.  Together with \eqref{eq:threshold-cone-section4} this
establishes \eqref{eq:threshold-global-outline} for $t\ge2$.
For $0<t\le2$, the uniform $L^1$ bound from
Lemma~\ref{lem:finite-propagation-green-section4} and the same
convolution give $|G_{\rm th}(t,x)|\le C_N\ang x^{-N}$.
\end{proof}

\subsection{Proof of Theorem~\ref{thm:main-green}}

\begin{proof}[Proof of Theorem~\ref{thm:main-green}]
Under the shell hypotheses, the exact decomposition is \eqref{eq:wave-components-section3}.  By
Propositions~\ref{prop:pole-dispersive-outline} and
\ref{prop:threshold-bound-outline},
\[
|G_{\rm wave}(t,x)|
\le C_N\ang t^{-3/2}
\frac{\ang t^N}{\ang{t,x}^N},
\qquad N\ge0.
\]
For $t+|x|\ge1$, Proposition~\ref{prop:green-LD-outline} gives
\eqref{eq:main-green-LD-envelope}, since its additional short-time
factor is at least one.  For $0<t+|x|\le1$, the spatial factor
on the right-hand side of \eqref{eq:main-green-LD-envelope} is
bounded above and below by positive constants.  The estimate
\eqref{eq:short-time-LD-section4} therefore gives the same bound,
with the singular term confined to $\kappa t\le|x|\le t$.
Finite propagation follows from
Lemma~\ref{lem:finite-propagation-green-section4}.
Under \eqref{H0}, Lemma~\ref{lem:A-endpoint-section2} allows the
same argument with $\kappa=0$, only the outer front, and the
short-time indicator $\one_{\{|x|\le t\}}$.
\end{proof}

\section{The complete electric field}
\label{sec:electric}

We first collect estimates for the pole primitives and the
free-transport source.  Two integrations in time then give the exact
electric decomposition.  We estimate its wave, corrected
low-frequency, and fixed-gap components separately and combine them
to prove Theorem~\ref{thm:main-electric}.
For a filled-ball background, Lemma~\ref{lem:A-endpoint-section2}
provides the center regularity, and only the outer front remains.
We use the Poisson operator
\[
\mathcal P_x:=-\nabla_x(-\Delta_x)^{-1},
\qquad
\widehat{\mathcal P_x f}(\xi)=-\frac{\ri\xi}{|\xi|^2}\widehat f(\xi).
\]
The subscript specifies the spatial variable on which the operator
acts; after rescaling to $z$, we write $\mathcal P_z$.

\subsection{Pole primitives and free-transport estimates}

For $j=1,2,3$, define the temporal pole primitives
\begin{equation}\label{eq:pole-primitives-section4}
 \widehat K_j(t,\xi)=\chi_\ell(k)\sum_{\sigma=\pm}
 \frac{\varrho_\sigma(k)}{(\ri\sigma\tau_{\rm p}(k))^j}
 \re^{\ri\sigma t\tau_{\rm p}(k)},\qquad k=|\xi|.
\end{equation}
For $t>0$, ordinary time differentiation gives
$\partial_tK_1=G_{\rm KG}$ and $\partial_tK_j=K_{j-1}$ for $j=2,3$.
On the fixed support of $\chi_\ell$, $\tau_{\rm p}$ has a positive lower
bound.  Together with the smooth dependence of $\tau_{\rm p}$ and
$r_{\rm p}$ on $k^2$ near zero in
Proposition~\ref{prop:pole-branch-section2}, this shows that the amplitudes
$\chi_\ell\varrho_\sigma(\ri\sigma\tau_{\rm p})^{-j}$ are compactly
supported and Cartesian smooth, including at the origin.
The phase Hessians were checked in
Proposition~\ref{prop:pole-dispersive-outline}; the Cartesian
conclusion of Lemma~\ref{lem:radial-dispersive-section4} therefore
gives, for the fixed indices $j=1,2,3$ and each fixed $N\ge0$,
\begin{equation}\label{eq:pole-primitives-bound-section4}
 |K_j(t,x)|\le C_N\ang t^{-3/2}
                  \frac{\ang t^N}{\ang{t,x}^N}.
\end{equation}
The electric argument uses this estimate, the scalar threshold
convolution, and the low-frequency symbols
\eqref{eq:R-minus-P-factor-section3}--\eqref{eq:Rcal-jump-conormal-section3}
of Proposition~\ref{prop:low-frequency-factor-section3}\textup{(iii)}.

Define the free-transport tensor moments by
\begin{equation}
\label{eq:free-moments-section5}
M_\ell(t,x)
:=\int_{\R^3}v^{\otimes\ell}h_{\rm in}(x-tv,v)\,\rd v,
\qquad \ell=0,1,2,3,
\end{equation}
where $M_0=\rho^{\rm free}$.  The operator $\operatorname{div}_x$
contracts the derivative with the last tensor index.  Direct differentiation gives
\begin{equation}
\label{eq:moment-hierarchy-section5}
\begin{gathered}
\partial_tM_\ell+\operatorname{div}_x M_{\ell+1}=0,
\qquad \ell=0,1,2,\\
\partial_t^j\rho^{\rm free}=(-1)^j\nabla_x^j\!:\!M_j,
\qquad j=1,2,3.
\end{gathered}
\end{equation}

The following proposition estimates the free-transport moments,
their Poisson transforms, and the terms localized away from zero
spatial frequency.
\begin{proposition}
\label{lem:free-moments-section5}
Let $M\ge11$.  The following estimates hold for
$t\ge0$.
\begin{enumerate}[label=\textup{(\roman*)}]
\item \emph{Free-transport moments.}
For $0\le\ell\le3$, $\ell+j\le3$, $j+|\alpha|\le3$, and
$0\le L\le8$,
\begin{align}
 |\partial_t^j\partial_x^\alpha M_\ell(t,x)|
 &\le C\ang t^{L-3-j-|\alpha|}\ang{t,x}^{-L}
                  [h_{\rm in}]_{M},
                         \label{eq:free-moment-pointwise-section5}\\
 \|\partial_t^j\partial_x^\alpha M_\ell(t)\|_1
 &\le C\ang t^{-j-|\alpha|}[h_{\rm in}]_{M}.
                         \label{eq:free-moment-L1-section5}
\end{align}
\item \emph{Poisson transforms.}
For $j=2,3$, $0\le\ell\le3$, and $|\beta|=j$, each tensor component
$(M_\ell)_I$, $I\in\{1,2,3\}^\ell$, satisfies
\begin{align}
 |\mathcal P_x\partial_x^\beta(M_\ell)_I(t,x)|
 &\le C\ang{t,x}^{-j-2}[h_{\rm in}]_{M},
                         \label{eq:poisson-moment-pointwise-section5}\\
 \|\mathcal P_x\partial_x^\beta(M_\ell)_I(t)\|_1
 &\le C\ang t^{1-j}[h_{\rm in}]_{M}.
                         \label{eq:poisson-moment-L1-section5}
\end{align}
For $\ell=0$, the empty index denotes the scalar $M_0$.
\item \emph{Estimates away from zero frequency.}
Let $a\in C^\infty(\R^3)$ be fixed and suppose that, for some $c_0>0$,
\[
 a(\xi)=0\quad (|\xi|\le c_0),\qquad
 |\partial_\xi^\beta a(\xi)|\le C_{a,\beta}\ang\xi^{-|\beta|},
 \quad |\beta|\le6.
\]
Then, for $0\le\ell\le3$,
\begin{equation}\label{eq:band-source-section5}
 |a(D_x)M_\ell(t,x)|+|a(D_x)E^{\rm free}(t,x)|
 \le C_a\ang t^{-4}\ang x^{-6}[h_{\rm in}]_{M}.
\end{equation}
Here $C_a$ depends on $c_0$ and the displayed multiplier bounds.
\end{enumerate}
\end{proposition}
\begin{proof}
\emph{(i) Free-transport moments.}
Since $M+6\ge17$, \eqref{eq:data-norm} gives
\[
 |\partial_x^\mu\partial_v^\nu h_{\rm in}(y,v)|
 \le [h_{\rm in}]_{M}\ang{y,v}^{-M-6},
 \qquad |\mu|+|\nu|\le M.
\]
For $t\ge1$, iterating \eqref{eq:moment-hierarchy-section5}
gives $\partial_t^jM_\ell=(-1)^j\nabla_x^j\!:\!M_{\ell+j}$, and the
change of variables $y=x-tv$ yields
\[
 M_{\ell+j}(t,x)=t^{-3}\int_{\R^3}
  \left(\frac{x-y}{t}\right)^{\!\otimes(\ell+j)}
  h_{\rm in}\left(y,\frac{x-y}{t}\right)\,\rd y.
\]
Each spatial derivative contributes a factor $t^{-1}$ and either
differentiates the data in $v$ or lowers the degree of the velocity
polynomial.  Writing $z=x/t$ and $v=z-y/t$, we therefore obtain
\[
 |\partial_t^j\partial_x^\alpha M_\ell(t,x)|
 \le Ct^{-3-j-|\alpha|}[h_{\rm in}]_{M}
       \int_{\R^3}\ang v^{\ell+j}\ang{y,v}^{-M-6}\,\rd y.
\]
The data derivatives have order at most $j+|\alpha|\le3$.
Since $\ang z\lesssim\ang{y,v}$ and $\ang{y,v}\ge\ang y$,
\[
 \ang z^L\int_{\R^3}\ang v^{\ell+j}\ang{y,v}^{-M-6}\,\rd y
 \lesssim\int_{\R^3}\ang y^{L+\ell+j-M-6}\,\rd y\le C,
 \qquad M+6-L-\ell-j\ge6>3.
\]
This proves \eqref{eq:free-moment-pointwise-section5} for $t\ge1$,
using $\ang{x/t}\asymp\ang{t,x}/\ang t$.

For $0\le t\le1$, differentiating the original moment integral
introduces $j$ spatial derivatives and a velocity polynomial of
degree at most $\ell+j$.  The inequality
$\ang x\lesssim\ang{x-tv,v}$ gives
\[
 \begin{aligned}
 \ang x^L|\partial_t^j\partial_x^\alpha M_\ell(t,x)|
 &\le C[h_{\rm in}]_{M}
       \int_{\R^3}\ang v^{\ell+j}\ang{x-tv,v}^{L-M-6}\,\rd v\\
 &\le C[h_{\rm in}]_{M}
       \int_{\R^3}\ang v^{L+\ell+j-M-6}\,\rd v
 \le C[h_{\rm in}]_{M}.
 \end{aligned}
\]
As $\ang t\asymp1$, this completes
\eqref{eq:free-moment-pointwise-section5}.  Taking $L=4$ and using
$\int_{\R^3}\ang{t,x}^{-4}\,\rd x=C\ang t^{-1}$ proves
\eqref{eq:free-moment-L1-section5}.

\emph{(ii) Poisson transforms.}
Fix $j=2,3$, $|\beta|=j$, and a component $(M_\ell)_I$.
For $t\ge1$, let $H_{\ell,t}(z)=t^3M_\ell(t,tz)$.
Applying \eqref{eq:free-moment-pointwise-section5} with time-derivative
index zero and $L=j+3$ gives
\[
 \sum_{|\alpha|\le j}
 \|\ang z^{j+3}\partial_z^\alpha(H_{\ell,t})_I\|_\infty
 \le C[h_{\rm in}]_{M}.
\]
Indeed, $\partial_z^\alpha H_{\ell,t}
=t^{3+|\alpha|}(\partial_x^\alpha M_\ell)(t,tz)$.
Since $\mathcal P_x\partial_x^\beta$ is homogeneous of order $j-1$,
Lemma~\ref{lem:coulomb-moments-v4} implies
\[
 \begin{aligned}
 |(\mathcal P_x\partial_x^\beta(M_\ell)_I)(t,tz)|
 =t^{-j-2}|(\mathcal P_z\partial_z^\beta(H_{\ell,t})_I)(z)|
 \le Ct^{-j-2}\ang z^{-j-2}[h_{\rm in}]_{M}.
 \end{aligned}
\]
For $0\le t\le1$, the same zero time-derivative index and
$L=j+3$ in (i) bound
the unscaled weighted $C^j$ norm of $(M_\ell)_I(t)$; applying
Lemma~\ref{lem:coulomb-moments-v4} directly proves
\eqref{eq:poisson-moment-pointwise-section5} in this range as well.
Since $j+2>3$, integrating this bound in $x$ gives
\eqref{eq:poisson-moment-L1-section5}.

\emph{(iii) Estimates away from zero frequency.}
Fourier inversion gives
\[
 \begin{aligned}
 \ang x^6|a(D_x)M_\ell(t,x)|
 \le C\|(1-\Delta_\xi)^3(a\widehat M_\ell)(t)\|_{L^1_\xi}
 \le C\sum_{|\beta|\le6}
       \|\partial_\xi^\beta(a\widehat M_\ell)(t)\|_{L^1_\xi}.
 \end{aligned}
\]
We estimate the last sum by
$C_a\ang t^{-4}[h_{\rm in}]_{M}$.
From \eqref{eq:free-moments-section5}, a typical term in
$\partial_\xi^\beta(a\widehat M_\ell)$ is, up to a fixed coefficient,
\[
 t^j(\partial_\xi^\mu a)(\xi)
 \iint_{\R^6}\re^{-\ri\xi\cdot(y+tv)}
 y^\sigma P_{\ell,j}(v)h_{\rm in}(y,v)\,\rd y\,\rd v,
 \qquad |\mu|+j+|\sigma|=|\beta|\le6,
\]
where $P_{\ell,j}$ has degree at most $\ell+j$.
Since $|\sigma|+\ell+j\le9$ and $M+6-9\ge8>6$, the data
assumption bounds the $L^1_{y,v}$ norms of all derivatives of
$y^\sigma P_{\ell,j}(v)h_{\rm in}$ of order at most $10\le M$
by $C[h_{\rm in}]_{M}$.

For $0\le t\le1$, put $k=|\xi|\ge c_0$ and $\theta=\xi/k$.
Integrating by parts four times in the direction
$\theta\cdot\nabla_y$ bounds the typical term by
\[
 C_a k^{-4}\ang\xi^{-|\mu|}[h_{\rm in}]_{M}.
\]
For $t\ge1$, integration by parts $4+j\le10$ times in the
velocity direction $\theta\cdot\nabla_v$ gives
\[
 C_a t^j(tk)^{-4-j}\ang\xi^{-|\mu|}[h_{\rm in}]_{M}
 =C_a t^{-4}k^{-4-j}\ang\xi^{-|\mu|}[h_{\rm in}]_{M}.
\]
Both frequency bounds are integrable on $|\xi|\ge c_0$ in three
dimensions.  Consequently,
\[
 \sum_{|\beta|\le6}
 \|\partial_\xi^\beta(a\widehat M_\ell)(t)\|_{L^1_\xi}
 \le C_a\ang t^{-4}[h_{\rm in}]_{M},
\]
which proves the $M_\ell$ bound in \eqref{eq:band-source-section5}.
Finally,
$\widehat E^{\rm free}=-\ri\xi|\xi|^{-2}\widehat M_0$.
The multiplier $-\ri a(\xi)\xi/|\xi|^2$ vanishes near zero and
satisfies the same six derivative bounds.  Applying the preceding
argument with $\ell=0$ proves the $E^{\rm free}$ bound.
\end{proof}

\subsection{The pole normal form and the exact field decomposition}

Recall the pole primitives \(K_j\) from
\eqref{eq:pole-primitives-section4}.  Two integrations in time
separate their initial contributions from the source history and
produce a boundary term that cancels the leading free density.

\begin{proposition}
\label{prop:pole-normal-form-section5}
For the free-transport density \(\rho^{\rm free}\),
\begin{equation}
\label{eq:pole-normal-form-section5}
\begin{aligned}
G_{\rm KG}*_{t,x}\rho^{\rm free}
={}K_1(t)*_x\rho^{\rm free}(0)
 +K_2(t)*_x\partial_t\rho^{\rm free}(0)
+K_2*_{t,x}\partial_t^2\rho^{\rm free}
 -K_1(0)*_x\rho^{\rm free}(t).
\end{aligned}
\end{equation}
Moreover,
\begin{equation}
\label{eq:K2-zero-section5}
K_2(0)\equiv0,
\end{equation}
and there is a radial symbol $\xi\mapsto b_0(|\xi|)$ in $C_c^\infty(\mathbb R^3)$ such that
\begin{equation}
\label{eq:pole-sum-rule-section5}
\chi_\ell(k)-\widehat K_1(0,k)=k^2b_0(k),
\qquad b_0(0)=4\pi m_4.
\end{equation}
\end{proposition}

\begin{proof}
Since \(\partial_tK_1=G_{\rm KG}\) and
\(\partial_tK_2=K_1\), two integrations by parts on \([0,t]\) give
\begin{align*}
G_{\rm KG}*_{t,x}\rho^{\rm free}
={}&K_1(t)*_x\rho^{\rm free}(0)
 -K_1(0)*_x\rho^{\rm free}(t)
 +K_1*_{t,x}\partial_t\rho^{\rm free}\\
={}&K_1(t)*_x\rho^{\rm free}(0)
 +K_2(t)*_x\partial_t\rho^{\rm free}(0)
 +K_2*_{t,x}\partial_t^2\rho^{\rm free}\\
&-K_1(0)*_x\rho^{\rm free}(t)
 -K_2(0)*_x\partial_t\rho^{\rm free}(t).
\end{align*}
The residue symmetry
\(\varrho_\sigma=\sigma\ri r_{\rm p}\) yields
\[
\widehat K_2(0,k)
=-\frac{\chi_\ell(k)}{\tau_{\rm p}(k)^2}
\bigl(\varrho_+(k)+\varrho_-(k)\bigr)=0,
\]
which proves \eqref{eq:pole-normal-form-section5} and
\eqref{eq:K2-zero-section5}.  Similarly,
\[
\widehat K_1(0,k)=\chi_\ell(k)\frac{2r_{\rm p}(k)}{\tau_{\rm p}(k)},
\qquad
\frac{2r_{\rm p}(k)}{\tau_{\rm p}(k)}
=1-4\pi m_4k^2+O(k^4)
\]
by Proposition~\ref{prop:pole-branch-section2}.
That proposition also gives smooth functions
\(\tau_{\rm p}(k)=\mathsf T(k^2)\) and
\(r_{\rm p}(k)=\mathsf R(k^2)\) near zero, with
\(\mathsf T(0)=1\), \(\mathsf R(0)=1/2\).  Thus
\[
B(s):=1-\frac{2\mathsf R(s)}{\mathsf T(s)},
\qquad B(0)=0,
\qquad \frac{B(s)}s=\int_0^1 B'(\theta s)\,\rd\theta.
\]
The last identity extends the quotient smoothly through zero.
Since \(\chi_\ell=1\) there, \(b_0(|\xi|)\) is Cartesian smooth at the
origin and \(b_0(0)=B'(0)=4\pi m_4\).  Away from zero its
smoothness follows from the pole branch, and its compact support
from that of \(\chi_\ell\).
\end{proof}

We now derive the decomposition announced in the Introduction.
For a multiplier $m(z,|\xi|)$, let $m^\vee$ denote its causal
inverse Laplace transform followed by inverse spatial Fourier
transformation.  Applying \(\chi_\ell(D_x)\) to
\eqref{eq:green-field-representation-intro} and separating the
pole response gives
\[
\chi_\ell(D_x)\rho
=\chi_\ell(D_x)\rho^{\rm free}
 +G_{\rm KG}*_{t,x}\rho^{\rm free}
 +[\chi_\ell(R-P)]^\vee*_{t,x}\rho^{\rm free}.
\]
Since \(\chi_\ell+\vartheta+\chi_h=1\), adding back the other frequency
components yields
\begin{equation}
\label{eq:density-pre-normal-form-section5}
\begin{aligned}
\rho
={}\chi_\ell(D_x)\rho^{\rm free}
 +G_{\rm KG}*_{t,x}\rho^{\rm free}
 +\vartheta(D_x)\rho
+[\chi_\ell(R-P)]^\vee*_{t,x}\rho^{\rm free}+\chi_h(D_x)\rho.
\end{aligned}
\end{equation}
The boundary term in
\eqref{eq:pole-normal-form-section5} cancels the leading free-density
contribution by \eqref{eq:pole-sum-rule-section5}.  Independently,
\eqref{eq:R-minus-P-slit-factor-section3} gives a factor \(k^2\)
in \(R-P\) for \(0<k\le k_0\).  Both factors compensate for the
Coulomb singularity, while the continuum term still requires
estimates of its source history.  We therefore define
\begin{align}
\rho_{\rm wave}
:={}&K_1(t)*_x\rho^{\rm free}(0)
 +K_2(t)*_x\partial_t\rho^{\rm free}(0)
 +K_2*_{t,x}\partial_t^2\rho^{\rm free}
 +\vartheta(D_x)\rho,
\label{eq:rho-wave-section5}\\
\rho_{\rm LD}
:={}&\chi_\ell(D_x)\rho^{\rm free}-K_1(0)*_x\rho^{\rm free}(t)
 +[\chi_\ell(R-P)]^\vee*_{t,x}\rho^{\rm free}+\chi_h(D_x)\rho,
\label{eq:rho-LD-section5}\\
E_{\rm wave}:={}&\mathcal P_x\rho_{\rm wave},
\qquad E_{\rm LD}:=\mathcal P_x\rho_{\rm LD}.
\label{eq:electric-decomposition-section5}
\end{align}
Then
\begin{equation}
\label{eq:electric-exact-sum-section5}
\rho=\rho_{\rm wave}+\rho_{\rm LD},
\qquad E=E_{\rm wave}+E_{\rm LD}.
\end{equation}
Equivalently, using \(\chi_\ell+\chi_h=1-\vartheta\) and
\(R_{\rm LD}=\chi_\ell(R-P)+\chi_h R\),
\begin{equation}
\label{eq:rho-LD-short-section5}
\rho_{\rm LD}
=\bigl[1-\vartheta(D_x)\bigr]\rho^{\rm free}
 -K_1(0)*_x\rho^{\rm free}(t)
 +G_{\rm LD}*_{t,x}\rho^{\rm free}.
\end{equation}
This is a pure spatial-frequency decomposition: the whole threshold
block remains in \(\vartheta(D_x)\rho\).  The symbols
\(\widehat K_1(0)\), \(\chi_\ell\), \(\chi_h\), and \(\vartheta\) are real
and even, while the pole pairs have Hermitian symmetry.  Hence both
densities and both electric fields are real for real data.

At the Laplace-multiplier level the corrected local term and the
pole-subtracted response in \eqref{eq:rho-LD-section5} form the block
\begin{equation}
\label{eq:corrected-low-multiplier-section5}
\bigl[\chi_\ell(Q-P)-\widehat K_1(0)\bigr]
\widetilde{\rho^{\rm free}}.
\end{equation}
For \(0<k\le k_0\) and \(\Re z>0\), equations
\eqref{eq:R-minus-P-slit-factor-section3} and
\eqref{eq:pole-sum-rule-section5} give
\begin{equation}
\label{eq:corrected-low-k2-section5}
\chi_\ell(k)\bigl[Q(z,k)-P(z,k)\bigr]-\widehat K_1(0,k)
=k^2\bigl[b_0(k)+\chi_\ell(k)\widetilde{\mathcal R}(z/k,k)\bigr].
\end{equation}

\subsection{The wave component}

We establish the pointwise estimate for the wave component.

\begin{proposition}
\label{prop:electric-wave-section5}
If \(M\ge11\), then
\begin{equation}
\label{eq:electric-wave-section5}
|E_{\rm wave}(t,x)|
\le C\frac{\ang t^{3/2}}{\ang{t,x}^3}[h_{\rm in}]_{M}.
\end{equation}
\end{proposition}

\begin{proof}
By \eqref{eq:rho-wave-section5}--\eqref{eq:electric-decomposition-section5},
the wave field consists of two initial pole terms, a forced pole term,
and the threshold contribution.  We use
Lemma~\ref{lem:radial-dispersive-section4} for the initial oscillatory
terms, including the one produced below by a third integration by parts
in time, and weighted source estimates for the remaining contributions.

For the first initial pole term, neutrality
\eqref{eq:neutrality} gives
\begin{equation}
\label{eq:neutral-quotient-section5}
\frac{\widehat{\rho^{\rm free}}(0,k\theta)}k
=\int_0^1\theta\cdot\nabla_\xi
\widehat{\rho^{\rm free}}(0,\lambda k\theta)\,\rd\lambda.
\end{equation}
Thus the actual amplitude of each oscillatory mode in
\(\mathcal P_x[K_1(t)*_x\rho^{\rm free}(0)]\) is
\[
a_{1,\sigma}(k,\theta)
=-\ri\theta\frac{\chi_\ell(k)\varrho_\sigma(k)}{\ri\sigma\tau_{\rm p}(k)}
\int_0^1\theta\cdot\nabla_\xi
\widehat{\rho^{\rm free}}(0,\lambda k\theta)\,\rd\lambda.
\]
We check the polar derivative condition of
Lemma~\ref{lem:radial-dispersive-section4}.  For example,
\[
D_k\int_0^1\theta\cdot\nabla_\xi
\widehat{\rho^{\rm free}}(0,\lambda k\theta)\,\rd\lambda
=\int_0^1\lambda k\,\theta\cdot
D_\xi^2\widehat{\rho^{\rm free}}(0,\lambda k\theta)\theta
\,\rd\lambda.
\]
Each angular derivative likewise either differentiates a smooth
spherical coefficient or introduces \(\lambda k\nabla_\xi\).
The chain rule therefore introduces no negative power of \(k\).
The radial prefactor has bounded \(D_k\) derivatives on its fixed
support by Proposition~\ref{prop:pole-branch-section2}.  Consequently
\[
\begin{aligned}
\|a_{1,\sigma}\|_{\mathrm{pol},6}
\le C\sum_{1\le|\alpha|\le7}
 \|\partial_\xi^\alpha\widehat{\rho^{\rm free}}(0)\|_\infty
\le C\iint_{\R^6}\ang x^7|h_{\rm in}(x,v)|\,\rd x\,\rd v
 \le C[h_{\rm in}]_{M}.
\end{aligned}
\]
Indeed \(M+6-7\ge10>6\).  Formula
\eqref{eq:polar-pole-source-section5} now proves the desired bound
for this term.

By \eqref{eq:moment-hierarchy-section5}, the second initial term
\(\mathcal P_x[K_2(t)*_x\partial_t\rho^{\rm free}(0)]\) has amplitude
\[
a_{2,\sigma}(k,\theta)
=-\theta\frac{\chi_\ell(k)\varrho_\sigma(k)}
                  {(\ri\sigma\tau_{\rm p}(k))^2}
 \bigl(\theta\cdot\widehat M_1(0,k\theta)\bigr).
\]
The same chain-rule calculation bounds its polar norm by
\(C[h_{\rm in}]_{M}\): six spatial moments and one velocity
moment suffice.  Applying \eqref{eq:polar-pole-source-section5}
again controls this contribution.

For the forced pole term, a third integration in time, as in
\cite[Section~4.2]{Ngu26}, gives
\begin{equation}
\label{eq:forced-pole-third-primitive-section5}
\begin{aligned}
K_2*_{t,x}\partial_t^2\rho^{\rm free}
={}K_3(t)*_x\partial_t^2\rho^{\rm free}(0)
 -K_3(0)*_x\partial_t^2\rho^{\rm free}(t)
+K_3*_{t,x}\partial_t^3\rho^{\rm free}.
\end{aligned}
\end{equation}
The moment estimates justify differentiation and the endpoint values
on this finite time interval.  After applying \(\mathcal P_x\), the
first term on the right has amplitude
\[
a_{3,\sigma}(k,\theta)
=\ri k\theta\frac{\chi_\ell(k)\varrho_\sigma(k)}
                     {(\ri\sigma\tau_{\rm p}(k))^3}
 \bigl(\theta\otimes\theta:\widehat M_2(0,k\theta)\bigr).
\]
Its polar derivatives are controlled as above, now using six spatial
and two velocity moments.  The joint moment degrees of the three
initial amplitudes are thus at most \(7,7,8\), all strictly
integrable under \eqref{eq:data-norm}.
The estimate \eqref{eq:polar-pole-source-section5} applies to this
third initial amplitude.  For the current-time boundary term,
\(K_3(0)\) is a Schwartz kernel.  The \(j=2\) case of
\eqref{eq:poisson-moment-pointwise-section5}, followed by
\eqref{eq:weighted-kernel-convolution-v4}, gives
\[
|K_3(0)*_x\mathcal P_x\partial_t^2\rho^{\rm free}(t,x)|
\le C\ang{t,x}^{-4}[h_{\rm in}]_{M}.
\]

For the remaining history term, keep
\(\mathcal P_x\partial_t^3\rho^{\rm free}\) intact.  The \(j=3\)
cases of \eqref{eq:poisson-moment-pointwise-section5} and
\eqref{eq:poisson-moment-L1-section5} give
\[
\begin{aligned}
|\mathcal P_x\partial_t^3\rho^{\rm free}(s,y)|
&\le C\ang{s,y}^{-5}[h_{\rm in}]_{M},\\
\|\mathcal P_x\partial_t^3\rho^{\rm free}(s)\|_1
&\le C\ang s^{-2}[h_{\rm in}]_{M}.
\end{aligned}
\]
Using \eqref{eq:pole-primitives-bound-section4} with weight four,
weakening the source's fifth spatial weight to the fourth, and
applying \eqref{eq:two-scale-convolution-section5}, we obtain,
since \(\ang{t-s}+\ang s\asymp\ang t\),
\[
\begin{aligned}
|K_3*_{t,x}(\mathcal P_x\partial_t^3\rho^{\rm free})(t,x)|
\le C\frac{\ang t^4}{\ang{t,x}^4}[h_{\rm in}]_{M}
 \int_0^t\ang{t-s}^{-3/2}\ang s^{-5}
      \min\{\ang{t-s}^3,\ang s^3\}\,\rd s.
\end{aligned}
\]
The time integral satisfies
\begin{equation}\label{eq:wave-time-convolution-section5}
\int_0^t\ang{t-s}^{-3/2}\ang s^{-5}
 \min\{\ang{t-s}^3,\ang s^3\}\,\rd s
\lesssim\ang t^{-3/2}.
\end{equation}
For \(t\ge1\), splitting at \(s=t/2\) bounds the two parts by
\[
C\ang t^{-3/2}\int_0^{t/2}\ang s^{-2}\,\rd s
\lesssim\ang t^{-3/2},
\qquad
C\ang t^{-5}\int_0^{t/2}\ang u^{3/2}\,\rd u
\lesssim\ang t^{-5/2},
\]
respectively, where \(u=t-s\) in the second part.  For \(t\le1\),
direct integration suffices.  Hence the history term is bounded by
\(C\ang t^{5/2}\ang{t,x}^{-4}[h_{\rm in}]_{M}\), which
implies \eqref{eq:electric-wave-section5} since
\(\ang t\le\ang{t,x}\).

Finally, choose a smooth frequency cutoff \(\chi_{\rm th}\)
vanishing near zero and equal to one on \(\operatorname{supp}\vartheta\).
Since \(G_{\rm th}=\vartheta(D_x)G\),
\begin{equation}
\label{eq:threshold-field-transfer-section5}
\mathcal P_x\vartheta(D_x)\rho
=\vartheta(D_x)E^{\rm free}
 +G_{\rm th}*_{t,x}[\chi_{\rm th}(D_x)E^{\rm free}].
\end{equation}
The two localized free fields obey
\(C\ang s^{-4}\ang y^{-6}[h_{\rm in}]_{M}\) by
\eqref{eq:band-source-section5}.  Proposition~\ref{prop:threshold-bound-outline}
and \eqref{eq:fixed-scale-transfer-v4}, with
\((a,N,m)=(3/2,4,4)\) in the notation of
Lemma~\ref{lem:weighted-convolution-v4}, therefore yield
\[
|G_{\rm th}*_{t,x}[\chi_{\rm th}(D_x)E^{\rm free}](t,x)|
\le C\frac{\ang t^{5/2}}{\ang{t,x}^4}[h_{\rm in}]_{M}.
\]
Together with the direct bound for \(\vartheta(D_x)E^{\rm free}\),
this controls the threshold term.  Summing the four contributions
proves \eqref{eq:electric-wave-section5}.
\end{proof}

\begin{remark}\label{rem:wave-dipole-v4}
The spatial exponent three in \eqref{eq:main-electric-wave}
cannot in general be increased under neutrality alone.
Choose real \(\varphi,\psi\in C_c^\infty(\R^3)\) with unit
integrals and set
\(h_{\rm in}(x,v)=-(D\cdot\nabla\varphi(x))\psi(v)\),
where \(D\in\R^3\setminus\{0\}\).  Then
\(\rho^{\rm free}(0)=-D\cdot\nabla\varphi\) is neutral and
has dipole moment \(\int y\rho^{\rm free}(0,y)\,\rd y=D\).
With \(N(x)=(4\pi|x|)^{-1}\), the Poisson formula and Taylor
expansion on the fixed support of the density give
\[
E(0,x)=-\int\nabla N(x-y)\rho^{\rm free}(0,y)\,\rd y
=\nabla^2N(x)D+O(|x|^{-4}).
\]
Thus, uniformly in \(\omega\in\mathbb S^2\),
\begin{equation}\label{eq:wave-dipole-tail-v4}
E(0,r\omega)
=\frac{3\omega(\omega\cdot D)-D}{4\pi r^3}+O(r^{-4}).
\end{equation}
The bound \eqref{eq:main-electric-front-v3} gives
\(E_{\rm LD}(0,r\omega)=O(r^{-4})\), so \(E_{\rm wave}\)
has the same leading term.  Since
\( |3\omega(\omega\cdot D)-D|^2=|D|^2+3|\omega\cdot D|^2\),
this datum admits neither a higher uniform spatial power nor a
strong \(L^1_x\) bound for \(E_{\rm wave}(0)\).
\end{remark}

\subsection{The low-frequency Landau-damped component}

We separate the smooth local term from the low-frequency response,
then estimate the latter by its endpoint transform and spherical
expansion.

Let $k_0$ be as in Proposition~\ref{prop:low-frequency-factor-section3}.
Choose $\chi_0\in C_c^\infty([0,k_0))$, equal to one near zero,
with support contained in the set where $\chi_\ell=1$.

\begin{proposition}
\label{prop:corrected-low-electric-section5}
Let
\begin{align}
E_{\rm LD}^{\rm low}
:=\mathcal P_x\Bigl(\chi_0(D_x)[\chi_\ell(D_x)-\widehat K_1(0,D_x)]\rho^{\rm free}
+[\chi_0(R-P)]^\vee*_{t,x}\rho^{\rm free}\Bigr).
\label{eq:low-electric-definition-section5}
\end{align}
If $M\ge11$, then, for every $t\ge0$,
\begin{equation}
\label{eq:low-electric-envelope-section5}
\begin{aligned}
 |E_{\rm LD}^{\rm low}(t,x)|
 \le \frac{C[h_{\rm in}]_M}{\ang{t,x}^{4}}
 \biggl[1+\sum_{\substack{a\in\{\kappa,1\}\\a>0}}
 \biggl\{\left(\frac{\ang{t,x}}{\ang{|x|-at}}\right)^{(2-\gamma_a)_+}
 +\one_{\{\gamma_a=2\}}
 \log\!\left(2+\frac{\ang{t,x}}{\ang{|x|-at}}\right)\biggr\}\biggr].
\end{aligned}
\end{equation}
\end{proposition}

\begin{proof}
\emph{The local term and the response representation.}
The two $k^2$ factors in \eqref{eq:corrected-low-k2-section5}
remove the Poisson singularity on the support of $\chi_0$.
By \eqref{eq:pole-sum-rule-section5}, the local term is
\[
 -[\chi_0b_0](D_x)\nabla_x\rho^{\rm free}
 =-\mathcal K_{\rm loc}*_x\nabla_x\rho^{\rm free},
 \qquad \mathcal K_{\rm loc}=\mathcal F^{-1}[\chi_0b_0]\in\mathcal S(\R^3).
\]
Indeed, the symbol is Cartesian smooth and compactly supported.
Taking $\ell=j=0$, $|\alpha|=1$, and $L=8$ in
\eqref{eq:free-moment-pointwise-section5}, then applying the weighted
spatial convolution bound \eqref{eq:weighted-kernel-convolution-v4},
gives
\[
 |\mathcal K_{\rm loc}*_x\nabla_x\rho^{\rm free}(t,x)|
 \le C\frac{\ang t^4}{\ang{t,x}^8}[h_{\rm in}]_{M}
 \le C\ang{t,x}^{-4}[h_{\rm in}]_{M}.
\]

For the response, put
\[
 j(k,s)=j_0(s)+k^2j_1(k,s),\qquad
 \mathfrak b(t;k,s,\theta)
 =\int_0^{kt}\re^{-\ri su}
       \widehat h_{\rm in}(k\theta,u\theta)\,\rd u.
\]
The free-transport identity
$\widehat{\rho^{\rm free}}(\tau,k\theta)
=\widehat h_{\rm in}(k\theta,k\tau\theta)$ yields
\begin{equation}\label{eq:source-history-v5}
 \frac{\mathfrak b(t;k,s,\theta)}k
 =\int_0^t \re^{-\ri ks\tau}
       \widehat{\rho^{\rm free}}(\tau,k\theta)\,\rd\tau,
 \qquad k>0.
\end{equation}
The endpoint and interior estimates
\eqref{eq:jump-conormal-section3} and
\eqref{eq:jump-interior-low-section3} in
Proposition~\ref{prop:low-frequency-factor-section3}\textup{(i)} imply
\[
 \sup_{0<k\le k_0}\int_{\mathcal I_\kappa}|j(k,s)|\,\rd s\le C.
\]
Part~\textup{(iv)} of the same proposition gives, for
$0<k\le k_0$ and $\Re z>0$,
\[
 R(z,k)-P(z,k)
 =k^3\int_{\mathcal I_\kappa}\frac{j(k,s)}{z-\ri ks}\,\rd s.
\]
The preceding bound makes $k^3\re^{-t\Re z}|j(k,s)|$ integrable in
$(t,s)\in[0,\infty)\times\mathcal I_\kappa$.  Thus Fubini's theorem
applied to the forward Laplace transform, followed by uniqueness, yields
\[
 \Lap^{-1}[\chi_0(k)(R-P)](t)
 =\chi_0(k)k^3\int_{\mathcal I_\kappa}
                 \re^{\ri kst}j(k,s)\,\rd s.
\]
Its time convolution with the free density and the Poisson multiplier
$-\ri\theta/k$ therefore give
\[
\begin{aligned}
 \widehat E_{\rm LD}^{\rm low,resp}(t,k\theta)
 &=-\ri\theta\chi_0(k)k^2\int_{\mathcal I_\kappa}
      \re^{\ri kst}j(k,s)
      \left[\int_0^t \re^{-\ri ks\tau}
      \widehat{\rho^{\rm free}}(\tau,k\theta)\,\rd\tau\right]\rd s\\
 &=-\ri\theta\chi_0(k)k\int_{\mathcal I_\kappa}
      \re^{\ri kst}j(k,s)\mathfrak b(t;k,s,\theta)\,\rd s.
\end{aligned}
\]
These finite-time integrals are absolutely convergent.  Spatial Fourier
inversion with $\xi=k\theta$ and $\rd\xi=k^2\rd k\,\rd\theta$ now yields
\begin{equation}
\label{eq:low-electric-real-integral-section5}
\begin{aligned}
 E_{\rm LD}^{\rm low,resp}(t,x)
 =-\frac{\ri}{(2\pi)^3}
 \int_{\mathbb S^2}\int_0^\infty\int_{\mathcal I_\kappa}
 \re^{\ri k(x\cdot\theta+ts)}\theta\chi_0(k)k^3j(k,s)\mathfrak b(t;k,s,\theta)\,\rd s\,\rd k\,\rd\theta.
\end{aligned}
\end{equation}
For $m+j+|\alpha|\le10$, Lemma~\ref{lem:truncated-real-source-appendix}
gives
\begin{equation}
\label{eq:truncated-source-symbol-section5}
 \left|D_k^m\partial_s^j\nabla_\theta^\alpha
       \mathfrak b(t;k,s,\theta)\right|
 \le C\min\{kt,1\}[h_{\rm in}]_{M},
 \qquad s\in\mathcal I_\kappa.
\end{equation}
Here and below $D_k=k\partial_k$ acts at fixed $t$, including on
the upper limit $kt$.

\medskip\noindent
\emph{Low frequencies.}
Absolute integration in \eqref{eq:low-electric-real-integral-section5}
first gives a uniform bound on $t+r\le6$.  Since $\ang{t,x}^{-4}$
is bounded below there, this proves \eqref{eq:low-electric-envelope-section5}
in that region.  In the remaining
region $t+r>6$, the very low frequencies satisfy
\begin{equation}
\label{eq:low-electric-small-k-section5}
 \int_0^{C/\ang{t,x}}k^3\min\{kt,1\}\,\rd k
 \lesssim t\ang{t,x}^{-5}\lesssim\ang{t,x}^{-4}.
\end{equation}
We used the boundedness of $j$ and
\eqref{eq:truncated-source-symbol-section5}; the field estimate carries
the factor $[h_{\rm in}]_{M}$.
Choose $c>0$ sufficiently large, depending only on the fixed endpoint
charts, and choose a smooth $\psi$ with
\[
 \psi(\lambda)=0\quad(\lambda\le c),\qquad
 \psi(\lambda)=1\quad(\lambda\ge2c).
\]
Insert $\psi(\ang{t,x}k)$ into the remaining integral.  Its complement
is covered by \eqref{eq:low-electric-small-k-section5}, and all
logarithmic derivatives of this cutoff are uniformly bounded.
If $c/\ang{t,x}\ge k_0$, the remaining frequency integral vanishes.
We may therefore assume $c/\ang{t,x}<k_0$ below.

\medskip\noindent
\emph{Away from the fronts.}
For $r>2t$ and every $s\in\mathcal I_\kappa$, the phase
\[
 \Phi_s(\xi)=x\cdot\xi+ts|\xi|,
 \qquad |\nabla_\xi\Phi_s|\ge r-t\ge r/2,
\]
has no stationary point.  The Cartesian amplitude in
\eqref{eq:low-electric-real-integral-section5} is of order $k$;
its derivatives are bounded by
$C_\sigma k^{1-|\sigma|}[h_{\rm in}]_{M}$ for every Cartesian
multi-index $\sigma\in\mathbb N_0^3$ with $|\sigma|\le5$, by
Proposition~\ref{prop:low-frequency-factor-section3} and
\eqref{eq:truncated-source-symbol-section5}.
Five integrations by parts using
\[
 L_s=\frac{\nabla_\xi\Phi_s\cdot\nabla_\xi}
                {\ri|\nabla_\xi\Phi_s|^2},
 \qquad L_s\re^{\ri\Phi_s}=\re^{\ri\Phi_s},
\]
thus give
\[
 C[h_{\rm in}]_{M}\,r^{-5}
       \int_{c/\ang{t,x}}^{k_0}k^{-2}\,\rd k
 \lesssim\ang{t,x}^{-4}[h_{\rm in}]_{M}.
\]
The cutoff derivatives obey the same bounds and remove the radial
boundary terms.  This treats the full $s$-integral in the exterior region.

Suppose next that $r\le2t$.  Since $t+r>6$, one has
$t>2$ and $t\asymp\ang{t,x}$.  Use a smooth partition of
$\mathcal I_\kappa$ into endpoint charts and a compact interior part.
For the interior part, five integrations by parts in $s$ yield
\[
 C[h_{\rm in}]_{M}\,t^{-5}
       \int_{c/\ang{t,x}}^{k_0}k^{-2}\,\rd k
 \lesssim\ang{t,x}^{-4}[h_{\rm in}]_{M}.
\]
For each positive endpoint $a$, choose a smaller chart
$|s-a|<\eta_a$ inside the corresponding chart of
Proposition~\ref{prop:low-frequency-factor-section3}, with
$0<\eta_a<a/8$.  If
$|r-at|>2\eta_at$, then throughout that chart
\[
 |\nabla_\xi\Phi_s|
 \ge |r-ts|\ge |r-at|-t|s-a|\ge\eta_at.
\]
The preceding five integrations in $\xi$ give the same bound
$Ct^{-5}\int_{c/\ang{t,x}}^{k_0}k^{-2}\,\rd k$
times the data norm.  For a negative endpoint the identical argument
uses $|s|$ in place of $s$.  These estimates include $r=0$.

\medskip\noindent
\emph{Endpoint transforms and spherical contributions.}
It remains to consider $|r-at|\le2\eta_at$, where
$r\asymp at\asymp t\asymp\ang{t,x}$ for the fixed $a>0$.
Our choice of $c$ ensures $tk\ge1$ and $rk\ge1$ on the remaining
frequency support.  We first transform the endpoint variable, then
apply the spherical expansion to the resulting angular amplitude.

Write $s=a+\epsilon_a u$ in its inward coordinate.  By
Proposition~\ref{prop:low-frequency-factor-section3}\textup{(i)},
\begin{equation}\label{eq:electric-endpoint-j-section5}
 |D_k^{j_1}\partial_u^{j_2}j(k,a+\epsilon_a u)|
 \le C_{j_1,j_2}u^{\gamma_a-j_2},\qquad j_1+j_2\le6,
 \quad 0<k\le k_0.
\end{equation}
Choose $\chi_{\rm ep}\in C_c^\infty([0,\eta_a))$ equal to one for
$0\le u\le\eta_a/2$.  Set $\beta=\min\{\gamma_a,5/2\}$ and define
\[
 f_a(t;k,u,\theta)
 =j(k,a+\epsilon_a u)\,\theta\,
       \mathfrak b(t;k,a+\epsilon_a u,\theta)\chi_{\rm ep}(u).
\]
Its normalized endpoint transform is
\begin{equation}\label{eq:electric-whole-transform-v3}
 Q_a(t;k,\theta)
 =(tk)^{\beta+1}\int_0^{\eta_a}
       \re^{\ri\epsilon_atku}f_a(t;k,u,\theta)\,\rd u.
\end{equation}
We apply Lemma~\ref{lem:finite-endpoint-v4} with $L=4>\beta+1$.
For the derivative ranges below, the transfer identity
\eqref{eq:endpoint-transfer-v3} requires $0\le j_1\le m$ and
$0\le j_2\le4+m-j_1$.  Thus $j_1+j_2\le6$ and $j_1+j_2+|\alpha|\le10$,
so the required bounds follow from
\eqref{eq:electric-endpoint-j-section5} and
\eqref{eq:truncated-source-symbol-section5}, the latter supplied by
Lemma~\ref{lem:truncated-real-source-appendix}.
The large- and small-parameter estimates of
Lemma~\ref{lem:finite-endpoint-v4} give
\begin{equation}\label{eq:electric-whole-budget-v3}
\begin{aligned}
 \|D_k^m\nabla_\theta^\alpha Q_a(t;k,\cdot)\|_{L^\infty_\theta}
 \le C\min\{1,(tk)^{\beta+2}\}[h_{\rm in}]_{M},\quad
 \text{if }0\le m\le2,\ |\alpha|\le2,
       \quad\text{or if }m=0,\ |\alpha|\le6.
\end{aligned}
\end{equation}
The additional factor $tk$ for $tk\le1$ comes from
$\min\{tk,1\}$ in \eqref{eq:truncated-source-symbol-section5}.

With $\omega=x/r$, the contribution of the positive endpoint is
\[
\begin{aligned}
 E_a^{\rm end}(t,r\omega)
 =-\frac{\ri}{(2\pi)^3t^{\beta+1}}
 \int_0^\infty \re^{\ri atk}\psi(\ang{t,x}k)\chi_0(k)k^{2-\beta}
 \times\left(\int_{\mathbb S^2}\re^{\ri kr\omega\cdot\theta}
             Q_a(t;k,\theta)\,\rd\theta\right)\rd k.
\end{aligned}
\]
Formula \eqref{eq:spherical-carriers-v4} gives
$E_a^{\rm end}=E_a^{(1)}+E_a^{(2)}+\mathcal R_a$, where
\begin{equation}\label{eq:electric-front-reduction-section5}
\begin{aligned}
 E_a^{(1)}(t,r\omega)
 =-\frac{1}{4\pi^2rt^{\beta+1}}
 \int_0^\infty &\psi(\ang{t,x}k)\chi_0(k)k^{1-\beta}\\
 &\times\bigl[\re^{\ri k(r+at)}Q_a(t;k,\omega)
       -\re^{-\ri k(r-at)}Q_a(t;k,-\omega)\bigr]\,\rd k,
\end{aligned}
\end{equation}
and
\begin{equation}\label{eq:electric-second-carrier-v4}
\begin{aligned}
 E_a^{(2)}(t,r\omega)
 =\frac{\ri}{8\pi^2r^2t^{\beta+1}}
 \int_0^\infty &\psi(\ang{t,x}k)\chi_0(k)k^{-\beta}\\
 &\times\bigl[\re^{\ri k(r+at)}\Delta_{\mathbb S^2}Q_a(t;k,\omega)+\re^{-\ri k(r-at)}\Delta_{\mathbb S^2}Q_a(t;k,-\omega)\bigr]\,\rd k.
\end{aligned}
\end{equation}
For $E_a^{(1)}$, the endpoint transform supplies $(tk)^{-\beta-1}$,
leaving $k^{1-\beta}$ in the radial Fourier integral.
Lemma~\ref{lem:truncated-radial-v4}, with $J=2$ and $\alpha=2-\beta$,
bounds the contribution with phase $-k(r-at)$ by
\[
 C[h_{\rm in}]_{M}
 \begin{cases}
 \ang{t,x}^{-\gamma_a-2}\ang{r-at}^{\gamma_a-2}
       +\ang{t,x}^{-4},&1\le\gamma_a<2,\\[1mm]
 \ang{t,x}^{-4}\log\!\left(2+
       \dfrac{\ang{t,x}}{\ang{r-at}}\right),&\gamma_a=2,\\[1mm]
 \ang{t,x}^{-4},&\gamma_a>2.
 \end{cases}
\]
For the phase $k(r+at)$, the same lemma and
$\ang{r+at}\asymp\ang{t,x}$ give
$C\ang{t,x}^{-4}[h_{\rm in}]_{M}$.
Together, these estimates give the single-endpoint contribution
on the right-hand side of \eqref{eq:low-electric-envelope-section5}.

For $E_a^{(2)}$, the radial power is $k^{-\beta}$.
The same lemma with $\alpha=1-\beta$ yields
\[
 |E_a^{(2)}(t,r\omega)|\le C[h_{\rm in}]_{M}
 \begin{cases}
 \ang{t,x}^{-4}\log\!\left(2+
       \dfrac{\ang{t,x}}{\ang{r-at}}\right),&\gamma_a=1,\\[1mm]
 \ang{t,x}^{-4},&\gamma_a>1.
 \end{cases}
\]
At $\gamma_a=1$ the logarithm is bounded by
$C\ang{t,x}/\ang{r-at}$, so this term satisfies the same upper bound
obtained above for $E_a^{(1)}$. At $\gamma_a=2$ the logarithm in the first contribution
remains explicitly in \eqref{eq:low-electric-envelope-section5}.

The spherical remainder uses only the six angular derivatives in
\eqref{eq:electric-whole-budget-v3}.  Its absolute estimate is
\begin{equation}\label{eq:electric-third-angle-v4}
 |\mathcal R_a(t,r\omega)|
 \le \frac{C[h_{\rm in}]_{M}}{r^3t^{\beta+1}}
       \int_{c/\ang{t,x}}^{k_0}k^{-\beta-1}\,\rd k
 \le C\ang{t,x}^{-4}[h_{\rm in}]_{M}.
\end{equation}
This estimate does not differentiate the spherical remainder along
$\lambda=kr$.  The negative endpoint has phases $k(r-at)$ and
$-k(r+at)$, with the two spherical poles interchanged, and obeys the
same bounds.  Summing the endpoint contributions and the preceding
local, low-frequency, and nonfront estimates proves
\eqref{eq:low-electric-envelope-section5}.  The bounded region was
already covered by absolute integration.
\end{proof}

\subsection{The fixed-gap component}

We now treat all frequencies outside the preceding low-frequency
block.  Define
\begin{align}
R_{\rm fg}
&:=(1-\chi_0)\chi_\ell(R-P)+\chi_h R,
\label{eq:Rfg-section5}\\
a_{\rm fg}(k)
&:=(1-\chi_0(k))\bigl[\chi_\ell(k)-\widehat K_1(0,k)\bigr]+\chi_h(k).
\label{eq:afg-section5}
\end{align}
Both vanish on a neighborhood of \(k=0\).  The corresponding part of
the Landau-damped field is exactly
\begin{equation}
\label{eq:Efg-section5}
E_{\rm LD}^{\rm fg}
=a_{\rm fg}(D_x)E^{\rm free}
+G_{\rm fg}*_{t,x}E^{\rm free},
\qquad G_{\rm fg}:=(R_{\rm fg})^\vee.
\end{equation}

The fixed-gap kernel transfers its front bounds to the electric field.

\begin{lemma}
\label{lem:fixed-gap-transfer-section5}
For \(t\ge1\),
\begin{equation}\label{eq:Gfg-section5}
 |G_{\rm fg}(t,x)|\le \frac{C}{\ang{t,x}^{6}}
 \left[1+\sum_{\substack{a\in\{\kappa,1\}\\a>0}}
 \left(\frac{\ang{t,x}}{\ang{|x|-at}}\right)^{(4-\gamma_a)_+}\right].
\end{equation}
For every fixed \(m\ge0\),
\begin{equation}\label{eq:Gfg-short-weighted-v5}
 \sup_{0<t\le1}\|\ang x^mG_{\rm fg}(t)\|_{L_x^1}\le C_m.
\end{equation}
For \(M\ge11\) and \(t\ge0\), the
corresponding electric component satisfies
\begin{equation}\label{eq:Efg-bound-section5}
\begin{aligned}
 |E_{\rm LD}^{\rm fg}(t,x)|
 \le \frac{C[h_{\rm in}]_M}{\ang{t,x}^{4}}
 \biggl[1+\sum_{\substack{a\in\{\kappa,1\}\\a>0}}
 \biggl\{\left(\frac{\ang{t,x}}{\ang{|x|-at}}\right)^{(2-\gamma_a)_+}+\one_{\{\gamma_a=2\}}
 \log\!\left(2+\frac{\ang{t,x}}{\ang{|x|-at}}\right)\biggr\}\biggr].
\end{aligned}
\end{equation}
\end{lemma}
\begin{proof}
The proof of Proposition~\ref{prop:green-LD-outline} estimates the
frequency and endpoint pieces before summation.  Inserting the fixed
multiplier \(1-\chi_0\) preserves those symbol bounds and gives
\eqref{eq:Gfg-section5}.  For short propagation times, the exact
decomposition gives
\begin{align}
 G_{\rm fg}=(1-\chi_0(D_x))G_{\rm LD}
 =G-\chi_0(D_x)G-G_{\rm KG}+\chi_0(D_x)G_{\rm KG}
                         -\vartheta(D_x)G,
                         \label{eq:Gfg-short-decomposition-v5}
\end{align}
where \(\chi_0\vartheta=0\).  Lemma~
\ref{lem:finite-propagation-green-section4} gives
\(\sup_{0<t\le1}\|\ang x^mG(t)\|_1\le C_m\).
Proposition~\ref{prop:pole-dispersive-outline} gives the same bound
for \(G_{\rm KG}\), by taking \(N>m+3\) in
\eqref{eq:KG-weighted-section4}.  The inverse Fourier transforms of
\(\chi_0\) and \(\vartheta\) are Schwartz functions.  For any
\(m\ge0\), Peetre's inequality and Fubini give
\[
 \|\ang x^m(\mathcal K*_xF)\|_1
 \le C_m\|\ang x^m\mathcal K\|_1\|\ang x^mF\|_1,
\]
since \(\ang x^m\le C_m\ang{x-y}^m\ang y^m\).
Applying this inequality to \eqref{eq:Gfg-short-decomposition-v5}
proves \eqref{eq:Gfg-short-weighted-v5}.

The scalar bound \eqref{eq:Gfg-section5} also implies
\begin{equation}\label{eq:Gfg-integrable-interfaces-v3}
 \|G_{\rm fg}(u)\|_\infty\le C\ang u^{-3},\qquad
 \|G_{\rm fg}(u)\|_1\le C,\qquad u\ge1.
\end{equation}
Indeed, when $|x|\le2u$, the right-hand side of
\eqref{eq:Gfg-section5} is $O(\ang u^{-3})$ because
$\gamma_a\ge1$ and $\ang{u,x}\asymp\ang u$.  When $|x|>2u$,
$\ang{|x|-au}\asymp\ang{u,x}$ for every $0<a\le1$, so it is
$O(\ang{u,x}^{-6})$.  In both regions it is bounded by
$C\ang u\,\ang{u,x}^{-4}$, whose spatial integral is uniformly
bounded.  Choose a smooth high-pass symbol
\(\chi_{\rm fg}\), vanishing near zero and equal to one on the
frequency supports of \(a_{\rm fg}\) and \(R_{\rm fg}\),
with symbol bounds through order six, as required in Proposition~\ref{lem:free-moments-section5}.
It equals one at all sufficiently high frequencies.  Equation
\eqref{eq:band-source-section5} yields
\[
 |a_{\rm fg}(D_x)E^{\rm free}(s,y)|
 +|\chi_{\rm fg}(D_x)E^{\rm free}(s,y)|
 \le C\ang s^{-4}\ang y^{-6}[h_{\rm in}]_{M}.
\]
Thus \eqref{eq:Efg-section5} gives the bound
\[
\begin{split}
 |E_{\rm LD}^{\rm fg}(t,x)|\le C[h_{\rm in}]_{M}
 \biggl[\ang t^{-4}\ang x^{-6}
 +\int_0^t\ang s^{-4}\int_{\R^3}
 |G_{\rm fg}(t-s,x-y)|\ang y^{-6}\,\rd y\,\rd s\biggr].
\end{split}
\]
The first term is at most
\(C\ang{t,x}^{-4}[h_{\rm in}]_{M}\); we estimate the
integral in three regions.

\emph{Bounded space--time.}
If \(t+r\le6\), \eqref{eq:Gfg-short-weighted-v5} and
\eqref{eq:Gfg-integrable-interfaces-v3} give
\(\sup_{u>0}\|G_{\rm fg}(u)\|_1\le C\).
The integral is uniformly bounded, and \(\ang{t,x}\asymp1\),
proving \eqref{eq:Efg-bound-section5} in this region.

\emph{The region \(r\le2t\).}
We may assume \(t+r>6\), so that \(t>2\) and
\(\ang{t,x}\asymp\ang t\).  On \(0\le s\le t/2\),
\(t-s\ge t/2\) and hence
\(\ang{t-s,x-y}\gtrsim\ang{t,x}\).
The bulk part of \eqref{eq:Gfg-section5} contributes at most
\(C\ang{t,x}^{-6}[h_{\rm in}]_{M}\).
For an endpoint with \(1\le\gamma_a<2\), the triangle inequality
\[
 |r-at|\le |y|+\bigl||x-y|-a(t-s)\bigr|+as
\]
implies
\[
 \ang{|x-y|-a(t-s)}^{-(4-\gamma_a)}
 \le C\ang{r-at}^{-(2-\gamma_a)}
       \ang s^{2-\gamma_a}\ang y^{2-\gamma_a}.
\]
By \eqref{eq:Gfg-section5}, the corresponding contribution is
therefore bounded by
\begin{align*}
 C[h_{\rm in}]_{M}\ang{t,x}^{-\gamma_a-2}
       \ang{r-at}^{-(2-\gamma_a)}
       \int_0^{t/2}\ang s^{-2-\gamma_a}\,\rd s
       \int_{\R^3}\ang y^{-4-\gamma_a}\,\rd y
      \le C[h_{\rm in}]_{M}
       \ang{t,x}^{-\gamma_a-2}\ang{r-at}^{\gamma_a-2}.
\end{align*}
Both integrals converge, including at \(\gamma_a=1\).
For \(\gamma_a\ge2\), the scalar front term is bounded by
\(C\ang{t-s,x-y}^{-4}\), so its contribution is at most
\(C\ang{t,x}^{-4}[h_{\rm in}]_{M}\).

On \(t/2\le s\le t\), first take \(t-s\ge1\).  The
\(L^\infty\) bound in \eqref{eq:Gfg-integrable-interfaces-v3}
and the source's spatial integral give
\[
 C[h_{\rm in}]_{M}\ang t^{-4}
       \int_1^{t/2}\ang u^{-3}\,\rd u
 \le C\ang{t,x}^{-4}[h_{\rm in}]_{M}.
\]
For \(0<t-s\le1\), \eqref{eq:Gfg-short-weighted-v5} with
\(m=6\) and \eqref{eq:weighted-kernel-convolution-v4} give
\[
 \int_{\R^3}|G_{\rm fg}(t-s,x-y)|\ang y^{-6}\,\rd y
 \le C\ang x^{-6}.
\]
The remaining time interval has length at most one, and its
contribution is at most
\(C\ang t^{-4}\ang x^{-6}[h_{\rm in}]_{M}\).
This proves \eqref{eq:Efg-bound-section5} throughout \(r\le2t\).

\emph{The region \(r>2t\).}
We again assume \(t+r>6\), so that \(r>4\) and
\(\ang{t,x}\asymp r\).  When \(t-s\ge1\), split the spatial
integral at \(|y|=r/4\).  On \(|y|\le r/4\), for every
\(0<a\le1\),
\[
 |x-y|-a(t-s)\ge r-|y|-t\ge r/4,
 \qquad \ang{t-s,x-y}\asymp r.
\]
Thus \eqref{eq:Gfg-section5} gives
\(|G_{\rm fg}(t-s,x-y)|\le Cr^{-6}\).
Integration against \(\ang s^{-4}\ang y^{-6}\) contributes
\(Cr^{-6}[h_{\rm in}]_{M}\).
On \(|y|>r/4\), use \(\ang y^{-6}\le Cr^{-6}\) and the
\(L^1\) bound in \eqref{eq:Gfg-integrable-interfaces-v3} to obtain
the same estimate.  For \(0<t-s\le1\), the preceding weighted
convolution bound gives \(C\ang x^{-6}\), and the time weight is
integrable.  Therefore
\( |E_{\rm LD}^{\rm fg}(t,x)|\le Cr^{-6}[h_{\rm in}]_{M}\),
which proves \eqref{eq:Efg-bound-section5} in the exterior region.
\end{proof}

\subsection{Completion of Theorem~\ref{thm:main-electric}}

\begin{proof}[Proof of Theorem~\ref{thm:main-electric}]
The decomposition \eqref{eq:electric-exact-sum-section5},
Propositions~\ref{prop:electric-wave-section5}
and~\ref{prop:corrected-low-electric-section5}, and
Lemma~\ref{lem:fixed-gap-transfer-section5} give
\eqref{eq:main-electric-wave} and \eqref{eq:main-electric-front-v3}.
To obtain \eqref{eq:main-electric-gamma-v3}, first suppose that
$r\le2t$, so that $\ang{t,x}\asymp\ang t$ and
$\ang{r-at}\ge1$.  Each power in the front bound is at most
$C\ang t^{(2-\gamma_{\min})_+}$.  If an endpoint has
$\gamma_a=2$, its logarithm is bounded by $C\log(2+\ang t)$;
when $\gamma_{\min}<2$ this is absorbed by
$C\ang t^{2-\gamma_{\min}}$.
For $r>2t$, one has $\ang{r-at}\asymp\ang{t,x}$ for every
$0<a\le1$, so the front bound is $C\ang{t,x}^{-4}[h_{\rm in}]_{M_*}$.
These two regions prove \eqref{eq:main-electric-gamma-v3}, including
its logarithmic case $\gamma_{\min}=2$.

The derivative requirements of \eqref{eq:band-source-section5},
\eqref{eq:electric-whole-budget-v3}, and
Lemmas~\ref{lem:velocity-ray-source-v6}
and~\ref{lem:truncated-real-source-appendix} hold with \(M_*=11\).
Under \eqref{H0}, Lemma~\ref{lem:A-endpoint-section2} justifies the
same estimates with only the outer front $a=1$.

It remains to prove the spatial norm estimates.
Rescaling $x=\ang t\,y$ in \eqref{eq:main-electric-wave} gives
the wave bounds in \eqref{eq:main-electric-global} and
\eqref{eq:E-weak-L1-section6}, since $\ang y^{-3}\in L^p$ for
$p>1$ and $\ang y^{-3}\in L^{1,\infty}$.
For $\gamma_{\min}\ge1$, the time factor in square brackets in
\eqref{eq:main-electric-gamma-v3} is bounded by $C\ang t$.
The same rescaling, with $\ang y^{-4}\in L^1\cap L^\infty$,
therefore yields
\[
 \|E_{\rm LD}(t)\|_{L_x^p}
 \le C_p\ang t^{-3+3/p}[h_{\rm in}]_{M_*},
 \qquad 1\le p\le\infty.
\]
Since $E=E_{\rm wave}+E_{\rm LD}$ and $\ang t\ge1$, the triangle
inequality proves \eqref{eq:main-electric-global}.
The $p=1$ case of the last bound, the inclusion
$L^1\subset L^{1,\infty}$, and the weak $L^1$ quasi-triangle
inequality prove \eqref{eq:E-weak-L1-section6}.
\end{proof}

\section{Scattering and relativistic applications}
\label{sec:consequences}
We deduce scattering from the electric-field decay in
Theorem~\ref{thm:main-electric}, then establish the relativistic
applications through the momentum--velocity correspondence.

\subsection{Scattering}

\begin{proof}[Proof of Corollary~\ref{cor:main-scattering}]
Using the almost-everywhere gradient convention of
Section~\ref{subsec:volterra-construction-v3}, the characteristic form of
\eqref{eq:linear-VP} is, for almost every
\((x,v)\),
\begin{equation}
\label{eq:characteristic-formula-section6}
h(t,x+tv,v)
=h_{\rm in}(x,v)
-\int_0^tE(s,x+sv)\cdot\nabla_v f^0(v)\,\rd s .
\end{equation}
By the case $p=\infty$ of \eqref{eq:main-electric-global},
\[
\int_0^\infty\|E(s)\|_{L_x^\infty}\,\rd s
\leq C[h_{\rm in}]_{{11}}
\int_0^\infty\ang s^{-3/2}\,\rd s
\leq C[h_{\rm in}]_{{11}}.
\]
Hence the integral in \eqref{eq:h-infinity} converges absolutely,
uniformly in \(x\) and essentially uniformly in \(v\), and defines an
\(L^\infty_{x,v}\) equivalence class.  Subtracting
\eqref{eq:h-infinity} from \eqref{eq:characteristic-formula-section6}
and using Theorem~\ref{thm:main-electric} gives
\begin{align*}
\operatorname*{ess\,sup}_{x,v}
\bigl|h(t,x+tv,v)-h_\infty(x,v)\bigr|
\leq
\|\nabla_v f^0\|_{L_v^\infty}
\int_t^\infty\|E(s)\|_{L_x^\infty}\,\rd s
\leq C\ang t^{-1/2}[h_{\rm in}]_{{11}}.
\end{align*}
\end{proof}

\subsection{Relativistic applications}
\label{subsec:relativistic-correspondence}

We prove Theorems~\ref{thm:relativistic-green} and
\ref{thm:relativistic-electric} for the background constructed in
\eqref{eq:relativistic-general-profile}.  We first check its regularity
and tail, then identify the linearized equation and its response.
Throughout this subsection $v\in\R^3$ is momentum and
$u=v/\ang v$ is velocity.

Set \(F(\varpi)=\mu(\varpi^2/2)\) and
\[
 W(\varpi):=-\int_\varpi^1(1-\sigma^2)F'(\sigma)\,\rd\sigma,
 \qquad 0\le\varpi<1,
\]
so that \(\Phi(v)=W(|v|/\ang v)\).  Integration by parts gives
\[
 W(\varpi)=(1-\varpi^2)F(\varpi)
               -2\int_\varpi^1\sigma F(\sigma)\,\rd\sigma.
\]
For a genuine shell satisfying \eqref{H1}--\eqref{H2}, this would
equal \(-2\int_\kappa^1\sigma F(\sigma)\,\rd\sigma<0\) for
\(0\le\varpi\le\kappa\).  Thus this matching construction
cannot produce a nonnegative relativistic background from the
shell class.  We work with \eqref{H0} and the assumed
nonnegativity of \(\Phi\).

The change of variables \(e=\sigma^2/2\) gives
\[
 \Phi(v)=-\int_{\frac{|v|^2}{2(1+|v|^2)}}^{1/2}
                          (1-2e)\mu'(e)\,\rd e.
\]
By \eqref{H0}, this is a smooth function of \(|v|^2\) near zero,
so \(\Phi\) is smooth in Cartesian coordinates at \(v=0\).
The chain rule yields
\[
 \nabla_v\Phi(v)=v\ang v^{-6}
          \mu'\!\left(\frac{|v|^2}{2\ang v^2}\right).
\]
Near the outer boundary write
$F(\varpi)=(1-\varpi)^{\gamma_1}b(\varpi)$, where $b(1)>0$.
Direct integration in \eqref{eq:relativistic-general-profile} gives
\[
 W(\varpi)=\frac{2\gamma_1b(1)}{\gamma_1+1}
 (1-\varpi)^{\gamma_1+1}\bigl(1+O(1-\varpi)\bigr).
\]
Using $1-\varpi=\ang v^{-2}/(1+\varpi)$, we obtain
\[
 \Phi(v)=\frac{\gamma_1b(1)}{2^{\gamma_1}(\gamma_1+1)}
 \ang v^{-2\gamma_1-2}\bigl(1+O(\ang v^{-2})\bigr),
 \qquad |v|\longrightarrow\infty.
\]
Together with the assumed nonnegativity, this proves that
$\Phi$ is an integrable radial background for every $\gamma_1\ge1$.

The inverse velocity map and its Jacobian are
\[
 v(u)=\frac{u}{\sqrt{1-|u|^2}},\qquad
 J(u)=\det D_uv(u)=(1-|u|^2)^{-5/2},\qquad
 \rd v=J(u)\,\rd u,\qquad |u|<1.
\]
Indeed, $D_uv$ has two tangential eigenvalues
$(1-|u|^2)^{-1/2}$ and one radial eigenvalue
$(1-|u|^2)^{-3/2}$.  Together with the gradient formula above,
this gives
\begin{equation}\label{eq:relativistic-coupling-section6}
 J(u)\nabla_v\Phi(v(u))
 =u\mu'(|u|^2/2)=\nabla_uf^0(u),\qquad |u|<1.
\end{equation}
At $|u|=1$, the weak gradient of the zero extension
has no surface term because $F(1)=0$, including when $\gamma_1=1$.
In particular,
$\int|\nabla_v\Phi|\,\rd v
=\int|\nabla_uf^0|\,\rd u<\infty$.

Define $g(t,x,u)=J(u)h_{\rm rel}(t,x,v(u))$ in the unit ball and
extend it by zero outside. Multiplying \eqref{eq:linear-relativistic-VP-intro} by $J$ and using \eqref{eq:relativistic-coupling-section6}, we obtain
\[
 \partial_tg+u\cdot\nabla_xg+E_{\rm rel}\cdot\nabla_uf^0=0,
 \qquad
 \int_{\R^3}g(t,x,u)\,\rd u
 =\int_{\R^3}h_{\rm rel}(t,x,v)\,\rd v.
\]
Thus $g$ satisfies \eqref{eq:linear-VP} with equilibrium $f^0$
and the same density and electric field.
The free densities agree as well:
\[
 \int_{\R^3}h_{{\rm rel,in}}(x-tv/\ang v,v)\,\rd v
 =\int_{\R^3}g_{\rm in}(x-tu,u)\,\rd u.
\]
Consequently, for $\xi\ne0$,
\begin{equation}\label{eq:relativistic-kernel-section6}
\begin{aligned}
 \widehat K_{\rm rel}(t,\xi)
 =\frac{\ri\xi}{|\xi|^2}\cdot
   \int_{\R^3}\re^{-\ri t\xi\cdot v/\ang v}
                    \nabla_v\Phi(v)\,\rd v
 =\frac{\ri\xi}{|\xi|^2}\cdot
   \int_{|u|<1}\re^{-\ri t\xi\cdot u}
                    \nabla_uf^0(u)\,\rd u
 =-t\widehat f^0(t\xi).
\end{aligned}
\end{equation}
The last equality is weak integration by parts in $u$.
This is the classical Fourier kernel \eqref{eq:kernel-K}, with continuous
value $\widehat K_{\rm rel}(t,0)=-t$ because $\int f^0=1$.
Finite-time uniqueness established after
\eqref{eq:volterra-series-section2} identifies the same causal
response $G$.
Theorem~\ref{thm:main-green} now proves
Theorem~\ref{thm:relativistic-green}.

Under the assumptions on the zero extension of
\eqref{eq:relativistic-data-intro} in
Theorem~\ref{thm:relativistic-electric}, the equality of free density
and response identifies the complete fields and every component of
the normal form \eqref{eq:rho-wave-intro}--\eqref{eq:electric-components-intro}.
Theorem~\ref{thm:main-electric} therefore gives its stated field
bounds with $[g_{\rm in}]_{11}$ on the right-hand sides.
The data condition is imposed on the zero extension of
$g_{\rm in}$ in \eqref{eq:data-norm}.  Real, neutral, smooth data
compactly supported in $(x,v)$ give a sufficient class.
For the scattering state $g_\infty$ from
Corollary~\ref{cor:main-scattering}, put
\[
 h_{{\rm rel},\infty}(x,v)
 :=J(v/\ang v)^{-1}g_\infty(x,v/\ang v).
\]
The diffeomorphism $v\leftrightarrow u$ preserves null sets and
$J^{-1}\le1$, so
\begin{equation}\label{eq:relativistic-scattering-section6}
 \operatorname*{ess\,sup}_{x,v}
 \bigl|h_{\rm rel}(t,x+tv/\ang v,v)-h_{{\rm rel},\infty}(x,v)\bigr|
 \le C\ang t^{-1/2}[g_{\rm in}]_{11},\qquad t\ge0.
\end{equation}
This finishes the proof of Theorem~\ref{thm:relativistic-electric}.

For $F(\varpi)=c_\gamma(1-\varpi^2)^\gamma$, direct integration gives
$W(\varpi)=\frac{\gamma}{\gamma+1}c_\gamma
(1-\varpi^2)^{\gamma+1}$.  Substitution of $\varpi=|v|/\ang v$
therefore gives \eqref{eq:relativistic-profile-section6} and proves
Corollary~\ref{cor:relativistic-power-profiles}.
The normalization fixes the mass of the effective classical
background, which need not equal the relativistic mass.

\appendix

\section{Oscillatory integral estimates}
\label{app:endpoints}

This appendix collects dispersive, spherical, endpoint, and radial
oscillatory estimates for the Green kernel and the electric field.
We use standard stationary- and nonstationary-phase methods
\cite{Hor03,Stein93}, recording the finite derivative bounds and
parameter uniformity needed below.

\subsection*{Cartesian and polar dispersive estimates}

\begin{lemma}
\label{lem:radial-dispersive-section4}
Let $\omega$ be smooth and even at zero.  Suppose the radial and
tangential Hessian eigenvalues $\omega''(k)$ and $\omega'(k)/k$ have
the same sign and absolute value at least $c_0>0$ on a neighborhood
of $0\le k\le k_1$.
Write $k=|\xi|$ and $\theta=\xi/|\xi|$ for $\xi\ne0$.
Both amplitudes $A(\xi)$ and $a(k,\theta)$ below are supported in
$k\le k_1$.  For every fixed integer $N\ge0$ and
$A\in C_c^{N+2}(\R^3)$,
\begin{equation}\label{eq:radial-dispersive-section4}
 \left|\int_{\R^3}\re^{\ri (x\cdot\xi+t\omega(|\xi|))}A(\xi)\,\rd\xi\right|
 \le C_N\|A\|_{C^{N+2}}\ang t^{-3/2}
                    \frac{\ang t^N}{\ang{t,x}^N}.
\end{equation}
For the second estimate, $a(k,\theta)$ need not extend smoothly in
Cartesian coordinates through $\xi=0$; its value at that single
point does not affect the integral.  If it has finite polar seminorm
\[
 \|a\|_{\mathrm{pol},6}
 =\max_{j+|\beta|\le6}\sup_{k>0,\theta}
       |(k\partial_k)^j\nabla_\theta^\beta a(k,\theta)|<\infty,
\]
then
\begin{equation}\label{eq:polar-pole-source-section5}
 \left|\int_{\R^3}\re^{\ri (x\cdot\xi+t\omega(|\xi|))}a(k,\theta)\,\rd\xi\right|
 \le C\|a\|_{\mathrm{pol},6}\ang t^{-3/2}
                    {\frac{\ang t^3}{\ang{t,x}^3}}.
\end{equation}
The constants depend on $c_0$, the fixed frequency neighborhood, and bounds for finitely many derivatives of $\omega$ on that neighborhood.
\end{lemma}
\begin{proof}
Put $\phi(\xi)=\omega(|\xi|)$ and $r=|x|$.
Write $I(t,x)$ for the Cartesian integral in
\eqref{eq:radial-dispersive-section4}.

\medskip\noindent
\emph{Phase geometry and spatial regions.}
For $k=|\xi|>0$ and $\theta=\xi/k$, differentiation gives
\[
 D^2\phi(\xi)=\omega''(k)\theta\otimes\theta
       +\frac{\omega'(k)}k(I-\theta\otimes\theta).
\]
Thus the radial eigenvalue is $\omega''(k)$ and each of the two
tangential eigenvalues is $\omega'(k)/k$.  Even smoothness at zero
gives $\omega'(0)=0$ and
\[
 \lim_{k\downarrow0}\frac{\omega'(k)}k=\omega''(0),
 \qquad D^2\phi(0)=\omega''(0)I.
\]
The assumed common sign is constant on the connected frequency
interval.  If it is negative, conjugating the integral replaces
$(x,\phi,A)$ by $(-x,-\phi,\overline A)$, with the same norms;
the polar case is identical.  We may therefore assume
$D^2\phi\ge c_0 I$.
Choose fixed radii $k_1<k_2<k_3$ within the neighborhood in the
hypothesis, and fix
\[
 C_*>2\left(\sup_{|\xi|\le k_3}|\nabla\phi(\xi)|+1\right).
\]
We distinguish the propagation region $r\le C_*\ang t$ from
the exterior region $r\ge C_*\ang t$, for both amplitude classes.
In the propagation region, $0\le t\le1$ also bounds $r$, so
absolute integration proves both estimates there.  For $t\ge1$,
set $v=x/t$ and $\Phi_v(\xi)=v\cdot\xi+\phi(\xi)$.
Then $|v|\le\sqrt2 C_*$ in this region, and the phase derivatives
used below are uniformly bounded on $\overline B_{k_3}$.
For $\xi,\eta\in B_{k_3}$ the segment formula gives
\[
 (\nabla\Phi_v(\xi)-\nabla\Phi_v(\eta))\cdot(\xi-\eta)
 =\int_0^1 (\xi-\eta)^TD^2\phi(\eta+s(\xi-\eta))
                (\xi-\eta)\,\rd s
 \ge c_0|\xi-\eta|^2.
\]
Consequently there is at most one critical point $\xi_v$ in this
ball, and, whenever $\xi_v\in\overline B_{k_2}$,
\[
 |\nabla\Phi_v(\xi)|\ge c_0|\xi-\xi_v|,
 \qquad \xi\in B_{k_3}.
\]
This also treats a critical point lying just outside the amplitude
support.  If there is no critical point in $\overline B_{k_2}$,
radiality and $\omega''\ge c_0$ imply
$|v|>\omega'(k_2)$.  On $k\le k_1$ we then have
\[
 |\nabla\Phi_v(\xi)|\ge |v|-\omega'(k_1)
 \ge\omega'(k_2)-\omega'(k_1)>0.
\]
The last lower bound is fixed independently of $v$.

\medskip\noindent
\emph{Cartesian amplitudes.}
We first prove \eqref{eq:radial-dispersive-section4}, including its
spatial weights.  Suppose that $\xi_v\in\overline B_{k_2}$.
For $h$ in a fixed small ball, Taylor's formula gives
\[
 \Phi_v(\xi_v+h)-\Phi_v(\xi_v)
 =\frac12 h^TH_v(h)h,\qquad
 H_v(h)=2\int_0^1(1-s)D^2\phi(\xi_v+sh)\,\rd s.
\]
Here $H_v\ge c_0I$ and $\|H_v\|_{C^3}$ is uniformly bounded by
$\|\phi\|_{C^5(B_{k_3})}$.  The matrix square root is smooth on the
fixed compact set of positive definite matrices containing these
$H_v(h)$.  Thus
\[
 F_v(h)=H_v(h)^{1/2}h,\qquad
 DF_v(0)=D^2\phi(\xi_v)^{1/2}
\]
has a uniform $C^3$ bound.  Moreover,
\[
 |DF_v(0)h|^2=h^TD^2\phi(\xi_v)h\ge c_0|h|^2,
\]
so the smallest singular value of $DF_v(0)$ is at least
$\sqrt{c_0}$.  Choose a radius $\delta>0$,
independent of $v$, so that $\xi_v+B_\delta\subset B_{k_3}$ and
\[
 \sup_{|h|\le\delta}\|DF_v(h)-DF_v(0)\|
 \le\tfrac12\sqrt{c_0}.
\]
The quantitative inverse function theorem gives a local inverse
$G_v=F_v^{-1}$ with uniform $C^3$ bounds.  In particular,
$J_v(y)=|\det DG_v(y)|$ has uniform $C^2$ bounds, and
\[
 \Phi_v(\xi_v+G_v(y))=\Phi_v(\xi_v)+\tfrac12|y|^2.
\]
These bounds are uniform in the moving critical point because
$D^m\Phi_v=D^m\phi$ for $m\ge2$.
Choose a fixed smooth cutoff $\chi_{\rm st}(h)$ supported in $B_\delta$
and equal to one on a smaller ball.  The transformed amplitude
\[
 b_v(y)=A(\xi_v+G_v(y))\chi_{\rm st}(G_v(y))J_v(y),
\]
extended by zero, has support in a fixed ball and satisfies
$\sum_{|\alpha|\le2}\|\partial^\alpha b_v\|_2\le C\|A\|_{C^2}$.

For such an amplitude $b$, Fourier inversion and the Gaussian
integral give
\[
 \left|\int_{\R^3}\re^{\ri t|y|^2/2}b(y)\,\rd y\right|
 \le Ct^{-3/2}\|\widehat b\|_1
 \le Ct^{-3/2}\sum_{|\alpha|\le2}\|\partial^\alpha b\|_2.
\]
Indeed, for $\eta>0$ the regularized Gaussian satisfies
\[
 \int_{\R^3}\re^{(\ri t/2-\eta)|y|^2+\ri y\cdot\zeta}\,\rd y
 =\left(\frac{\pi}{\eta-\ri t/2}\right)^{3/2}
   \exp\left(-\frac{|\zeta|^2}{4(\eta-\ri t/2)}\right).
\]
Its modulus is at most $Ct^{-3/2}$, uniformly in $\eta$ and
$\zeta$, since $|\eta-\ri t/2|\ge t/2$ and
$\Re(\eta-\ri t/2)^{-1}>0$.
Fourier inversion for $b$ and Fubini's theorem therefore give
the first bound for the regularized integral; dominated convergence
as $\eta\downarrow0$ gives the displayed integral.
The second inequality
follows by Cauchy--Schwarz from
$\ang\zeta^{-2}\in L^2(\R^3)$ and Plancherel's identity.
The stationary contribution is therefore bounded by
$C\|A\|_{C^2}t^{-3/2}$.  This argument uses exactly the finite
chain $\phi\in C^5$, $F_v,G_v\in C^3$, $J_v\in C^2$, and
$b_v\in H^2$.

On the complement of the critical ball, the gradient has a
uniform positive lower bound.  Two integrations by parts with
$\nabla\Phi_v\cdot\nabla_\xi/(\ri t|\nabla\Phi_v|^2)$ give
$C\|A\|_{C^2}t^{-2}$.  The same argument treats the whole
integral when no critical point lies in $\overline B_{k_2}$.
This proves the Cartesian estimate in the propagation region.

In the exterior region $r\ge C_*\ang t$, our choice of $C_*$
gives $|x+t\nabla\phi|\ge r/2$ on the amplitude support.  Set
\[
 e=x/r,\quad \lambda=t/r,\qquad
 b_{\rm ex}(\xi)=\frac1{\ri r}
   \frac{e+\lambda\nabla\phi(\xi)}
        {|e+\lambda\nabla\phi(\xi)|^2},
 \qquad L_{\rm ex}=b_{\rm ex}\cdot\nabla_\xi.
\]
Each required derivative of $b_{\rm ex}$ is $O(r^{-1})$,
uniformly in $e$ and $\lambda$.
Since $L_{\rm ex}\re^{\ri(x\cdot\xi+t\phi)}
=\re^{\ri(x\cdot\xi+t\phi)}$, $N+2$ integrations by parts yield
\[
 |I(t,x)|\le C_N\|A\|_{C^{N+2}}r^{-N-2}
 \le C_N\|A\|_{C^{N+2}}\ang t^{N-3/2}r^{-N}.
\]
Compact support removes all boundary terms.  This argument also
covers $0\le t\le1$ in the exterior region.  Since
$\ang{t,x}\asymp\ang t$ when $r\le C_*\ang t$ and
$\ang{t,x}\asymp r$ otherwise, the two estimates prove
\eqref{eq:radial-dispersive-section4} for all $t\ge0$.

\medskip\noindent
\emph{Polar amplitudes.}
Under the polar assumptions, the amplitude need not extend continuously to $\xi=0$.
We match the localization scales to the inverse powers of $k$
created by Cartesian differentiation.  Normalize $\|a\|_{\mathrm{pol},6}=1$.
On $k>0$, the polar-coordinate differentiation identities imply
\[
 |\partial_\xi^\alpha a(k,\theta)|\le C_\alpha k^{-|\alpha|},
 \qquad |\alpha|\le6.
\]
Here one uses
$\partial_{\xi_j}=\theta_j\partial_k+k^{-1}\mathcal T_j$,
where $\mathcal T_j$ is a smooth tangential vector field on the
sphere, together with
$k^j\partial_k^j=\prod_{\ell=0}^{j-1}(k\partial_k-\ell)$.
Only the joint polar orders in the stated seminorm occur.

We first work in the propagation region, where it remains to treat
$t\ge1$.  Put $\epsilon=t^{-1/2}$.  An origin ball of radius a
fixed multiple of $\epsilon$ contributes
$O(\epsilon^3)=O(t^{-3/2})$ by absolute integration.
Suppose a critical point lies in $\overline B_{k_2}$ and put
$k_v=|\xi_v|\le k_2$.
Fix $0<c<1/4$ with $ck_2<\delta$, small enough for the local
Morse construction for the uniformly nondegenerate rescaled phases
below.  This choice depends only on the fixed phase bounds and
$c_0$, and is independent of $t,v,k_v$.

If $k_v\le4\epsilon$, choose the origin cutoff equal to one for
$k\le8\epsilon$ and zero for $k\ge16\epsilon$.
It covers the critical point, and on the remaining support
$k/2\le|\xi-\xi_v|\le3k/2$.
If $k_v>4\epsilon$, choose the origin cutoff equal to one for
$k\le\epsilon$ and zero for $k\ge2\epsilon$.  Use a second cutoff
equal to one for $|\xi-\xi_v|\le c k_v/2$ and supported in
$|\xi-\xi_v|<c k_v$.  This ball lies in $B_\delta(\xi_v)$,
and on it $k\asymp k_v$ and $k>2\epsilon$.
Thus the two cutoff neighborhoods are disjoint, including when
$\xi_v$ lies just outside the amplitude support.

For the critical contribution in the second case, set
$\xi=\xi_v+k_vy$ and write
\[
 \begin{aligned}
 I_{\rm crit}(t,x)
 &=k_v^3\re^{\ri t\Phi_v(\xi_v)}
   \int_{\R^3}\re^{\ri(tk_v^2)\Psi_v(y)}A_v(y)\,\rd y,\\
 \Psi_v(y)
 &:=k_v^{-2}\bigl[\Phi_v(\xi_v+k_vy)-\Phi_v(\xi_v)\bigr].
 \end{aligned}
\]
Here $A_v$ includes the critical cutoff, which has fixed support
in $y$.  On this support the chain rule gives
\[
 |D_y^\alpha[a(|\xi_v+k_vy|,
                 (\xi_v+k_vy)/|\xi_v+k_vy|)]|
 \le C k_v^{|\alpha|}|\xi_v+k_vy|^{-|\alpha|}\le C,
 \qquad |\alpha|\le2,
\]
so $\|A_v\|_{C^2}\le C$ uniformly.
Also $\Psi_v(0)=0$, $\nabla\Psi_v(0)=0$, and
\[
 D_y^2\Psi_v(y)=D^2\phi(\xi_v+k_vy)\ge c_0 I,\qquad
 D_y^j\Psi_v(y)=k_v^{j-2}D^j\phi(\xi_v+k_vy),\quad 2\le j\le5.
\]
Since $k_v\le k_2$, the $C^5$ bounds for $\Psi_v$ are uniform.
The local Morse--Gaussian argument in the Cartesian proof applies
to this phase family, even though the rescaled phases need not be
radial.  With the large parameter $tk_v^2>16$, it gives
\[
 |I_{\rm crit}(t,x)|
 \le Ck_v^3(tk_v^2)^{-3/2}\|A_v\|_{C^2}=Ct^{-3/2}.
\]

It remains to estimate the part outside the origin and critical
balls.  Put $d=|\xi-\xi_v|$ and
\[
 b_v(\xi)=\frac{\nabla\Phi_v(\xi)}{\ri t|\nabla\Phi_v(\xi)|^2},
 \qquad L=b_v\cdot\nabla_\xi.
\]
This operator satisfies
\[
 L\re^{\ri t\Phi_v(\xi)}=\re^{\ri t\Phi_v(\xi)}.
\]
The gradient lower bound and bounded phase derivatives give
$|\partial_\xi^j b_v|\le Ct^{-1}d^{-1-j}$ for $j\le2$ on the
fixed frequency region.  If $f$ denotes the amplitude with both
cutoffs included, then
\[
 |(L^*)^2 f|\le Ct^{-2}
 \bigl(d^{-2}|\nabla^2 f|+d^{-3}|\nabla f|+d^{-4}|f|\bigr).
\]
Indeed, the formal integration-by-parts adjoint is
$L^*f=-b_v\cdot\nabla f-(\operatorname{div}b_v)f$, and
\[
\begin{aligned}
 (L^*)^2f={}&\sum_{i,j}(b_v)_i(b_v)_j\partial_{ij}f
 +\bigl((b_v\cdot\nabla)b_v
              +2(\operatorname{div}b_v)b_v\bigr)\cdot\nabla f
 +\bigl(b_v\cdot\nabla(\operatorname{div}b_v)
              +(\operatorname{div}b_v)^2\bigr)f.
\end{aligned}
\]
The origin-cutoff derivatives of order $j$ are $O(\epsilon^{-j})$
where $k\asymp\epsilon$.  Those of the critical cutoff are
$O(k_v^{-j})$ where $k\asymp k_v$.  Thus all derivatives of $f$
through order two obey the same $k^{-j}$ bounds as the amplitude.  Hence
\[
 |(L^*)^2 f|\le Ct^{-2}
       \bigl(d^{-4}+k^{-1}d^{-3}+k^{-2}d^{-2}\bigr).
\]
If $k_v\le4\epsilon$, its integral is bounded by
$Ct^{-2}\int_{\{8\epsilon\le k\le k_1\}}k^{-2}\,\rd k
\le Ct^{-2}\epsilon^{-1}$.
For $k_v>4\epsilon$, split into
$k<k_v/2$, $k_v/2\le k\le2k_v$, and $k>2k_v$,
always excluding the cutoff balls.
On the first region $d\asymp k_v$, and the integral without
$t^{-2}$ is at most
\[
 C\int_0^{k_v/2}
 (k^2k_v^{-4}+kk_v^{-3}+k_v^{-2})\,\rd k\le Ck_v^{-1}.
\]
On the middle region $d\ge c'k_v$ and $k\asymp k_v$, so the
integrand is $O(k_v^{-4})$ on a set of volume $O(k_v^3)$.
On the last region $d\asymp k$ and the integral is bounded by
$C\int_{\{2k_v\le k\le k_1\}}k^{-2}\,\rd k\le Ck_v^{-1}$.
Empty regions contribute zero.  Their sum is $O(k_v^{-1})$, hence
$O(\epsilon^{-1})$.  Two integrations by parts give
$Ct^{-2}\epsilon^{-1}=Ct^{-3/2}$
in both cases.  The cutoffs remove the origin and the critical point;
compact support removes the outer boundary terms.

When there is no critical point in $\overline B_{k_2}$, use the
origin cutoff at scale $\epsilon$ and the uniform gradient bound.
The twice differentiated amplitude has integral bounded by
$C\int_0^{k_1}(1+k^{-1}+k^{-2})k^2\,\rd k$,
including the cutoff transition, and the remaining integral is
$O(t^{-2})$.  This completes the polar time-decay estimate in
$r\le C_*\ang t$.

For the polar exterior region $r\ge C_*\ang t$, use the same
operator $L_{\rm ex}$ and a smooth cutoff $\chi_{\rm ex}(rk)$
vanishing for $rk\le1$ and equal to one for $rk\ge2$.
The omitted ball has volume $O(r^{-3})$.  On the remaining
support the differentiation identities above give
\[
 |\partial_\xi^\alpha(a\chi_{\rm ex}(rk))|
 \le C_\alpha k^{-|\alpha|},\qquad |\alpha|\le4.
\]
Every term after four adjoint integrations contains four factors
of $b_{\rm ex}$ or its derivatives, and at most four derivatives
of the amplitude.  Since $k\le k_1$, lower inverse powers of $k$
are also bounded by a constant times $k^{-4}$, and hence
\[
 |(L_{\rm ex}^*)^4(a\chi_{\rm ex}(rk))|\le Cr^{-4}k^{-4},
 \qquad
 r^{-4}\int_{1/r}^{k_1}k^{-4}k^2\,\rd k\le Cr^{-3}.
\]
The cutoffs remove the inner and outer boundary terms; an empty
radial interval contributes zero.  Thus the exterior bound is
$Cr^{-3}$, which is stronger than $C\ang t^{3/2}r^{-3}$.
For $0\le t\le1$, this exterior argument and absolute
integration on bounded $r$ again cover the whole space.
Combining with the propagation-region estimate and restoring
$\|a\|_{\mathrm{pol},6}$ proves
\eqref{eq:polar-pole-source-section5}.
All estimates use at most the six joint polar orders in the
hypothesis; the exterior argument uses four.
\end{proof}

\subsection*{Spherical expansion with a quantitative remainder}

The next lemma gives the first two terms of the stationary-phase
expansion on the sphere, with a finite-order remainder;
cf.\ \cite[Chapter~VIII]{Stein93}.

\begin{lemma}
\label{lem:spherical-carriers-v4}
{
For $A\in C^6(\mathbb S^2)$, $\omega\in\mathbb S^2$, and $\lambda\ge1$,
\begin{align}\label{eq:spherical-carriers-v4}
 \int_{\mathbb S^2}\re^{\ri \lambda\omega\cdot\theta}A(\theta)\,\rd\theta
 ={}&\frac{2\pi}{\ri \lambda}
       \bigl(\re^{\ri \lambda}A(\omega)-\re^{-\ri \lambda}A(-\omega)\bigr)-\frac\pi{\lambda^2}
       \bigl(\re^{\ri \lambda}\Delta_{\mathbb S^2}A(\omega)
              +\re^{-\ri \lambda}\Delta_{\mathbb S^2}A(-\omega)\bigr)
       +\mathcal E_3,\\*
 |\mathcal E_3|\le& C\lambda^{-3}\|A\|_{C^6(\mathbb S^2)}.\notag
\end{align}
The constants are uniform in $\omega$.  If $A=A(z,\theta)$ depends
on an external parameter $z$, then, with $\lambda$ and $\omega$ fixed,
\[
 |\partial_z^\alpha\mathcal E_3(\lambda,\omega;z)|
 \le C\lambda^{-3}
       \|\partial_z^\alpha A(z,\cdot)\|_{C^6(\mathbb S^2)}
\]
whenever the norm on the right is finite.
}
\end{lemma}
\begin{proof}
Choose an orthonormal basis $e_1,e_2$ of $\omega^\perp$ and define
the latitude average
\[
 B_\omega(u)=\int_0^{2\pi}
 A\bigl(u\omega+\sqrt{1-u^2}(e_1\cos\varphi+e_2\sin\varphi)\bigr)
 \,\rd\varphi.
\]
The spherical integral is
$I(\lambda,\omega)=\int_{-1}^1\re^{\ri\lambda u}B_\omega(u)\,\rd u$.
Under the $C^3$ regularity verified below, two integrations by parts give
\[
\begin{split}
 I(\lambda,\omega)
 &=\left[\re^{\ri\lambda u}
       \left(\frac{B_\omega(u)}{\ri\lambda}
             -\frac{B_\omega'(u)}{(\ri\lambda)^2}\right)\right]_{-1}^1
       +\mathcal E_3,\\
 \mathcal E_3
 &=\frac1{(\ri\lambda)^2}
      \int_{-1}^1\re^{\ri\lambda u}B_\omega''(u)\,\rd u.
\end{split}
\]
It remains to compute the boundary values and control this remainder.

Near either pole use the tangent-plane chart
\[
 F_\pm(y)=A\bigl(\pm\sqrt{1-|y|^2}\,\omega+y_1e_1+y_2e_2\bigr),
 \qquad |y|<\rho_0<1,
\]
with fixed $\rho_0$.  For a function $F$ on this plane let
\[
 H_F(q)=\frac1{2\pi}\int_0^{2\pi}
 F(\sqrt q\,\nu_\varphi)\,\rd\varphi,
 \qquad \nu_\varphi=(\cos\varphi,\sin\varphi).
\]
The vector $\nu_\varphi$ is the outward unit normal to the circle.
Differentiation and the divergence theorem yield, for $q>0$,
\[
\begin{split}
 H_F'(q)
 =\frac1{4\pi\sqrt q}\int_0^{2\pi}
       \nabla F(\sqrt q\,\nu_\varphi)\cdot\nu_\varphi\,\rd\varphi
 =\frac1{4\pi q}\int_{|y|\le\sqrt q}\Delta_yF(y)\,\rd y
  =\frac14\int_0^1 H_{\Delta_yF}(sq)\,\rd s.
\end{split}
\]
The last expression is continuous at $q=0$ and gives
$H_F'(0)=\Delta_yF(0)/4$.  Iterating this identity, for example,
\[
 H_F''(q)=\frac1{16}\int_0^1\int_0^1
 s\,H_{\Delta_y^2F}(s\tau q)\,\rd\tau\,\rd s,
\]
and once more for the third derivative, gives
$\|H_F^{(j)}\|_\infty\le C_j\|F\|_{C^{2j}}$ for $0\le j\le3$.
Since $B_\omega(u)=2\pi H_{F_\pm}(1-u^2)$ near $u=\pm1$,
the chain rule controls the first three derivatives at both endpoints.
Away from them, differentiation under the latitude integral gives
the same bound.  Thus
\[
 \|B_\omega\|_{C^3([-1,1])}
 \le C\|A\|_{C^6(\mathbb S^2)},
\]
uniformly in $\omega$.

For the parametrizations defining $F_\pm$, the induced metric is
\[
 g_{ij}(y)=\delta_{ij}+\frac{y_i y_j}{1-|y|^2}.
\]
Thus $g(0)=I$ and $\partial_{y_k}g_{ij}(0)=0$; the local
Laplace--Beltrami formula gives
$\Delta_yF_\pm(0)=\Delta_{\mathbb S^2}A(\pm\omega)$.
Using $q=1-u^2$ now gives
\[
 B_\omega(\pm1)=2\pi A(\pm\omega),\qquad
 B_\omega'(1)=-\pi\Delta_{\mathbb S^2}A(\omega),\qquad
 B_\omega'(-1)=\pi\Delta_{\mathbb S^2}A(-\omega).
\]
These values give the first two terms in
\eqref{eq:spherical-carriers-v4}.  One further integration by parts gives
\[
 \mathcal E_3=\frac1{(\ri\lambda)^3}
 \left([\re^{\ri\lambda u}B_\omega''(u)]_{-1}^1
       -\int_{-1}^1\re^{\ri\lambda u}B_\omega'''(u)\,\rd u\right),
\]
and hence the stated $C\lambda^{-3}\|A\|_{C^6}$ bound.
Applying the same linear calculation to $\partial_z^\alpha A$, with
$\lambda$ and $\omega$ fixed, proves the parameter estimate.
\end{proof}

\subsection*{Endpoint transforms and radial integrals}

We use the classical endpoint method for Fourier integrals with an
algebraic endpoint; see \cite{Olver97}.  The lemma also records the
mixed-parameter bounds used below.

\begin{lemma}
\label{lem:finite-endpoint-v4}
Fix $\epsilon\in\{-1,1\}$, $\beta>-1$, and an integer
$L>\beta+1$.  Let $f=f(k,\theta,u)$ be supported in $0\le u\le u_0$
and smoothly cut off at $u_0$, with
\[
 |\partial_u^j f(k,\theta,u)|\le C_j u^{\beta-j},
 \qquad 0\le j\le L.
\]
Define
\[
 \mathcal E_f(k,\theta,\zeta)
 :=\int_0^{u_0}\re^{\ri\epsilon\zeta u}f(k,\theta,u)\,\rd u,
 \qquad
 T_f(k,\theta,\zeta):=\zeta^{\beta+1}\mathcal E_f(k,\theta,\zeta),
 \qquad \zeta>0.
\]
Then the following estimates hold.
\begin{enumerate}[label=(\roman*)]
\item For every $\zeta>0$,
\[
 |\mathcal E_f(k,\theta,\zeta)|
 \le C\min\{1,\zeta^{-\beta-1}\},\qquad
 |T_f(k,\theta,\zeta)|\le C\min\{1,\zeta^{\beta+1}\}.
\]
\item Put $D_k=k\partial_k$ and $D_u=u\partial_u$.
Let $m\ge0$ be an integer and $\alpha$ an angular multi-index.
Suppose additionally that
\[
 |D_k^{j_1}\nabla_\theta^\alpha D_u^{j_2} f(k,\theta,u)|
 \le C_{j_1,j_2,\alpha}u^\beta
\]
for
\begin{equation}\label{eq:endpoint-mixed-allocation-B}
 0\le j_1\le m,\qquad 0\le j_2\le L+m-j_1.
\end{equation}
Then
\begin{equation}\label{eq:endpoint-small-parameter-v3}
 \bigl|(D_k+\zeta\partial_\zeta)^m
       \nabla_\theta^\alpha T_f(k,\theta,\zeta)\bigr|
 \le C_{m,\alpha}\min\{1,\zeta^{\beta+1}\},\qquad \zeta>0.
\end{equation}
In particular, for fixed $t>0$,
\begin{equation}\label{eq:endpoint-ray-symbol-v196}
 \bigl|D_k^m\nabla_\theta^\alpha[T_f(k,\theta,tk)]\bigr|
 \le C_{m,\alpha}\min\{1,(tk)^{\beta+1}\}.
\end{equation}
\end{enumerate}
The constants depend on the indicated derivative bounds.
A common multiplicative parameter factor in these bounds is
retained in the conclusions.
\end{lemma}
\begin{proof}
\emph{Basic transform bounds.}
We prove (i) for any amplitude $g$ satisfying the stated endpoint
bounds; the same argument will then apply to differentiated
amplitudes.  If $0<\zeta\le1$, its defining integral gives
\[
 |\mathcal E_g(k,\theta,\zeta)|
 \le C\int_0^{u_0}u^\beta\,\rd u\le C,
\]
since $\beta>-1$.  For $\zeta\ge1$, take
$\chi_*\in C_c^\infty([0,\infty))$ equal to one on $[0,1]$
and zero on $[2,\infty)$, and split using $\chi_*(\zeta u)$.
The lower part is bounded by
$C\int_0^{2/\zeta}u^\beta\,\rd u\le C\zeta^{-\beta-1}$.
On the complement, $L$ integrations by parts give
\[
 C\zeta^{-L}\int_{\{\zeta^{-1}\le u\le u_0\}}
       u^{\beta-L}\,\rd u
 \le C\zeta^{-\beta-1},
\]
because $L>\beta+1$.  If $h\ge1$ derivatives fall on
$\chi_*(\zeta u)$, they have size $O(\zeta^h)$ on
$u\asymp\zeta^{-1}$, and their contribution is bounded by
\[
 C\zeta^{-L}\zeta^h\int_{u\asymp\zeta^{-1}}
          u^{\beta-L+h}\,\rd u
 \le C\zeta^{-\beta-1}.
\]
All integrals are restricted to the support of $g$, with empty
regions contributing zero.  The cutoffs remove both boundary
terms.  Multiplication by $\zeta^{\beta+1}$ proves the bound for
$T_g$, and hence (i).

\medskip\noindent
\emph{Parameter derivatives and the ray estimate.}
Integration by parts in $u$ gives
\[
 \zeta\partial_\zeta\mathcal E_f=-\mathcal E_{(1+D_u)f},
 \qquad
 \zeta\partial_\zeta T_f=T_{(\beta-D_u)f}.
\]
The boundary term vanishes because $uf=O(u^{\beta+1})\to0$
at zero and $f$ is cut off at $u_0$.
Commuting the parameter derivatives and iterating yields
\begin{equation}\label{eq:endpoint-transfer-v3}
 (D_k+\zeta\partial_\zeta)^m\nabla_\theta^\alpha T_f
 =\sum_{j_1=0}^m\binom{m}{j_1}\,
 T_{D_k^{j_1}\nabla_\theta^\alpha(\beta-D_u)^{m-j_1}f}.
\end{equation}
Every transform on the right uses the same exponent $\beta$.
For its amplitude, $L$ further ordinary endpoint derivatives use
at most $L+m-j_1$ operations of $D_u$, since
\[
 \partial_u^j=u^{-j}\prod_{h=0}^{j-1}(D_u-h),\qquad
 \left|\partial_u^j
   D_k^{j_1}\nabla_\theta^\alpha(\beta-D_u)^{m-j_1}f\right|
 \le C_{m,\alpha,j}u^{\beta-j},\quad 0\le j\le L.
\]
Thus \eqref{eq:endpoint-mixed-allocation-B} supplies precisely
these bounds, using at most $L+m+|\alpha|$ simultaneous radial,
endpoint, and angular operations.  Applying (i) to each amplitude
proves \eqref{eq:endpoint-small-parameter-v3}; in the small-parameter
range this uses the absolute bound on each defining
$\mathcal E$-integral and its prefactor $\zeta^{\beta+1}$.
Finally, the chain rule at fixed $t$ gives
\[
 D_k^m[T_f(k,\theta,tk)]
 =\left.(D_k+\zeta\partial_\zeta)^mT_f(k,\theta,\zeta)
   \right|_{\zeta=tk}.
\]
Angular derivatives commute with these operations, proving
\eqref{eq:endpoint-ray-symbol-v196}.
\end{proof}

The next lemma controls truncated radial oscillatory integrals
arising after the endpoint Fourier transform.

\begin{lemma}
\label{lem:truncated-radial-v4}
{
Let $R\ge1$, $\ang\ell\lesssim R$, and let $\psi(Rk)$ be a smooth
cutoff equal to zero for $Rk\le c$ and one for $Rk\ge2c$.
Suppose $|(k\partial_k)^j q(k)|\le C_j$ for $0<k\le1$ and $j\le J$.
For $k\ge1$ assume either support in a fixed compact interval, with
derivatives through order $J$ uniformly bounded there, or
$|\partial_k^j q(k)|\le C_j k^{-B-j}$ with $B>\max(\alpha,0)$.
If $J$ is an integer satisfying $J>\max\{\alpha,0\}$, then
\begin{equation}\label{eq:truncated-radial-v4}
 \left|\int_0^\infty \re^{-\ri k\ell}\psi(Rk)k^{\alpha-1}q(k)\,\rd k\right|
 \lesssim
 \begin{cases}
 \ang\ell^{-\alpha}+R^{-\alpha},&\alpha>0,\\
 \log(2+R/\ang\ell),&\alpha=0,\\
 R^{-\alpha},&\alpha<0.
 \end{cases}
\end{equation}
The constants depend on the fixed cutoff, $\alpha,J,B$, and the
constants in the assumed derivative bounds, but not on $R$ or $\ell$.
}
\end{lemma}
\begin{proof}
Put $a(k)=k^{\alpha-1}q(k)$.  The identities expressing ordinary
$k$-derivatives through $k\partial_k$ give
\[
 |\partial_k^j a(k)|\le C_j k^{\alpha-1-j},
 \qquad 0<k\le1,\quad 0\le j\le J.
\]
At high frequency the compact support or decay assumption yields
$\sum_{j=0}^J\|\partial_k^j a\|_{L^1(1,\infty)}\le C$.
We integrate by parts on the full positive axis below and only
then split the derivative integrals into low and high frequencies;
there is no artificial boundary at $k=1$.  The smooth cutoffs
remove the boundary at zero, and the assumptions remove it at
infinity.

If $|\ell|\le1$, then $\ang\ell\asymp1$.  Absolute integration gives
\[
 \left|\int_0^\infty \re^{-\ri k\ell}\psi(Rk)a(k)\,\rd k\right|
 \le C\left(1+\int_{\{c/R\le k\le1\}}k^{\alpha-1}\,\rd k\right).
\]
The resulting bounds are $C$, $C\log(2+R)$, and $CR^{-\alpha}$
for $\alpha>0$, $\alpha=0$, and $\alpha<0$, respectively.
All possibly empty intervals in the proof contribute zero.

Now let $|\ell|\ge1$, so $|\ell|\asymp \ang\ell$, and fix
$\chi_*\in C_c^\infty([0,\infty))$ equal to one on $[0,1]$
and zero on $[2,\infty)$.
For $\alpha>0$, removing $\psi(Rk)$ costs $CR^{-\alpha}$ by
absolute integration on $k\le2c/R$.  Split the resulting integral
using $\chi_*(\ang\ell k)$.  The lower part is $O(\ang\ell^{-\alpha})$.
After $J$ integrations by parts on the upper part, a term with
$h\ge1$ derivatives on $\chi_*(\ang\ell k)$ is bounded by
\[
 C \ang\ell^{-J}\int_{k\asymp \ang\ell^{-1}}
       \ang\ell^h k^{\alpha-1-J+h}\,\rd k
 \le C \ang\ell^{-\alpha}.
\]
If no derivative falls on this cutoff, the bound is
\[
 C \ang\ell^{-J}\int_{\{\ang\ell^{-1}\le k\le1\}}
       k^{\alpha-1-J}\,\rd k+C \ang\ell^{-J}
 \le C \ang\ell^{-\alpha},
\]
since $J>\alpha$.  Transitions meeting $k\ge1$ are bounded
uniformly when $\ang\ell$ is bounded.  This proves the two claimed terms.
In particular, $\alpha=1$ uses two integrations; one would leave
$\ang\ell^{-1}\log \ang\ell$ under the stated derivative bounds.

For $\alpha=0$, retain $\psi(Rk)$ and split with the same
$\chi_*(\ang\ell k)$.  The lower part is bounded by
$C\log(2+R/\ang\ell)$ through absolute integration on
$c/R\le k\le2/\ang\ell$.  On the upper part, one integration by parts
has three types of terms.  A derivative of $\chi_*(\ang\ell k)$ gives
\[
 C \ang\ell^{-1}\int_{k\asymp \ang\ell^{-1}}\ang\ell k^{-1}\,\rd k\le C.
\]
A derivative of $\psi(Rk)$ is supported on $k\asymp R^{-1}$;
if that support meets $1-\chi_*(\ang\ell k)$, then $R\lesssim \ang\ell$.
Together with $\ang\ell\lesssim R$, this gives $R\asymp \ang\ell$ and
\[
 C \ang\ell^{-1}\int_{k\asymp R^{-1}}Rk^{-1}\,\rd k
 \le C R/\ang\ell\le C.
\]
Finally, a derivative of $a(k)=k^{-1}q(k)$ contributes at most
\[
 C \ang\ell^{-1}\int_{\{\ang\ell^{-1}\le k\le1\}}k^{-2}\,\rd k
       +C \ang\ell^{-1}\le C.
\]
The high-frequency derivative integrals are included in the last
constant, and bounded transition regions obey the same estimates.
Hence the upper part is $O(1)$, absorbed by $\log(2+R/\ang\ell)$.

For $\alpha<0$, absolute integration alone gives
$C\int_{\{c/R\le k\le1\}}k^{\alpha-1}\,\rd k+C
\le CR^{-\alpha}$, completing the proof.
\end{proof}

\section{Transverse inversion and bounded-variation tools}
\label{app:terminal-rates}

We record the localized transverse inversion used to pass from the
joint spectral measure to the complex dispersion roots.  The two
remaining tools control absolute damping amplitudes and
one-dimensional phases with a group-velocity sublevel bound.

The first lemma localizes the classical Cauchy boundary-value formula;
see \cite[Part~I, Chapter~2]{Mus58}.

\begin{lemma}
\label{lem:transverse-cauchy-v18}
Let $\zeta\in C_c^\infty((0,\infty))$ be real.  Suppose
$q=a+\ri b$, where $a$ belongs to a fixed compact subset of
$(0,\infty)$ and $0\le b\le b_0$.
Assume that $\zeta=1$ on a fixed neighborhood of all the possible
values of $a$.
For any fixed nonnegative integer $m$, define
\begin{equation}\label{eq:transverse-transform-v18}
 \mathcal T_m(L,q)=\frac1{2\pi}\int_\R \re^{\ri Lk}\zeta(k)k^m
  \left(\frac1{k^2-q^2}-\frac1{k^2-\bar q^{\,2}}\right)\,\rd k.
\end{equation}
When $b=0$, the fractions $(k^2-q^2)^{-1}$ and
$(k^2-\bar q^{\,2})^{-1}$ mean
$(k^2-a^2-\ri0)^{-1}$ and $(k^2-a^2+\ri0)^{-1}$,
respectively.  For each fixed $L$, the transform extends continuously
from $b>0$ to this boundary, with
\begin{equation}\label{eq:transverse-real-boundary-v199}
 \mathcal T_m(L,a)=\frac{\ri}{2}a^{m-1}\re^{\ri La},
 \qquad L\in\R.
\end{equation}
Uniformly for these $q$, and for every fixed integer $N\ge1$,
\begin{equation}\label{eq:transverse-residues-v18}
 \mathcal T_m(L,q)=
 \begin{cases}
  \dfrac{\ri}{2}q^{m-1}\re^{\ri Lq}+O_{m,N}(|L|^{-N}),&L\ge1,\\[1.5mm]
  \dfrac{\ri}{2}\bar q^{\,m-1}\re^{\ri L\bar q}+O_{m,N}(|L|^{-N}),&L\le-1.
 \end{cases}
\end{equation}
Moreover,
\begin{equation}\label{eq:transverse-small-parameter-v18}
 \sup_{L\in\R}|\mathcal T_m(L,q)|\le C_m.
\end{equation}
\end{lemma}

\begin{proof}
On the positive support of $\zeta$, division at the positive root gives
\[
 \frac{k^m}{k^2-q^2}
 =\frac{q^{m-1}}{2(k-q)}+H_m(k,q).
\]
The apparent singularity of $H_m$ at $k=q$ is removable, and the
other denominator $k+q$ stays away from zero.  For $m=0$, explicitly,
$H_0(k,q)=-[2q(k+q)]^{-1}$, with $q$ bounded away from zero.
Thus $\zeta H_m$ has uniformly bounded $L^1$ norms of every fixed
$k$-derivative, with constants allowed to depend on $m$.
Its Fourier transform is $O_{m,N}(|L|^{-N})$;
the same holds with $q$ replaced by $\bar q$.

For $L\ne0$, interpret the following decomposition and its Fourier
integrals in the oscillatory sense:
\[
 \frac{\zeta(k)}{k-q}
 =\frac1{k-q}+\frac{\zeta(k)-1}{k-q}.
\]
The support of the second term stays a fixed positive distance from $a$. This term tends to zero at infinity, and all its derivatives of positive order are uniformly integrable. Repeated integration by parts therefore
gives $O_N(|L|^{-N})$, uniformly as $b\downarrow0$.
Closing the full Cauchy integral in the appropriate half-plane gives
\[
 \frac1{2\pi}\int_\R\frac{\re^{\ri Lk}}{k-q}\,\rd k
       =\ri \re^{\ri Lq}\one_{\{L>0\}},\qquad
 \frac1{2\pi}\int_\R\frac{\re^{\ri Lk}}{k-\bar q}\,\rd k
       =-\ri \re^{\ri L\bar q}\one_{\{L<0\}}.
\]
The boundary values at $b=0$ are understood in the same
distributional sense.  Subtraction proves
\eqref{eq:transverse-residues-v18}, including the positive
$\ri/2$ coefficient for both signs of $L$.

To prove the uniform bound, let $b>0$.  The difference of the
two kernels is
\[
 \frac1{k^2-q^2}-\frac1{k^2-\bar q^{\,2}}
 =\frac{4\ri ab}{(k^2-a^2+b^2)^2+(2ab)^2}.
\]
On the fixed positive support of $\zeta$, the substitution
$X=k^2-a^2+b^2$ has Jacobian $\rd X=2k\,\rd k$, uniformly
bounded away from zero.  Taking absolute values gives
\[
 |\mathcal T_m(L,q)|
 \le C_m\int_\R\frac{2ab}{X^2+(2ab)^2}\,\rd X
 \le C_m,\qquad L\in\R.
\]
For $b=0$, the Sokhotski--Plemelj formula gives
\[
 \frac1{X\mp\ri0}=\pv\frac1X\pm\ri\pi\delta_0(X),
 \qquad
 \frac1{k^2-a^2-\ri0}-\frac1{k^2-a^2+\ri0}
 =2\pi\ri\delta(k^2-a^2).
\]
Since $\zeta$ is supported on the positive axis and $\zeta(a)=1$,
\[
 \mathcal T_m(L,a)
 =\ri\int_\R \re^{\ri Lk}\zeta(k)k^m\delta(k^2-a^2)\,\rd k
 =\frac{\ri}{2}a^{m-1}\re^{\ri La}.
\]
This proves \eqref{eq:transverse-real-boundary-v199}.
The same substitution $X=k^2-a^2+b^2$ identifies the kernels for
$b>0$ as a Poisson approximate identity applied to the smooth
compactly supported factor $\re^{\ri Lk}\zeta(k)k^{m-1}$.
For fixed $L$, it therefore also proves continuity as $q$ approaches
the real boundary.  The bound remains uniform for all $L$, including
$m=0$ because $a$ stays away from zero, and proves
\eqref{eq:transverse-small-parameter-v18}.
\end{proof}

The following lemma controls the variation of the damping factors
on the continuum side of the outer spectral edge.

\begin{lemma}
\label{lem:terminal-damping-bv-v2}
Assume $\gamma_1\ge1$ and fix $\varepsilon_{\rm tr}>0$
sufficiently small.
Put $\nu_{\rm c}(h)=\nu(1+h)$ for $0<h<\varepsilon_{\rm tr}$.
For all $t\ge0$ and $0\le v\le2$, the functions
\[
 \re^{-t|1-v(1+h)|\nu_{\rm c}(h)},\qquad \re^{-t(1+v(1+h))\nu_{\rm c}(h)}
\]
have uniformly bounded supremum norm and total variation.
\end{lemma}

\begin{proof}
Write $\nu_{\rm c}(h)=h^{\gamma_1}b(h)$.  The exact factor
\eqref{eq:damping-exact-factor-section2} gives
\[
 b(h)=\frac{\pi^2\mu_{\rm out}((1+h)^{-1})}
            {(1+h)^{\gamma_1+3}\omega(1+h)}.
\]
The one-sided endpoint expansions in
Proposition~\ref{prop:terminal-jet-section2}(ii)--(iii) give
\[
 \omega(1+h)\longrightarrow\omega_*>0,\qquad
 h\omega'(1+h)\longrightarrow0,\qquad
 h^2\omega''(1+h)\longrightarrow0.
\]
The logarithmic derivative estimates when $\gamma_1=1$ satisfy
these same weighted limits.  Since the numerator and the factor
$(1+h)^{-\gamma_1-3}$ are smooth at zero, the quotient rule yields
\[
 b(h)\longrightarrow b_0:=\frac{\pi^2\mu_{\rm out}(1)}{\omega_*}>0,
 \qquad hb'(h)\longrightarrow0,\qquad
 h^2b''(h)\longrightarrow0.
\]
The product rule now gives
\[
 \begin{split}
 \nu_{\rm c}'&=h^{\gamma_1-1}(\gamma_1b+hb'),\\
 \nu_{\rm c}''&=h^{\gamma_1-2}
       [\gamma_1(\gamma_1-1)b+2\gamma_1hb'+h^2b''].
 \end{split}
\]
Consequently $\nu_{\rm c}'>0$ near zero and
\[
 \frac{\nu_{\rm c}\nu_{\rm c}''}{(\nu_{\rm c}')^2}
 =\frac{b[\gamma_1(\gamma_1-1)b+2\gamma_1hb'+h^2b'']}
        {(\gamma_1b+hb')^2}
 \longrightarrow\frac{\gamma_1-1}{\gamma_1},
\]
including the limit zero when $\gamma_1=1$.
Hence, on one fixed interval,
\[
 Z(h)=h+\frac{\nu_{\rm c}(h)}{\nu_{\rm c}'(h)},\qquad
 Z'(h)=2-\frac{\nu_{\rm c}(h)\nu_{\rm c}''(h)}{\nu_{\rm c}'(h)^2}
 \longrightarrow1+\frac1{\gamma_1}>0.
\]
After decreasing $\varepsilon_{\rm tr}$, $Z$ is strictly increasing.
If $v=0$, the outgoing exponential is monotone.
For $v>0$ set $h_v=(1-v)/v$.  When $h<h_v$,
\[
 \partial_h[(h_v-h)\nu_{\rm c}(h)]=\nu_{\rm c}'(h)[h_v-Z(h)],
\]
which changes sign at most once.  When $h>h_v$,
$\partial_h[(h-h_v)\nu_{\rm c}(h)]=\nu_{\rm c}+(h-h_v)\nu_{\rm c}'>0$.
If $h_v\le0$, only the latter region occurs.  Thus
$\re^{-tv|h_v-h|\nu_{\rm c}(h)}$ is monotone on at most three intervals and
takes values in $[0,1]$; its variation is uniformly bounded.
The function $(1+v(1+h))\nu_{\rm c}(h)$ has derivative
$v\nu_{\rm c}+(1+v+vh)\nu_{\rm c}'>0$.
Hence the incoming damping factor is nonincreasing in $h$.
\end{proof}

We use the standard sublevel-set argument underlying van der Corput
estimates; cf.\ \cite[Chapter~VIII]{Stein93}.  The form below allows
bounded-variation amplitudes.

\begin{lemma}
\label{lem:sublevel-vdc-appendix}
Let $I$ be either of the intervals
$(1-\varepsilon_{\rm tr},1)$ and $(1,1+\varepsilon_{\rm tr})$,
or their union.  Suppose that $g\ge c_0>0$ and $g^{-1}\in BV(I)$.
Assume that $V$ is $C^1$ and monotone on each interval of a fixed
finite partition of $I$, and that, uniformly for $v\in\R$,
\[
 |\{\rho\in I:|V(\rho)-v|\le\epsilon\}|
 \le C_0\epsilon,\qquad 0<\epsilon\le1.
\]
Let $\phi_v$ be real and locally $C^1$ on these intervals, with
$\phi_v'=g(V-v)$.  For each $t\ge1$ and $v\in\R$, suppose the
amplitude $B_{t,v}$ satisfies
\[
 \|B_{t,v}\|_{L^\infty(I)}+\operatorname{Var}_I(B_{t,v})\le C_B.
\]
Then
\begin{equation}\label{eq:group-velocity-vdc-appendix}
 \left|\int_I \re^{\ri t\phi_v(\rho)}B_{t,v}(\rho)\,\rd\rho\right|
 \le Ct^{-1/2}.
\end{equation}
Here $C$ depends only on $c_0,C_0,C_B$, the variation of $g^{-1}$,
and the number of partition intervals, and is uniform over the
parameters for which the amplitude bound holds.
On disconnected $I$, variation is the sum over its components;
no continuity across $\rho=1$ is required.
\end{lemma}

\begin{proof}
Fix $t\ge1$ and $v\in\R$, and put $\epsilon=t^{-1/2}$.
The sublevel set
$E_{\epsilon,v}=\{\rho\in I:|V(\rho)-v|\le\epsilon\}$
contributes at most $C_B|E_{\epsilon,v}|\le C\epsilon$.
Split its complement at the fixed partition points.  By monotonicity
of $V$, the resulting intervals have a uniformly bounded number.
On each such interval $J$, the change of variables $s=V(\rho)-v$ gives
\[
 \begin{aligned}
 \|(V-v)^{-1}\|_{L^\infty(J)}&\le\epsilon^{-1},\\
 \operatorname{Var}_J((V-v)^{-1})
 &=\int_J\frac{|V'(\rho)|}{|V(\rho)-v|^2}\,\rd\rho
 =\int_{V(J)-v}\frac{\rd s}{s^2}
 \le2\epsilon^{-1}.
 \end{aligned}
\]
The image interval carries positive Lebesgue measure, also when
$V$ is decreasing.  Since $(\phi_v')^{-1}=g^{-1}/(V-v)$,
the $BV$ product inequality gives
\[
 \left\|\frac{B_{t,v}}{\phi_v'}\right\|_{L^\infty(J)}
 +\operatorname{Var}_J\left(\frac{B_{t,v}}{\phi_v'}\right)
 \le C\epsilon^{-1}.
\]

On a compact subinterval $J'\Subset J$, integration by parts in
the Stieltjes sense, using one-sided traces at its endpoints, yields
\[
 \int_{J'} \re^{\ri t\phi_v}B_{t,v}\,\rd\rho
 =\left[\frac{\re^{\ri t\phi_v}B_{t,v}}{\ri t\phi_v'}\right]_{\partial J'}
 -\frac1{\ri t}\int_{J'}\re^{\ri t\phi_v}
       \,\rd\left(\frac{B_{t,v}}{\phi_v'}\right).
\]
The supremum and variation bounds control both terms by
$C(t\epsilon)^{-1}$, independently of $J'$.
Letting $J'$ increase to $J$ gives the same bound on $J$.
Summing over the fixed finite collection of intervals and adding
the sublevel contribution gives
$C(\epsilon+(t\epsilon)^{-1})=Ct^{-1/2}$, proving
\eqref{eq:group-velocity-vdc-appendix}.
\end{proof}

\section{Convolution and real-frequency source estimates}
\label{app:electric-source}

We collect weighted convolution and Poisson estimates, followed by
Fourier-ray bounds and their application to finite source histories.

\begin{lemma}
\label{lem:weighted-convolution-v4}
Write $W_m(A,x)=\ang{x/A}^{-m}$ for $A\ge1$.  For $m\ge0$,
\begin{equation}\label{eq:weighted-kernel-convolution-v4}
 |\mathcal K|*_xW_m(A,\cdot)(x)
 \le C_m\|\ang y^m\mathcal K(y)\|_1W_m(A,x).
\end{equation}
For $m>3$ and $A,B\ge1$,
\begin{equation}\label{eq:two-scale-convolution-section5}
 W_m(A,\cdot)*_xW_m(B,\cdot)(x)
 \le C_m\min(A^3,B^3)W_m(A+B,x).
\end{equation}
For $a,b>1$,
\begin{equation}\label{eq:time-convolution-v4}
 \int_0^t\ang{t-s}^{-a}\ang s^{-b}\,\rd s
 \lesssim_{a,b}\ang t^{-\min(a,b)}.
\end{equation}
Consequently, if $|\mathcal K(u,x)|\lesssim\ang u^{-a}W_m(\ang u,x)$,
$|F(s,x)|\lesssim\ang s^{-N}\ang x^{-m}$, $m>3$ and $N\ge a>1$,
then
\begin{equation}\label{eq:fixed-scale-transfer-v4}
 |\mathcal K*_{t,x}F(t,x)|\lesssim\ang t^{-a}W_m(\ang t,x).
\end{equation}
\end{lemma}
\begin{proof}
\emph{Spatial convolution.}
The inequality
$\ang{x/A}\le\sqrt2\ang{(x-y)/A}\ang{y/A}$, together with
$A\ge1$, proves \eqref{eq:weighted-kernel-convolution-v4} by
integration against $|\mathcal K(y)|$.
For \eqref{eq:two-scale-convolution-section5}, we may assume
$A\le B$ by symmetry.  If $|x|\le2(A+B)$, then
$W_m(A+B,x)$ has a fixed positive lower bound, and
\[
 (W_m(A,\cdot)*_xW_m(B,\cdot))(x)
 \le\|W_m(A,\cdot)\|_1\|W_m(B,\cdot)\|_\infty
 \lesssim A^3W_m(A+B,x).
\]
If $r=|x|>2(A+B)$, split the convolution into $|x-y|\ge r/2$
and $|x-y|<r/2$.  On the latter region $|y|>r/2$, so
\[
\begin{split}
 (W_m(A,\cdot)*_xW_m(B,\cdot))(x)
 &\lesssim(B/r)^m\int_{\R^3}W_m(A,y)\,\rd y
        +(A/r)^m\int_{\R^3}W_m(B,x-y)\,\rd y\\
 &\lesssim A^3(B/r)^m+B^3(A/r)^m
 \lesssim A^3W_m(A+B,x).
\end{split}
\]
The last step uses $m>3$, $A\le B$, and $A+B\asymp B$.

\emph{Time convolution.}
For $t\ge1$, the two halves of the time integral obey
\[
\begin{split}
 \int_0^{t/2}\ang{t-s}^{-a}\ang s^{-b}\,\rd s
 &\lesssim \ang t^{-a}\int_0^{t/2}\ang s^{-b}\,\rd s\lesssim \ang t^{-a},\\
 \int_{t/2}^t\ang{t-s}^{-a}\ang s^{-b}\,\rd s
 &\lesssim \ang t^{-b}\int_0^{t/2}\ang u^{-a}\,\rd u\lesssim \ang t^{-b}.
\end{split}
\]
Here $a,b>1$ ensure integrability.  Their sum is
bounded by $C\ang t^{-\min(a,b)}$.  For $0\le t\le1$, direct integration
gives the same estimate, proving \eqref{eq:time-convolution-v4}.

\emph{The fixed source scale.}
For \eqref{eq:fixed-scale-transfer-v4}, apply
\eqref{eq:two-scale-convolution-section5} with
$A=\ang{t-s}$ and $B=1$.  The volume factor is one, and
\[
 W_m(\ang{t-s}+1,x)\lesssim W_m(\ang t,x),
 \qquad \ang{t-s}+1\le2\ang t.
\]
The remaining time integral is bounded by $C\ang t^{-a}$ using
\eqref{eq:time-convolution-v4} with $b=N\ge a>1$.
\end{proof}

The following weighted Poisson estimate uses the Newtonian kernel;
for background on potential kernels, see \cite{Stein70}.

\begin{lemma}
\label{lem:coulomb-moments-v4}
For $j=2,3$, $|\beta|=j$, and a scalar function $H$,
\begin{equation}\label{eq:direct-coulomb-v4}
 |\mathcal P_z\partial_z^\beta H(z)|
 \le C_j\ang z^{-j-2}
        \sum_{|\alpha|\le j}\|\ang y^{j+3}\partial^\alpha H(y)\|_\infty.
\end{equation}
\end{lemma}
\begin{proof}
Write
\[
 N_j(H)=\sum_{|\alpha|\le j}
       \|\ang y^{j+3}\partial^\alpha H(y)\|_\infty,
 \qquad \mathcal K(z)=-\nabla_z(4\pi|z|)^{-1}.
\]
The bound $|\mathcal K(z)|\lesssim|z|^{-2}$ makes the kernel locally
integrable.  For $|z|\le2$, splitting into $|z-y|<1$ and its complement
therefore gives
\[
\begin{split}
 |\mathcal P_z\partial_z^\beta H(z)|
 \le C\|\partial^\beta H\|_\infty
          \int_{|w|<1}|w|^{-2}\,\rd w+C\|\partial^\beta H\|_1
 \le C_jN_j(H)\le C_j\ang z^{-j-2}N_j(H),
\end{split}
\]
where $\ang y^{-j-3}$ is integrable in $\R^3$.

For $R=|z|>2$, choose a smooth cutoff $\chi_R$ equal to one on
$|y|\le R/4$, supported in $|y|\le R/2$, and satisfying
$|\partial^\alpha\chi_R|\le C_\alpha R^{-|\alpha|}$.
Splitting the source gives
\[
 \mathcal P_z\partial_z^\beta H(z)
 =\partial^\beta\mathcal K*_z(\chi_RH)(z)
  +\mathcal K*_z\partial^\beta\bigl((1-\chi_R)H\bigr)(z).
\]
On the central source support, $|z-y|\ge R/2$, so integration by parts
transfers all $j$ derivatives to the smooth kernel.  Consequently,
\[
 |\partial^\beta\mathcal K*_z(\chi_RH)(z)|
 \le C_jR^{-j-2}\|H\|_1
 \le C_jR^{-j-2}\|\ang y^{j+3}H\|_\infty.
\]
For the exterior source the derivatives remain on the product.
The term with no derivative on $\chi_R$ is supported in $|y|\ge R/4$;
terms with $h\ge1$ derivatives on $\chi_R$ lie in
$R/4\le|y|\le R/2$ and contribute the bounded factors $R^{-h}$.
Thus the Leibniz rule gives
\[
 |\partial^\beta((1-\chi_R)H)(y)|
 \le C_j\mathbf1_{\{|y|\ge R/4\}}\ang y^{-j-3}N_j(H).
\]
Writing $z=Re$, $|e|=1$, and $y=Rw$, we obtain
\[
\begin{split}
 |\mathcal K*_z\partial^\beta((1-\chi_R)H)(z)|
 &\le C_jN_j(H)\int_{|y|\ge R/4}|z-y|^{-2}\ang y^{-j-3}\,\rd y\\
 &\le C_jN_j(H)R^{-2}R^{-j-3}R^3
       \int_{|w|\ge1/4}|e-w|^{-2}|w|^{-j-3}\,\rd w\\
 &\le C_jR^{-j-2}N_j(H).
\end{split}
\]
The last integral is uniform for $|e|=1$: its domain avoids $w=0$,
its singularity at $w=e$ is locally integrable, and its tail is
bounded by $C\int_2^\infty r^{-j-3}\,\rd r$.
Combining the two source regions proves \eqref{eq:direct-coulomb-v4}.
\end{proof}

\subsection*{Fourier-ray and finite-history estimates}

We first estimate weighted Fourier integrals along velocity rays,
then apply them to the finite source history in
\eqref{eq:source-history-v5}.

\begin{lemma}
\label{lem:velocity-ray-source-v6}
Suppose $M\ge11$.  Uniformly for
$0\le k\le k_0$, $\theta\in\mathbb S^2$,
$|\beta|+|\nu|\le10$, and integers $0\le L\le11$, we have
\begin{align}
 \sup_{u\in\R}\ang u^L
 \left|\partial_\xi^\beta\partial_\eta^\nu
             \widehat h_{\rm in}(k\theta,u\theta)\right|
 &\le C[h_{\rm in}]_{M},
 \label{eq:fourier-source-decay-appendix}\\
 \left\|u^L\partial_\xi^\beta\partial_\eta^\nu
             \widehat h_{\rm in}(k\theta,u\theta)\right\|_{L^2(\R_u)}
 &\le C[h_{\rm in}]_{M}.
 \label{eq:fourier-ray-L2-v6}
\end{align}
Consequently, for $0\le d\le10$ and every $U\ge0$,
\begin{equation}\label{eq:fourier-ray-L1-v6}
 \int_0^U u^d\left|\partial_\xi^\beta\partial_\eta^\nu
             \widehat h_{\rm in}(k\theta,u\theta)\right|\,\rd u
 \le C\min\{U,1\}[h_{\rm in}]_{M}.
\end{equation}
\end{lemma}
\begin{proof}
The Fourier representation is
\[
 \partial_\xi^\beta\partial_\eta^\nu
       \widehat h_{\rm in}(k\theta,u\theta)
 =\iint_{\R^6}\re^{-\ri(k\theta\cdot x+u\theta\cdot v)}
       (-\ri x)^\beta(-\ri v)^\nu h_{\rm in}(x,v)\,\rd x\,\rd v.
\]
Here $|\beta|+|\nu|\le10$ and $M+6\ge17$.
For $|u|\le1$, taking absolute values and using \eqref{eq:data-norm}
 bounds this integral by
\[
 C[h_{\rm in}]_{M}
 \iint_{\R^6}\ang{x,v}^{|\beta|+|\nu|-M-6}\,\rd x\,\rd v
 \le C[h_{\rm in}]_{M},\qquad M+6-|\beta|-|\nu|\ge7>6.
\]
For $|u|\ge1$, integrating $L$ times in the velocity direction
$\theta\cdot\nabla_v$ gives the exact identity
\[
\begin{split}
 (\ri u)^L\partial_\xi^\beta\partial_\eta^\nu
                    \widehat h_{\rm in}(k\theta,u\theta)=\iint_{\R^6}\re^{-\ri(k\theta\cdot x+u\theta\cdot v)}
       (\theta\cdot\nabla_v)^L
       \bigl[(-\ri x)^\beta(-\ri v)^\nu h_{\rm in}(x,v)\bigr]
       \,\rd x\,\rd v.
\end{split}
\]
The Leibniz rule uses at most $L\le11\le M$ physical derivatives
and only lowers the polynomial degree.  The same integrable majorant
therefore proves \eqref{eq:fourier-source-decay-appendix}.

For \eqref{eq:fourier-ray-L2-v6}, write $v=w\theta+z$, $z\in\theta^\perp$,
and define
\[
\begin{split}
 G_L(w)=\int_{\theta^\perp}\int_{\R^3}
   &\re^{-\ri k\theta\cdot x}(\theta\cdot\nabla_v)^L
       \bigl[(-\ri x)^\beta(-\ri v)^\nu h_{\rm in}(x,v)\bigr]
       \big|_{v=w\theta+z}\,\rd x\,\rd z.
\end{split}
\]
Its one-dimensional Fourier transform is exactly
\[
 \widehat G_L(u)=(\ri u)^L\partial_\xi^\beta\partial_\eta^\nu
                 \widehat h_{\rm in}(k\theta,u\theta).
\]
The five-dimensional marginal of the joint weight gives
\[
\begin{split}
 |G_L(w)|
 &\le C[h_{\rm in}]_{M}
       \int_{\R^3\times\theta^\perp}
       (1+w^2+|x|^2+|z|^2)^{(|\beta|+|\nu|-M-6)/2}\,\rd x\,\rd z\\
 &\le C[h_{\rm in}]_{M}\ang w^{|\beta|+|\nu|-M-1}
  \le C[h_{\rm in}]_{M}\ang w^{-2}.
\end{split}
\]
Indeed, the scaling $(x,z)=\ang w(x',z')$ leaves a finite integral
because $M+6-|\beta|-|\nu|>5$.  Rotations preserve the joint weight, so the
constant is uniform in $\theta$.  This bounds $\|G_L\|_{L^2(\R_w)}$;
the displayed Fourier identity and one-dimensional Plancherel prove
\eqref{eq:fourier-ray-L2-v6}.  The moment degree $|\beta|+|\nu|\le10$ and the physical
derivative order $L\le11$ enter these estimates independently.

For $0\le U\le1$, \eqref{eq:fourier-source-decay-appendix} with $L=0$
gives \eqref{eq:fourier-ray-L1-v6}.  For $U>1$, the part on $[0,1]$
is bounded in the same way.  Choose the integer $L=\lceil d\rceil+1\le11$;
then Cauchy--Schwarz gives
\[
\begin{split}
 \int_1^U u^d\left|\partial_\xi^\beta\partial_\eta^\nu
                    \widehat h_{\rm in}(k\theta,u\theta)\right|\,\rd u\le\left(\int_1^\infty u^{-2}\,\rd u\right)^{1/2}
       \left\|u^L\partial_\xi^\beta\partial_\eta^\nu
                    \widehat h_{\rm in}(k\theta,u\theta)\right\|_{L^2(\R_u)}
 \le C[h_{\rm in}]_{M},
\end{split}
\]
by \eqref{eq:fourier-ray-L2-v6}.  For the integer powers used below,
this is $L=d+1$.
\end{proof}

We now estimate the source history and the derivatives of its moving
upper limit.

\begin{lemma}
\label{lem:truncated-real-source-appendix}
For $k\in[0,k_0]$, $t\ge0$, and $s$ in a fixed compact interval, put
\[
 \mathfrak b(t;k,s,\theta):=\int_0^{kt}\re^{-\ri su}
                 \widehat h_{\rm in}(k\theta,u\theta)\,\rd u.
\]
For $M\ge11$ and $m+j+|\alpha|\le10$,
\begin{equation}\label{eq:real-source-finite-bound-appendix}
 |(k\partial_k)^m\partial_s^j\nabla_\theta^\alpha
                  \mathfrak b(t;k,s,\theta)|
 \le C\min\{kt,1\}[h_{\rm in}]_{M}.
\end{equation}
\end{lemma}
\begin{proof}
All radial derivatives are taken with $t$ fixed.  The first one gives
\begin{equation}\label{eq:moving-source-derivative-v6}
\begin{split}
 k\partial_k\mathfrak b(t;k,s,\theta)
 ={}\int_0^{kt}\re^{-\ri su}k\theta\cdot\nabla_\xi
              \widehat h_{\rm in}(k\theta,u\theta)\,\rd u
 +kt\,\re^{-\ri skt}\widehat h_{\rm in}(k\theta,kt\theta).
\end{split}
\end{equation}
Each $s$ derivative adds a factor $-\ri u$.  For a unit tangent
vector $e\in T_\theta\mathbb S^2$, angular differentiation gives
\[
 \nabla_{\theta,e}\widehat h_{\rm in}(k\theta,u\theta)
 =k\,e\cdot\nabla_\xi\widehat h_{\rm in}(k\theta,u\theta)
  +u\,e\cdot\nabla_\eta\widehat h_{\rm in}(k\theta,u\theta).
\]
Repeated angular derivatives have uniformly smooth coefficients in
finitely many sphere charts; derivatives of those coefficients give
lower-order terms.  The resulting terms split into two types.

\emph{Integral terms.}
These are bounded-coefficient multiples of
\[
 k^\ell\int_0^{kt}\re^{-\ri su}u^d
 \partial_\xi^\beta\partial_\eta^\nu
         \widehat h_{\rm in}(k\theta,u\theta)\,\rd u,
 \qquad 0\le d\le j+|\alpha|,\quad
 |\beta|+|\nu|\le m+|\alpha|.
\]
Here $\ell\ge0$ and both derivative bounds are at most ten.
Since $k\le k_0$, \eqref{eq:fourier-ray-L1-v6} bounds each term
by $C\min\{kt,1\}[h_{\rm in}]_{M}$.

\emph{Boundary terms.}
When a radial derivative reaches the upper limit, it supplies one
factor $kt$.  Further derivatives act on this factor, the exponential,
or the Fourier transform.  Every boundary term has the form
\[
 \re^{-\ri skt}s^{j_0} k^\ell(kt)^b
 \partial_\xi^\beta\partial_\eta^\nu
       \widehat h_{\rm in}(k\theta,kt\theta),
\]
up to a uniformly bounded coefficient, where $j_0,\ell\ge0$ and
\[
 1\le b\le m+j+|\alpha|\le10,
 \qquad |\beta|+|\nu|\le m+|\alpha|-1\le9.
\]
Since $s$ stays in a fixed compact interval and $k\le k_0$,
\eqref{eq:fourier-source-decay-appendix} with $L=b+1\le11$ gives
\[
 \left|\re^{-\ri skt}s^{j_0}k^\ell(kt)^b
     \partial_\xi^\beta\partial_\eta^\nu
       \widehat h_{\rm in}(k\theta,kt\theta)\right|
 \le C[h_{\rm in}]_{M}(kt)^b\ang{kt}^{-b-1}
 \le C\min\{kt,1\}[h_{\rm in}]_{M}.
\]
Combining the two types proves \eqref{eq:real-source-finite-bound-appendix}.
At $k=0$ or $t=0$ the source and the displayed differentiated terms
have value zero by continuous extension for fixed $t$.
\end{proof}

\subsection*{\centering Acknowledgment}
The research of Quoc-Hung Nguyen was supported by the CAS Project
for Young Scientists in Basic Research (Grant No.~YSBR-031) and by
the National Natural Science Foundation of China
(Grant Nos.~1251101538 and 12595282). Zhilong Xue is supported by the Postdoctoral Fellowship Program of CPSF
(Grant No.~GZB20260747) and the China Postdoctoral Science Foundation
(Grant No.~2026M793420).

\subsection*{\centering AI-use disclosure}

The authors used OpenAI's ChatGPT (including GPT-5.4, GPT-5.5,
GPT-5.6, and GPT-6) for language editing, organization,
bibliographic checks, mathematical discussions, and checks of
intermediate calculations.  They reviewed and revised the AI-assisted
material, independently verified the mathematical arguments and
calculations, and take full responsibility for the final manuscript.

\subsection*{Conflict of interest}

The authors declare that they have no conflict of interest.

\subsection*{Data availability}
No data were used or generated in this study.

\par\medskip
\begingroup\small
\noindent\textsc{Haina Li}\par\nobreak
\noindent School of Mathematics and Statistics,
Beijing Institute of Technology, Beijing 100081, China.\par\nobreak
\noindent\textit{Email address:}
\href{mailto:lihaina@bit.edu.cn}{\texttt{lihaina@bit.edu.cn}}.\par

\medskip
\noindent\textsc{Quoc-Hung Nguyen}\par\nobreak
\noindent State Key Laboratory of Mathematical Sciences,
Academy of Mathematics and Systems Science, Chinese Academy of Sciences,
Beijing 100190, China; Institute of Mathematics, Academy of Mathematics
and Systems Science, Chinese Academy of Sciences,
Beijing 100190, China.\par\nobreak
\noindent\textit{Email address:}
\href{mailto:qhnguyen@amss.ac.cn}{\texttt{qhnguyen@amss.ac.cn}}.\par

\medskip
\noindent\textsc{Zhilong Xue}\par\nobreak
\noindent Institute of Mathematics, Academy of Mathematics and Systems Science,
Chinese Academy of Sciences, Beijing 100190, China.\par\nobreak
\noindent\textit{Email address:}
\href{mailto:zhilongxue@amss.ac.cn}{\texttt{zhilongxue@amss.ac.cn}}.\par
\endgroup

\end{document}